\RequirePackage{ifpdf}
\ifpdf 
\documentclass[pdftex]{sigma}
\else
\documentclass{sigma}
\fi

\numberwithin{equation}{section}
\numberwithin{theorem}{section}

\let\leq\leqslant
\let\ge\geqslant
\let\geq\geqslant

\def\vvn#1>{\vadjust{\vsk#1>}}

\def\XXX/{{\slshape XXX}}

\newcommand{\nc}{\newcommand}
\nc{\on}{\operatorname}
\nc{\ch}{\mbox{ch}}
\nc{\Z}{{\mathbb Z}}
\nc{\C}{{\mathbb C}}
\nc{\R}{{\mathbb R}}
\nc{\pone}{{\mathbb C}{\mathbb P}^1}
\nc{\p}{\partial_u}
\nc{\F}{{\mathcal F}}
\nc{\arr}{\rightarrow}
\nc{\larr}{\longrightarrow}
\nc{\al}{\alpha}
\nc{\ri}{\rangle}
\nc{\lef}{\langle}
\nc{\W}{{\mathcal W}}
\nc{\la}{\lambda}
\nc{\ep}{\epsilon}
\nc{\eps}{\varepsilon}

\nc{\Om}{\Omega}
\nc{\su}{\widehat{{\mathfrak sl}}_2}
\nc{\sw}{{\mathfrak s}{\mathfrak l}}

\nc{\g}{{\mathfrak g}}
\nc{\h}{{\mathfrak h}}
\nc{\n}{{\mathfrak n}}
\nc{\N}{\widehat{\n}}
\nc{\G}{\widehat{\g}}
\nc{\De}{\Delta_+}
\nc{\gt}{\widetilde{\g}}
\nc{\Ga}{\Gamma}
\nc{\one}{{\mathbf 1}}
\nc{\z}{{\mathfrak Z}}
\nc{\zz}{{\mathcal Z}}
\nc{\Hh}{{\mathcal H}_\beta}
\nc{\qp}{q^{\frac{k}{2}}}
\nc{\qm}{q^{-\frac{k}{2}}}
\nc{\La}{\Lambda}
\nc{\wt}{\widetilde}
\nc{\qn}{\frac{[m]_q^2}{[2m]_q}}
\nc{\cri}{_{\on{cr}}}
\nc{\kk}{h^\vee}
\nc{\sun}{\widehat{\sw}_N}
\nc{\hh}{\widehat{\mathfrak h}}
\nc{\HH}{{\mathcal H}_{q,t}}
\nc{\ca}{\wt{{\mathcal A}}_{h,k}(\sw_2)}

\nc{\gl}{\widehat{{\mathfrak g}{\mathfrak l}}_2}
\nc{\el}{\ell}
\nc{\bi}{\bibitem}
\nc{\om}{\omega}
\nc{\WW}{\W_\beta}
\nc{\scr}{{\mathbf S}}
\nc{\ab}{{\mathbf a}}
\nc{\rr}{r}
\nc{\ol}{\overline}
\nc{\con}{qt^{-1} + q^{-1}t}
\nc{\den}{q^{\el-1} t^{-\el+1}+ q^{-\el+1} t^{\el-1}}
\nc{\ds}{\displaystyle}
\nc{\B}{B}
\nc{\A}{{\mathbb A}}
\nc{\GG}{{\mathcal G}}
\nc{\UU}{{\mathcal U}}
\nc{\MM}{{\mathcal M}}
\nc{\CC}{{\mathcal C}}
\nc{\GL}{{}^L G}
\nc{\dzz}{\frac{dz}{z}}
\nc{\Res}{\on{Res}}
\nc{\rep}{{\mathcal R}ep \;}
\nc{\uqg}{U_q \G}
\nc{\uqgg}{U_q \g}
\nc{\Fq}{{\mathbb F}_q}

\nc{\stimes}{\ltimes}
\nc{\K}{\hat{\mathcal K}}
\nc{\Ql}{\ol{\mathbb Q}_\ell}

\nc{\ga}{\gamma}
\nc{\PL}{{}^L P}
\nc{\E}{\mc E}
\nc{\mc}{\mathcal}
\nc{\mbf}{\mathbf}
\nc{\bb}{{\mathfrak b}}
\nc{\OO}{{\mc O}}
\nc{\Po}{{\mc P}}
\nc{\V}{{ V_{\bs\La}}}
\nc{\yy}{{\mc Y}}
\nc{\M}{\mathcal M}
\nc{\Coh}{{{\mathcal C}oh}}
\nc{\Cohn}{\Coh_n}
\nc{\f}{{\mathcal F}}
\nc{\si}{_E}
\nc{\Gaf}{{\mathbb G}_{a,\Fq}}
\nc{\KK}{{\mathfrak k}}

\nc{\PO}{{\mathbb P^1}}
\nc{\PR}{{\mathbb P^r}}

\nc{\Wr}{{ {\rm Wr}}}

\nc{\bs}{\boldsymbol}

\nc{\U}{\mathcal{U}}

\nc{\sing}{{\rm Sing}\,}

\nc{\Ll}{\La - \al(\bs l)}

\nc{\nL}{L_{\om_r}^{\otimes n}[\mu]}
\nc{\nnL}{L_{\om_r}^{\otimes n}}

\nc{\snL}{ \sing L_{\om_r}^{\otimes n}[\mu]}
\nc{\btz}{ \om(\bs t; \bs z)}
\nc{\SL}{ \mathfrak{sl}}
\nc{\GR}{ {G(r+1,d)}}

\def\glN{\mathfrak{gl}_N}

\let\der\partial

\let\si\sigma

\nc{\Te}{\mathfrak{T}}

\nc{\End}{ {\rm End}}
\nc{\Wrd}{{\rm Wr}^{\rm(d)}}

\nc{\Y}{Y(\glN)}
\nc{\MLz}{M_{\bs \La}(\bs z)}

\nc{\D}{\mathfrak {D}}

\nc{\DM}{\D_{Q, M(\bs z)}(u,\tau)}

\nc{\Dt}{\D_{\tilde {\bs t}}(u,\tau)}
\nc{\Se}{\mathfrak S}
\nc{\Xe}{\mathcal X}
\nc{\Tee}{\mathcal S}
\nc{\glx}{\glN[x]}
\nc{\dd}{\mathfrac d{du}}

\nc{\Kk}{{\bs K}}
\nc{\Sk}{{S_{\bs k}}}
\newcommand{\AT}{ A_{p,\Phi}}
\newcommand{\nash}{{\V[\La_\infty]}}
\newcommand{\snash}{{\Sing \V[\La_\infty]}}
\newcommand{\Sing}{{\rm {Sing\,}}}
\newcommand{\Hess}{{\rm {Hess}^{(a)}}}
\newcommand{\HessF}{{\rm Hess}\, F}
\newcommand{\OS}{\mathcal {A}}
\def\FF{{\mathcal F}}
\def\a{{\alpha}}

\nc{\glr}{{\frak{gl}_{r+1}}}
\nc{\slr}{{\frak{sl}_{r+1}}}
\nc{\glrs}{{\frak{gl}_{r+1}[s]}}
\nc{\lc}{{,\dots ,}}
\nc{\dl}{{\delta}}
\nc{\rdet}{{{\rm rdet}}}
\newtheorem{thm}{Theorem}[section]
\newtheorem{cor}[thm]{Corollary}
\newtheorem{Lemma}[thm]{Lemma}

\newtheorem{conj}[thm]{Conjecture}
\newtheorem{assumption}[thm]{Assumption}
\newtheorem{data}[thm]{Data}
 
\newtheorem{rem}[thm]{Remark}

\newcommand{\Pee}{{\mathbb P}}

\newcommand{\T}{{\mathcal{T}}}
\nc{\s}{{\sing}}

\begin{document}

\allowdisplaybreaks

\renewcommand{\thefootnote}{$\star$}

\renewcommand{\PaperNumber}{032}

\FirstPageHeading

\ShortArticleName{Quantum Integrable Model of an Arrangement of Hyperplanes}

\ArticleName{Quantum Integrable Model\\ of an Arrangement  of Hyperplanes\footnote{This paper is a
contribution to the Special Issue ``Relationship of Orthogonal Polynomials and Special Functions with Quantum Groups and Integrable Systems''. The
full collection is available at
\href{http://www.emis.de/journals/SIGMA/OPSF.html}{http://www.emis.de/journals/SIGMA/OPSF.html}}}

\Author{Alexander~VARCHENKO}

\AuthorNameForHeading{A.~Varchenko}

\Address{Department of Mathematics, University of North Carolina
at Chapel Hill,\\
Chapel Hill, NC 27599-3250, USA}
\Email{\href{mailto:anv@email.unc.edu}{anv@email.unc.edu}}

\ArticleDates{Received July 19, 2010, in f\/inal form March 19, 2011;  Published online March 28, 2011}

\Abstract{The goal of this paper is to give a geometric construction of the Bethe algebra (of Hamiltonians) of a Gaudin model associated to
a simple Lie algebra.
More precisely, in this paper a quantum integrable model is assigned to a weighted arrangement of af\/f\/ine hyperplanes.
We show (under certain assumptions) that the algebra of Hamiltonians of the model is isomorphic to
the algebra of functions on the critical set of the corresponding master function. For a discriminantal
arrangement we show (under certain assumptions) that the symmetric part of the
algebra of  Hamiltonians is isomorphic to the Bethe algebra of the corresponding Gaudin model.
It is expected that this correspondence holds in general (without the assumptions).
As a byproduct of constructions we show that in a Gaudin model (associated to an arbitrary simple Lie algebra), the Bethe vector,
corresponding to an isolated critical point of the master function, is nonzero.}

\Keywords{Gaudin model; arrangement of hyperplanes}

\Classification{82B23; 32S22; 17B81; 81R12}

\begin{flushright}
\it To Sabir Gusein-Zade on the occasion of his 60-th birthday
\end{flushright}

{\small

\tableofcontents

}

\renewcommand{\thefootnote}{\arabic{footnote}}
\setcounter{footnote}{0}

\section{Introduction}

\subsection{Quantum integrable models and Bethe ansatz}

A quantum integrable model is a vector space $V$ and an
``interesting'' collection of commuting linear operators
$K_1,K_2,
\ldots :\ V\to V$. The operators are called
Hamiltonians or transfer matrices or conservation laws. The problem is
to f\/ind common eigenvectors and eigenvalues.

The Bethe ansatz  is a method to diagonalize commuting linear operators.
One invents a~vector-valued function $v(t)$ of some new parameters
$t=(t_1,\dots,t_k)$ and f\/ixes the parameters  so that $v(t)$ becomes a common eigenvector
of the Hamiltonians. One shows that $v(t)\in V$ is an eigenvector
if $t$ satisf\/ies some system of equations,
\begin{gather}
\label{BAE intr}
\Psi_j(t)=0,\qquad j = 1, \dots, k.
\end{gather}
The equations are called the Bethe ansatz equations. The vector $v(t)$ corresponding to a solution of
the equations is called a Bethe vector. The method is called the Bethe ansatz
method.

\subsection{Gaudin model}

One of the simplest and interesting  models is the quantum Gaudin
model introduced in~\cite{G1} and~\cite{G2}.
Choose a simple Lie algebra~$\g$,  an orthonormal basis $\{J_a\}$ of
$\g$ with respect to a nondegenerate $\g$-invariant bilinear form, and a collection of
 distinct complex numbers $x=(x_1,\dots,x_N)$. Then one def\/ines certain elements
of the $N$-th tensor power of the universal enveloping algebra of~$\g$, denoted by
$\Kk_1(x),\dots, \Kk_N(x) \in (U\g)^{\otimes N}$ and
called the Gaudin Hamiltonians,
\[
\Kk_b(x) = \sum_{b\ne c}\sum_a \frac {J_a^{(b)}J_a^{(c)}}{x_b-x_c},
\qquad
b = 1,\dots,N .
\]
Here we use the standard notation: if $J\in\g$, then $J^{(i)}
 = 1^{\otimes (i-1)}\otimes J \otimes 1^{\otimes (N-i)}$.

The Gaudin Hamiltonians commute with
each other and commute with the diagonal subalgebra
$(U\g)_{\rm diag} \subset (U\g)^{\otimes N}$,
\[
 [\Kk_b(x),\Kk_c(x)]=0, \qquad
[\Kk_b(x), (U\g)_{\rm diag}] = 0 .
\]
Let $V_\La$ denote the f\/inite-dimensional irreducible $\g$-module with highest
weight $\La$. Decompose a~tensor product $V_{\bs\La} =
\otimes_{b=1}^N V_{\La_b}$ into irreducible
$\g$-modules,
\begin{gather}
\label{mult}
V_{\bs\La} = \oplus_{ \La_\infty} V_{\La_\infty} \otimes  W_{\bs \La,\La_\infty} ,
\end{gather}
where $W_{\bs \La,\La_\infty}$ is the multiplicity space of a representation $V_{\La_\infty}$.
 The Gaudin Hamiltonians act on $V_{\bs\La}$, preserve decomposition \eqref{mult}, and by
Schur's lemma induce commuting linear operators on each multiplicity space $W=W_{\bs \La,\La_\infty}$,
\[
\Kk_1(x),\dots,\Kk_N(x) : \  W \to W .
\]
These commuting linear operators on a multiplicity space constitute
the quantum Gaudin model.

Thus, the  Gaudin model depends on $\g$, $(\;,\;)$,  highest weights
 $\La_1,\dots,\La_N, \La_\infty$ and complex numbers
$x_1,\dots,x_N$.

\subsection{Gaudin model as a semiclassical limit of KZ equations}
\label{Gaudin model and KZ equations}

On every multiplicity space $W=W_{\bs \La,\La_\infty}$
of the tensor product $V_{\bs\La} = \oplus_{\La_\infty} V_{\La_\infty} \otimes  W_{\bs \La,\La_\infty}$
  one has a system of Knizhnik--Zamolodchikov (KZ) dif\/ferential equations,
\[
\kappa \frac{\der I}{\der x_b}(z) = \Kk_b(x)I(x),
\qquad b=1,\dots,N,
\]
where $x=(x_1,\dots,x_N )$,   $I(x)\in W$ is the unknown function,
$\Kk_b(x)$ are the Gaudin Hamil\-to\-nians and $\kappa\in\C^\times$ is a parameter of the dif\/ferential equations.
KZ equations are equations for conformal blocks in the Wess--Zumino--Novikov--Witten conformal f\/ield theory.

For any value of  $\kappa$,  KZ equations def\/ine a f\/lat connection on the trivial bundle $\C^N\times W\to \C^N$
with singularities over the diagonals.
 The f\/latness conditions,
\[
\left[\kappa\frac{\der }{\der x_b} - \Kk_b(x), \kappa\frac{\der }{\der x_c} - \Kk_c(x)\right] = 0 ,
\]
in particular, imply  the conditions $\left[\Kk_b(x),\Kk_c(x)\right]=0$, which
are the commutativity conditions for the Gaudin Hamiltonians.

Thus, we observe two related problems:
\begin{enumerate}\itemsep=0pt
\item[(1)] Given a nonzero number $\kappa$, f\/ind solutions of
the KZ equations.
\item[(2)] Given $x$, f\/ind eigenvectors of the Gaudin Hamiltonians.
\end{enumerate}

Problem (2) is a semiclassical limit of problem (1) as $\kappa$ tends to zero. Namely,
assume  that the KZ equations have an asymptotic solution
of the form
\[
I(x) = e^{P(x)/\kappa}\big(w_0(x) + \kappa w_1(x) + \kappa^2 w_2(x) + \cdots \big)
\]
as $\kappa \to 0$. Here $P(x)$ is a scalar function and $w_0(x),w_1(x),w_2(x),\dots$ are some
$W$-valued functions of $x$.
Substituting this expression to the KZ equations and
equating the leading terms we get
\[
\frac{\der P}{\der x_b}(x) w_0(x) = \Kk_b(x) w_0(x) .
\]
Hence, for any $x$ and $b$, the vector
 $w_0(x)$ is an eigenvector of the Gaudin Hamiltonian $\Kk_b(x)$ with eigenvalue
$\frac{\der P}{\der x_b}(x)$.

Thus, in order to diagonalize the Gaudin Hamiltonians it is enough to construct asymptotic solutions
to the KZ equations.

\subsection{Bethe ansatz in the Gaudin model}

The Gaudin model has a Bethe ansatz  \cite{B,BF, RV}.
 The Bethe ansatz has special features.
Namely,
\begin{enumerate}\itemsep=0pt
\item[$(i)$]
In the Gaudin model, there exists a function $\Phi(t_1,\dots,t_k)$ (the master function) such that
the Bethe ansatz equations~\eqref{BAE intr}
 are the critical point equations for the master function,
\[
\frac{\der \Phi}{\der t_j}(t_1,\dots,t_k)=0,
\qquad
j=1,\dots,k.
\]

\item[$(ii)$]
In the Gaudin model, the vector space $W$ has a symmetric bilinear form
$S$ (the tensor Shapovalov form) and Gaudin Hamiltonians are symmetric operators.

\item[$(iii)$]
In the Gaudin model, the Bethe
vectors assigned to (properly understood) distinct critical points are orthogonal  and  the square of the norm of a Bethe
vector equals the Hessian of  the master
function at the corresponding critical point,
\[
S(v(t),v(t)) = \det \left(
\frac{\der^2\Phi }{\der t_i\der t_j} \right)(t) .
\]
In particular, the Bethe vector corresponding to a nondegenerate critical point is nonzero.
\end{enumerate}

These statements indicate a connection between
the  Gaudin Hamiltonians  and
a  mysterious master function (which
 is not present in the def\/inition of the Gaudin model).

 One of the  goals of this paper is to uncover this mystery
 and show that the Bethe ansatz can be interpreted as an elementary construction in the
 theory of arrangements, where the master function is a basic object.

\subsection{Gaudin model and arrangements}
\label{Gaudin model and arrangements}

For a weighted  arrangement of af\/f\/ine hyperplanes, we will  construct (under certain assumptions)
a quantum integrable model, that is, a vector space $W$ with
a symmetric bilinear form $S$ (called the contravariant form),
a collection of commuting symmetric linear operators on $W$ (called {\it
the naive geometric Hamiltonians}),
a master function
$\Phi(t)$ and vectors $v(t)$ (called the special vectors or called the values of the canonical element) which are eigenvectors of the linear operators if $t$ is a critical point of the master function.

Then for a given Gaudin model  $(W, S,\Kk_1(x),\Kk_2(x),\ldots : W\to W)$, one can f\/ind a suitable (discriminantal) arrangement and
identify the objects of the Gaudin model
with the correspon\-ding objects of the quantum integrable model of the
arrangement. After this identif\/ication, the master function
$\Phi (t)$ and the special vectors $v(t)$ of the arrangement provide the Gaudin model with a Bethe ansatz,
that is, with a method to diagonalize the Gaudin Hamiltonians.

\subsection{Bethe algebra}
\label{sec Bethe algebra intr}

Let $W_{\bs \La,\La_\infty}$ be the vector space of a Gaudin model. It turns out that the subalgebra of
$\End(W_{\bs \La,\La_\infty})$ generated by the Gaudin Hamiltonians can be extended to a larger commutative
subalgebra called the Bethe algebra. A general construction of the Bethe algebra for a~simple Lie algebra $\g$
is given in~\cite{FFR}. That construction is formulated in terms of the center of the universal enveloping algebra
of the corresponding af\/f\/ine Lie algebra $\hat \g$  at the critical level. As a result of that construction,
for any $x$  one obtains a commutative subalgebra $\B(x)\subset (U\g)^{\otimes N}$ which
commutes with the diagonal subalgebra $U\g\subset (U\g)^{\otimes N}$.
To def\/ine the Bethe algebra of~$V_{\bs\La}$ or of~$W_{\bs \La,\La_\infty}$ one considers the
 image of~$\B(x^0)$ in~$\End (V_{\bs\La})$ or in~$\End(W_{\bs \La,\La_\infty})$. The Gaudin Hamiltonians $\Kk_b(x)$
 are elements of the Bethe algebra of $V_{\bs\La}$ or of $W_{\bs \La,\La_\infty}$.

The construction of the Bethe algebra in \cite{FFR} is not explicit and it is not easy to study the Bethe algebra of
$W_{\bs \La,\La_\infty}$ for a particular Gaudin model.
For example, a standard dif\/f\/icult question is if the Bethe algebra of
$W_{\bs \La,\La_\infty}$ is a maximal commutative subalgebra of $\End(W_{\bs \La,\La_\infty})$,
cf.~\mbox{\cite{MTV7,FFRy}}.

\subsection{Algebra of geometric Hamiltonians}
\label{sec aLgebra of geome hams}

In this paper we
address the following problem. Is there a geometric construction of the Bethe algebra?
For the quantum integrable model $(W, S,   K_1,K_2,\ldots : W\to W)$ of a given weighted arrangement,
can one def\/ine a geometric ``Bethe algebra'' $A$, which is a maximal commutative subalgebra of $\End(W)$,
which contains the naive geometric Hamiltonians $K_1, K_2, \dots$,
 and which can be identif\/ied with the Bethe algebra of~\cite{FFR} for discriminantal arrangements?

In this paper (under certain assumptions) we construct  an algebra $A$ called {\it
the algebra of geometric Hamiltonians} and identify it (in certain cases) with the Bethe algebra of the Gaudin model.

\subsection[Quantum integrable model of a weighted arrangement (or dynamical theory \ \\ of arrangements)]{Quantum integrable model of a weighted arrangement\\
(or dynamical theory of arrangements)}
\label{sec qim}

To def\/ine the quantum integrable model of an arrangement we consider in an af\/f\/ine space $\C^k$ (with coordinates
$t=(t_1,\dots,t_k)$) an arrangement of
$n$ hyperplanes,  $k<n$. Each hyperplane is allowed to move parallelly to itself. The parallel shift of the $i$-th
hyperplane is measured by a~number~$z_i$ and  for every
$z=(z_1,\dots,z_n)\in \C^n$ we get in $\C^k$ an af\/f\/ine arrangement of hyperplanes $\A(z) = (H_j(z))$,
\[
H_j(z) = \big\{t\in \C^k\ |\ g_j(t) +z_j=0\big\},
\]
where $g_j(t)$ are given linear functions on $\C^k$.
We assign nonzero numbers $a=(a_j)$ to the hyperplanes of  $\A(z)$ (the numbers do not depend on $z$) and obtain
a family of parallelly translated weighted hyperplanes.

For  generic $z\in\C^n$, the arrangement $\A(z)$ has normal crossings only. The discriminant $\Delta\subset\C^n$
is the subset of all points $z$ such that $\A(z)$ is not with normal crossings only.

The master function $\Phi$ on $\C^n\times\C^k$ is the function
\[
\Phi(z,t) = \sum_j a_j \log (g_j(t)+z_j) .
\]
Let $\OS(\A(z)) = \oplus_{p=0}^k
\OS^p(\A(z))$ be the Orlik--Solomon algebra and
$\FF(\A(z)) = \oplus_{p=0}^k
\FF^p(\A(z))$ the dual vector space. We are interested in the top degree components $\OS^k(\A(z))$ and $\FF^k(\A(z))$.

The weights $a$ def\/ine on  $\FF^k(\A(z))$ a symmetric
bilinear form $S^{(a)}$ (called the contravariant form) and a degree-one element $\nu(a) = \sum_j a_j (H_j(z)) \in \OS^1(\A(z))$.
Denote $\sing \FF^k(\A(z)) \subset \FF^k(\A(z))$ the annihilator of the subspace
$\nu(a)\cdot\OS^{k-1}(\A(z))\subset \OS^k(\A(z))$.

For  $z^1,z^2\in \C^n-\Delta$, all combinatorial objects of the arrangements $\A(z^1)$ and $\A(z^2)$ can be
canonically identif\/ied.
In particular, the spaces $\FF^k(\A(z^1))$, $\FF^k(\A(z^2))$ as well as the spaces $\sing \FF^k(\A(z^1))$,
$\sing\FF^k(\A(z^2))$ can be canonically identif\/ied.
For  $z\in\C^n-\Delta$, we denote
$V=\FF^k(\A(z)), \sing V=\sing \FF^k(\A(z))$.

For any nonzero number $\kappa$, the hypergeometric integrals
\begin{gather*}
\int_{\gamma(z)} e^{\Phi(z,t)/\kappa} \omega ,
\qquad
\omega \in \OS^k(\A(z)) ,
\end{gather*}
def\/ine a Gauss--Manin (f\/lat) connection on the trivial bundle $\C^n\times \sing V\to \C^n$ with singularities over the discriminant.
The Gauss--Manin dif\/ferential equations for horizontal sections of
the connection on $\C^n\times \Sing V \to \C^n$ have the form
\begin{gather}
\label{dif eqn intr}
\kappa \frac{\der I}{\der z_j}(z) = K_j(z)I(z),
\qquad
j\in J,
\end{gather}
where $I(z)\in \Sing V$ is a horizontal section,
 $K_j(z): V\to V$, $j\in J,$ are suitable linear operators  preserving $\Sing V$ and
independent of~$\kappa$. For every $j$,
the operator $K_j(z)$ is a rational function of $z$ regular on $\C^n-\Delta$.
Each operator is symmetric
with respect to the contravariant form $S^{(a)}$.

These dif\/ferential equations are our source of quantum integrable models for weighted arrangements.
The quantum integrable models are the semiclassical limit of these
dif\/ferential equations similarly to the transition from  KZ equations to the Gaudin model in
Section~\ref{sec Bethe algebra intr}.

The f\/latness of the connection for all $\kappa$ implies
the commutativity of the operators,
\[
K_i(z)|_{\Sing V}K_j(z)|_{\Sing V}=K_j(z)|_{\Sing V}K_i(z)|_{\Sing V}
\qquad
\text{ for all} \ i,j\ \text{and}\   z\in \C^n-\Delta .
\]
Let $V^*$ be the space dual to $V$. If $M:V\to V$ is a linear operator, then
 $M^*:V^*\to V^*$ denotes the dual operator.
Let
$W \subset V^*$ be the image of $V$ under the map
$ V \to  V^*$ associated with the contravariant form and
$\sing W\subset W$ the image of $\sing V\subset V$.
The contravariant form induces on $W$ a nondegenerate
symmetric bilinear form, also denoted by $S^{(a)}$.

The operators $K_j(z)^*$ preserve the subspaces $\sing W\subset W\subset V^*$. The operators
$K_j(z)^*|_{W} :W\to W$ are symmetric with respect to the contravariant form. The operators
$K_j(z)^*|_{\sing W} :\sing W\to \sing W$, $j\in J$, commute.

For   $z\in \C^n-\Delta$, we def\/ine {\it the quantum integrable model assigned to}
$(\A(z), a)$ to be the collection
\[
\big(\sing W;\ S^{(a)}|_{\sing W};\ K_1(z)^*|_{\sing W},\dots, K_n(z)^*|_{\sing W} : \ \sing W\to \sing W\big).
\notag
\]
The unital subalgebra of $\End(\sing W)$ generated by operators
$K_1(z)^*|_{\sing W}$, \dots,
$K_n(z)^*|_{\sing W}$ will be called {\it the
algebra of geometric Hamiltonians of $(\A(z),a)$}.

It is clear that any weighted
 arrangement with normal crossings only can be realized as a~f\/iber $(\A(z), a)$ of such
a construction and, thus,  a weighted   arrangement with normal crossings only
is provided with a quantum integrable model.

If $z^0$ is a point of the discriminant,  the construction of the quantum integrable model assigned to the
arrangement $(\A(z^0),a)$ is more delicate. The operator valued functions $K_j(z)$ may have f\/irst order
poles at $z^0$. We write $K_j(z) = K_j^0(z) + K^1_j(z)$, where $K_j^0(z)$ is the polar part at~$z^0$ and~$K_j^1(z)$ the regular part. By suitably regularizing the operators $K_j^1(z^0)$, we make them
(under certain assumptions) preserve a suitable subspace
of~$\sing V$,  commute on that subspace and be symmetric with respect to the contravariant form.
The algebra of regularized operators~$K_j^1(z^0)$
on that subspace produces
the quantum integrable model assigned to the f\/iber~$(\A(z^0),a)$.
It may happen that some linear
combinations (with constant coef\/f\/icients) $K_\xi(z)=\sum_j\xi_jK_j(z)$ are regular at~$z^0$.
In this case  the operator $K_\xi(z^0)$  preserves that subspace and is an element of the
algebra of Hamiltonians. In this case the operator~$K_\xi(z^0)$  is called a naive geometric Hamiltonian.
(It is naive in the sense that we don't need to go through the regularization procedure to produce that element of
the algebra of Hamiltonians. In that sense, for $z\in \C^n-\Delta$ all elements of the algebra of geometric
Hamiltonians of the arrangement $(\A(z),a)$ are naive.)

For applications to the Gaudin model one needs a suitable equivariant version of the described construction, the corresponding family of parallelly
translated hyperplanes has a symmetry group and the Gaudin model
is identif\/ied with  the skew-symmetric part of the corresponding quantum integrable model of the arrangement.

\subsection{Bethe ansatz for the quantum integrable model of an arrangement}

 The Hamiltonians of the model are
(suitably regularized) right
hand sides of the Gauss--Manin dif\/ferential equations~\eqref{dif eqn intr}. The solutions to the equations are the
integrals  $\int_{\gamma(z)} e^{\Phi(z,t)/\kappa} \omega$. By taking the semiclassical limit of the integrals
 as $\kappa$ tends to zero,
we obtain  eigenvectors of
the Hamiltonians, cf.\ Section~\ref{Gaudin model and KZ equations}.
 The eigenvectors of the Hamiltonians are labeled
by the critical points of the phase $\Phi$ of the integrals due to the steepest descent method.

\subsection{Geometric interpretation of the algebra of  Hamiltonians}
\label{sec intr geom interpret}

It turns out that the solutions  $\int_{\gamma(z)} e^{\Phi(z,t)/\kappa} \omega$ to equations~\eqref{dif eqn intr}
produce more than just eigenvectors of the geometric Hamiltonians. They also produce a geometric interpretation
of the whole algebra of geometric Hamiltonians. It turns out that the algebra of geometric Hamiltonians is naturally isomorphic
 (under certain conditions) to the algebra of functions on the critical set of the master function~$\Phi$. The isomorphism is established  through the semiclassical limit of the integrals.
Moreover, this isomorphism identif\/ies the residue bilinear form on the algebra of functions and the contravariant form of the arrangement.

This geometric interpretation of the algebra of geometric Hamiltonians
 is motivated by the recent paper~\cite{MTV5} where a connection between the algebra of functions on the critical set of the
master function and the Bethe algebra of the $\glr$ Gaudin model is established, cf.~\cite{MTV4}.

\subsection{Byproducts of constructions}
\label{Byproducts of the constructions}

The general motive of this paper is the interplay between the combinatorially def\/ined linear objects of a weighted
arrangement  and the critical set of the corresponding master function (which is a nonlinear characteristics of
the arrangement). As byproducts of our considerations we get relations between linear and nonlinear characteristics of an arrangement.
For example, we prove that the sum of Milnor numbers of the critical points of the master function is not greater than the rank of the
contravariant form on $\sing V$.

As another example of such an interaction we show that in any Gaudin model (associated with  any simple Lie algebra) the Bethe vector
corresponding to an isolated critical points of the master function is nonzero. That result is known for nondegenerate critical points,
see~\cite{V4}, and for the Gaudin models associated with the Lie algebra $\glr$,
see~\cite{MTV5}.

\subsection{Exposition of the material}
\label{Exposition of the material}

In Section~\ref{Sec Arrangements}, basic facts of the theory of arrangements are collected.
The main objects are the space of f\/lags,  contravariant form,  master function,  canonical element.

In Section~\ref{sec par trans}, a family of parallelly translated hyperplanes is introduced.
In Section~\ref{sec Preservation}, remarks on the conservation of the number of
critical points of the master function under deformations are presented.

In Section~\ref{Hamiltonians of good fibers},
the Gauss--Manin dif\/ferential equations are considered.
 The  quantum integrable model of a weighted arrangement with normal crossings only
is introduced.  The ``key identity''~\eqref{formula H} is formulated.

In Section~\ref{sec asympt}, the asymptotic solutions to the Gauss--Manin dif\/ferential equations are discussed.
In Section~\ref{sec real posit bad}, the quantum integrable model of any f\/iber
$(\A(z^0),a)$, $z^0\in\Delta$, is def\/ined under assumptions of positivity of weights $(a_j)$ and reality
of functions $(g_j(t))$. A general conjecture is formulated.

In Section~\ref{Geometric interpretation of the algebra of Hamiltonians},
it is shown that the algebra of geometric Hamiltonians is isomorphic to the algebra of functions on the critical
set of the master function. That fact is proved for any $(\A(z),a)$, $z\in \C^n$
under Assumption \ref{i} of positivity of  $(a_j)$ and reality
of  $(g_j(t))$.
In Section~\ref{More on Hamiltonians of bad fibers}, more results in this direction are obtained, see Theorems~\ref{thm BAD} and~\ref{thm ham of unb arr}.

In Section~\ref{Arrangements with symmetries sec}, an equivariant version of the algebra of
geometric Hamiltonians is introduced and in Section~\ref{sec Gaudin} relations with the Gaudin model are described.

\newpage

\section{Arrangements}
\label{Sec Arrangements}
\subsection{An af\/f\/ine arrangement}
\label{An affine arrangement}
Let $k$ and $n$ be positive integers, $k<n$. Denote $J=\{1,\dots,n\}$.

Let $\A =(H_j)_{j\in J}$,  be an arrangement of $n$ af\/f\/ine hyperplanes in
$\C^k$. Denote
\[
U = \C^k - \cup_{j\in J} H_j   ,
\]
the complement.
An edge $X_\al \subset \C^k$ of the arrangement $\A$ is a nonempty intersection of some
hyperplanes  of $\A$. Denote by
 $J_\al \subset J$ the subset of indices of all hyperplanes containing $X_\al$.
Denote
$l_\a = \mathrm{codim}_{\C^k} X_\a$.

We always
assume that $\A$ is essential, that is, $\A$ has a vertex, an edge which is a point.

An edge is called dense if the subarrangement of all hyperplanes containing
it is irreducible: the hyperplanes cannot be partitioned into nonempty
sets so that, after a change of coordinates, hyperplanes in dif\/ferent
sets are in dif\/ferent coordinates. In particular, each hyperplane of
$\A$ is a dense edge.

\subsection{Orlik--Solomon algebra}

Def\/ine complex vector spaces $\OS^p(\A)$, $p = 0,  \dots, k$.
 For $p=0$ set $\OS^p(\A)=\C$. For  $p \geq 1$,
 $\OS^p(\A)$   is generated by symbols
$(H_{j_1},\dots ,H_{j_p})$ with ${j_i}\in J$, such that
\begin{enumerate}\itemsep=0pt
\item[$(i)$] $(H_{j_1},\dots ,H_{j_p})=0$
if $H_{j_1}$,\dots ,$H_{j_p}$ are not in general position, that is, if the
intersection $H_{j_1}\cap \dots  \cap H_{j_p}$ is empty or
 has codimension
 less than $p$;
\item[$(ii)$]
$ (H_{j_{\sigma(1)}},\dots ,H_{j_{\sigma(p)}})=(-1)^{|\sigma|}
(H_{j_1},\dots ,H_{j_p})
$
for any element $\sigma$ of the
symmetric group $ S_p$;
\item[$(iii)$]
$\sum_{i=1}^{p+1}(-1)^i (H_{j_1},\dots ,\widehat{H}_{j_i},\dots ,H_{j_{p+1}}) = 0
$
for any $(p+1)$-tuple $H_{j_1},\dots ,H_{j_{p+1}}$ of hyperplanes
in $\A$ which are
not in general position and such that $H_{j_1}\cap\dots \cap H_{j_{p+1}}\not = \varnothing$.
\end{enumerate}

The direct sum $\OS(\A) = \oplus_{p=1}^{N}\OS^p(\A)$ is an
 algebra with respect to
 multiplication
\[
(H_{j_1},\dots ,H_{j_p})\cdot(H_{j_{p+1}},\dots ,H_{j_{p+q}}) =
 (H_{j_1},\dots ,H_{j_p},H_{j_{p+1}},\dots ,H_{j_{p+q}})  .
\]
 The algebra is  called the Orlik--Solomon algebra of  $\A$.

\subsection{Weights}
\label{sec weights}

An arrangement $\A$ is weighted if a map $a: J\to \C$, $j\mapsto a_j,$ is given;
 $a_j$ is called the weight of~$H_j$.
For an edge $X_\al$, def\/ine its weight as
$a_\al = \sum_{j\in J_\al}a_j$.

We always assume that $a_j\neq 0$ for every $j\in J$.

 Def\/ine
\[
\nu(a)  =  \sum_{j\in J} a_j (H_j)  \in   \OS^1(\A)  .
\]
 Multiplication by $\nu(a)$ def\/ines a dif\/ferential
\[
d^{(a)}\ :\   \OS^p(\A)\ \to\ \OS^{p+1}(\A) ,
\qquad
 x \ \mapsto\ \nu(a)\cdot x   ,
\]
on  $\OS(\A)$.

\subsection{Space of f\/lags, see \cite{SV}}

 For an edge $X_\alpha$,   $l_\alpha=p$, a f\/lag starting at $X_\alpha$ is a sequence
\[
X_{\alpha_0}\supset
X_{\alpha_1} \supset \dots \supset X_{\alpha_p} = X_\alpha
\]
of edges such that
$ l_{\alpha_j} = j$ for $j = 0, \dots , p$.

 For an edge $X_\alpha$,
 we def\/ine $\overline{\FF}_{\alpha}$  as  the
complex vector space  with basis vectors
\[
\overline{F}_{{\alpha_0},\dots,{\alpha_p}=\alpha}
\]
 la\-bel\-ed by the elements of
the set of all f\/lags  starting at~$X_\alpha$.

 Def\/ine  $\FF_{\alpha}$ as the quotient of
$\overline{\FF}_{\alpha}$ by the subspace generated by all
the vectors of the form
\begin{gather}
\label{flagrelation} \sum\limits_{X_\beta,\,
X_{\alpha_{j-1}}\supset X_\beta\supset X_{\alpha_{j+1}}}
\overline {F}_{{\alpha_0},\dots,
{\alpha_{j-1}},{\beta},{\alpha_{j+1}},\dots,{\alpha_p}=\alpha} .
\notag
\end{gather}
Such a vector is determined  by  $j \in \{ 1, \dots , p-1\}$ and
an incomplete f\/lag $X_{\alpha_0}\supset\dots \supset
X_{\alpha_{j-1}} \supset X_{\alpha_{j+1}}\supset\dots \supset
X_{\alpha_p} = X_\alpha$ with $l_{\alpha_i} = i$.

Denote by ${F}_{{\alpha_0},\dots,{\alpha_p}}$ the image in $\FF_\alpha$ of the basis vector
$\overline{F}_{{\alpha_0},\dots,{\alpha_p}}$.  For $p=0,\dots,k$, set
\[
{\FF}^p(\A)  =  \oplus_{X_\alpha,  l_\alpha=p}\ {\FF}_{\alpha}  .
\]

\subsection{Duality, see \cite{SV}}
The vector spaces $\OS^p(\A)$ and $\FF^p(\A)$ are dual.
The pairing $ \OS^p(\A)\otimes\FF^p(\A) \to \C$ is def\/ined as follows.
{}For $H_{j_1},\dots ,H_{j_p}$ in general position, set
$F(H_{j_1},\dots ,H_{j_p})=F_{{\alpha_0},\dots,{\alpha_p}}$
where
\[
X_{\alpha_0}=\C^k,\quad X_{\alpha_1}=H_{j_1},\quad \dots , \quad
X_{\alpha_p}=
H_{j_1} \cap \dots \cap H_{j_p} .
\]
Then def\/ine $\langle (H_{j_1},\dots ,H_{j_p}), F_{{\alpha_0},\dots,{\alpha_p}}
 \rangle = (-1)^{|\sigma|},$
if $F_{{\alpha_0},\dots,{\alpha_p}}
= F(H_{j_{\sigma(1)}},\dots ,H_{j_{\sigma(p)}})$ for some $\sigma \in S_p$,
and $\langle (H_{j_1},\dots ,H_{j_p}), F_{{\alpha_0},\dots,{\alpha_p}} \rangle = 0$ otherwise.

Denote by $\delta^{(a)} :  \FF^p(\A)\to\FF^{p-1}(\A)$
 the map dual to
$d^{(a)} :   \OS^{p-1}(\A) \to \OS^{p}(\A)$.
An element $v \in \FF^k(\A)$ is called  singular  if
$\delta^{(a)} v = 0$. Denote by
\[
\Sing \FF^k(\A) \subset \FF^k(\A)
\]
the subspace of all singular vectors.

\subsection{Contravariant map and form, see \cite{SV}}

 Weights  $a$ determine a contravariant  map
\[
\mathcal S^{(a)} : \ \FF^p(\A) \to \OS^p(\A),
\qquad
  {F}_{{\alpha_0},\dots,{\alpha_p}} \mapsto
\sum \ a_{j_1} \cdots a_{j_p}\ (H_{j_1}, \dots , H_{j_p}) ,
\]
 where the sum is taken over all $p$-tuples $(H_{j_1},\dots ,H_{j_p})$ such that
\[
H_{j_1} \supset X_{\al_1},\quad \dots,\quad H_{j_p}\supset X_{\alpha_p}  .
\]
Identifying $\OS^p(\A)$ with $\FF^p(\A)^*$, we consider
the map  as a bilinear form,
\[
S^{(a)} :  \ \FF^p(\A) \otimes \FF^p(\A) \to \C .
\]
The bilinear form is
called the contravariant form.
The contravariant form  is symmetric.
For $F_1, F_2 \in \FF^p(\A)$,
\[
S^{(a)}(F_1,F_2) =
\sum_{\{j_1, \dots , j_p\} \subset J}   a_{j_1} \cdots a_{j_p}
  \langle (H_{j_1}, \dots , H_{j_p}), F_1 \rangle
  \langle (H_{j_1}, \dots , H_{j_p}), F_2 \rangle   ,
\]
where the sum is over all unordered $p$-element subsets.

The contravariant form was introduced in~\cite{SV}. It is an analog of the Shapovalov form
in representation theory. On relations between the contravariant  and  Shapovalov forms, see~\cite{SV} and Section~\ref{sec Gaudin}.

\subsection{Remarks on generic weights}
\label{sec rem on gener wei}

\begin{thm}[\cite{SV}]
\label{thm Shap nondeg}
If weights $a$ are such that none of the dense edges has  weight zero, then
the contravariant form is nondegenerate.
\end{thm}

\begin{thm}\label{thm Y, O ,SV} \sloppy
If weights $a$ are such that none of the dense edges has  weight zero, then
  $H^p(\OS^*(\A),d^{(a)}) = 0$ for $ p < k$
and  $\dim H^k(\OS^*,d^{(a)}) = |\chi(U)|$, where $\chi (U)$ is the Euler
characteristics of $U$. In particular, these statements hold if all weights are positive.
\end{thm}

\begin{proof}
The theorem is proved in \cite{Y, OT2}.
It is also a straightforward corollary of some results in~\cite{SV}.
Namely, in~\cite{SV} a f\/lag complex
$d : \FF^p(\A) \to \FF^{p+1}(\A)$ was considered with  the dif\/ferential def\/ined by formula~(2.2.1) in~\cite{SV}. By~\cite[Corollary~2.8]{SV}
 the cohomology spaces of that f\/lag complex are trivial in all
degrees less than the top degree.
In~\cite{SV}, it was also proved that the contravariant map
 def\/ines a  homomorphism of the
f\/lag complex to the complex  $d^{(a)} : \OS^p(\A) \to \OS^{p+1}(\A)$.
Now Theorem~\ref{thm Y, O ,SV}
is a corollary of Theorem~\ref{thm Shap nondeg}.
\end{proof}

\begin{cor}
\label{lem on dim of sing}
If weights $a$ are such that none of the dense edges has  weight zero, then
the dimension of  $\Sing\FF^k(\A)$ equals $|\chi(U)|$.
\end{cor}

\subsection{Orlik--Solomon algebra as an algebra of dif\/ferential forms}

For  $j\in J$,  f\/ix a def\/ining equation for the hyperplane $H_j$, $f_{j} = 0$,
where $f_j$ is a polynomial of degree one in variables $t_1,\dots,t_k$.
Consider the logarithmic dif\/ferential form
$\omega_j = df_j/f_j$ on~$\C^k$.
Let $\bar{\OS}(\A)$ be the $\C$-algebra of dif\/ferential forms
generated by 1 and $\omega_j$, $j\in J$.
The map ${\OS}(\A) \to \bar{\OS}(\A)$, $(H_j) \mapsto \omega_j$,
is an isomorphism. We identify~${\OS}(\A)$ and~$\bar{\OS}(\A)$.

\subsection{Critical points of the master function}
\label{master function}

Given weights $a: J\to \C$, def\/ine the (multivalued) master function
$\Phi : U \to \C$,
\begin{gather}
\label{def mast fun}
\Phi  =  \Phi_{\A,a}  =  \sum_{j\in J} a_j \log f_j .
\end{gather}
Usually the function $e^\Phi=\prod_j f_j^{a_j}$ is called the master
function, see \cite{V1,V2,V4}, but it is more convenient
to work   with  def\/inition \eqref{def mast fun}.

A point $t\in U$ is a critical point of $\Phi$ if\ $d\Phi\vert_t = 0$.
We can rewrite this equation as
$\nu(a)\vert_t=0$ since
\begin{gather}
\label{omega(a)}
\nu(a) = d\Phi .
\end{gather}

\begin{thm}[\cite{V2, OT1, Si}]
\label{thm V,OT,S}
For generic exponents $a$ all critical points of
$\Phi$ are nondegenerate and the number of critical points equals
  $|\chi(U)|$.
\end{thm}

Denote  $\C(t)_U$ the algebra of rational functions on $\C^k$ regular on $U$ and
\[
I_\Phi
 =\langle
\frac{\partial \Phi} {\partial t_i }
\ |\ i=1,\dots,k\ \rangle  \subset  \C(t)_U
\]
the ideal generated by f\/irst  derivatives of $\Phi$.
Let $
A_\Phi = \C(t)_U/ I_\Phi
$
 be
the algebra of functions on the critical set  and
$[ \ ]$: $\C(t)_U \to  A_\Phi$, $f\mapsto [f]$, the canonical homomorphism.

If all critical points are isolated,
then the critical set  is f\/inite and the algebra
$ A_\Phi$ is f\/inite-dimensional. In that case, $ A_\Phi$
is the direct sum of local algebras corresponding to points $p$ of
the critical set,
\begin{gather*}
 A_\Phi = \oplus_{p}  A_{p,\Phi}  .
\end{gather*}
The local algebra
$A_{p,\Phi}$ can be def\/ined as the quotient of the algebra of germs at $p$
 of holomorphic functions modulo the ideal
$ I_{p,\Phi}$ generated f\/irst derivatives of $\Phi$.
Denote by $ \frak m_p\subset A_{p,\Phi}$ the maximal ideal   generated by germs
of functions equal to zero at $p$.

\begin{Lemma}
\label{lem f_j generate}
The elements
$[1/f_j]$, $j\in J$, generate $ A_\Phi$.
\end{Lemma}

\begin{proof}
If $H_{j_1}, \dots, H_{j_k}$ intersect transversally, then
$1/f_{j_1}, \dots, 1/f_{j_k}$ form a coordinate system on $U$.
This remark proves the lemma.
\end{proof}

Def\/ine a rational function $\Hess : \C^k \to \C$, regular on $U$, by the formula
\[
\Hess (t)  = \det_{1\leq i,j \leq k}
\left( \frac{\partial^2\Phi} {\partial t_i \partial t_j} \right) (t)  .
\]
The function is called the Hessian of $\Phi$.

Let $p$ be an isolated critical point of $\Phi$.
Denote by $[\Hess]_p$ the image
of the Hessian in
$A_{p,\Phi}$. It is known that the image is nonzero and
the one-dimensional subspace $\C [\Hess]_p\subset \AT$ is the
annihilating ideal
 of the maximal ideal
$\frak m_p \subset A_{p,\Phi}$.

Let
 $\rho_{p} :  A_{p,\Phi} \to \C$, be the
Grothendieck residue,
\begin{gather*}
f   \mapsto   \frac 1{(2\pi i)^k}\,\Res_{p}
  \frac{ f}{\prod\limits_{i=1}^k  \frac{\der \Phi}{\der t_i}}
=\frac{1}{(2\pi i)^k}\int_{\Gamma_p}
\frac{f\ dt_1\wedge\dots\wedge dt_k}{\prod\limits_{s=1}^k \frac{\der \Phi}{\der t_i}}  ,
\end{gather*}
where
$\Gamma_p=\{(t_1,\dots,t_k)\ |\ |\frac{\der \Phi}{\der t_i}|=\epsilon_i,\ i=1,\dots,k\}$ is
the real $k$-cycle oriented by the condition
\[
d\arg \frac{\der \Phi}{\der t_1}\wedge\dots\wedge d\arg \frac{\der \Phi}{\der t_k}\ge0  ;
\]
here $\epsilon_s$ are suf\/f\/iciently small positive real numbers, see~\cite{GH}.

It is known that
$\rho_{p} : [\Hess]_p  \mapsto \mu_p$,
where $\mu_p = \dim_\C  A_{p,\Phi}$ is the Milnor number of the critical point
$p$.
Let $(\,,\,)_{p}$ be the residue
bilinear form,
\begin{gather}
\label{Gr form}
( f, g)_{p}  =  \rho_{p} (f g)  .
\end{gather}
That form is nondegenerate.

\subsection[Special vectors in $\FF^k(\A)$ and canonical element]{Special vectors in $\boldsymbol{\FF^k(\A)}$ and canonical element}
\label{Special}

A dif\/ferential top degree form $\eta \in {\OS}^k(\A)$ can be written as
\[
\eta  =  f  dt_1 \wedge \dots \wedge dt_k,
\]
where $f$ is a rational function on $\C^k$,  regular on $U$.

Def\/ine a rational map $v : \C^k \to \FF^k(\A)$ regular on $U$.
For $t\in U$,  set $v(t)$ to be the element of~$\FF^k(\A)$ such that
\[
\langle  \eta, v(t)  \rangle  =  f(t)
\qquad
{\mathrm{for\ any\ }} \ \eta \in \OS^k(\A)  .
\]
The map $v$ is called the specialization map and its value~$v(t)$
is called the special vector at $t\in U$, see~\cite{V4}.

Let $(F_m)_{m\in M}$ be a basis of $\FF^k(\A)$ and $(H^m)_{m\in M} \subset \OS^k(\A)$ the dual basis.
Consider the element $\sum_m H^m\otimes F_m \in \OS^k(\A) \otimes \FF^k(\A)$.
We have $H^m = f^mdt_1\wedge\dots\wedge dt_k$ for some $f^m \in \C(t)_U$.
The element
\begin{gather*}
E = \sum_{m\in M} f^m \otimes F_m   \in  \C(t)_U\otimes \FF^k(\A)
\end{gather*}
will be called the canonical element of the arrangement $\A$. It  does not depend on the choice of the basis $(F_m)_{m\in M}$.

For any $t\in U$, we have
\begin{gather*}
v(t) = \sum_{m\in M} f^m(t) F_m .
\end{gather*}

Denote by $[E]$ the image of the canonical element in $ A_\Phi\otimes\FF^k(\A)$.

\begin{Lemma}
\label{lem v sing}
We have $[E] \in  A_\Phi \otimes \Sing \FF^k(\A)$.
\end{Lemma}

\begin{proof}
By formula \eqref{omega(a)}, if $e \in d^{(a)}\cdot \OS(\A)^{k-1}$, then
$\langle 1\otimes e, E\rangle \in I_\Phi$. Hence
\begin{gather}
\label{frak}
\langle 1\otimes e, [E]\rangle = 0
\end{gather}
in $A_\Phi$.
Let $e_1,\dots,e_l$ be a basis of $d^{(a)}\cdot \OS^{k-1}(\A)$. Extend it to a basis
$e_1,\dots,e_l, e_{l+1}, \dots, e_{|M|}$ of
$\OS^{k}(\A)$. Let
$e^1,\dots,e^l, e^{l+1},\dots,e^{|M|}$ be the dual basis of $\FF^k(\A)$.
Then $e^{l+1},\dots,e^{|M|}$ is a basis of $\Sing \FF^k(\A)$.
Let $[E] = \sum\limits_{i=1}^{|M|} [g^i]\otimes e^i$ for some $[g^i]\in A_\Phi$.
By~\eqref{frak} we have
$[g^i]=0$ for all $i\leq l$.
\end{proof}

\begin{thm}[\cite{V4}]
\label{thm crit sing}
 A point $t\in U$ is a critical point of $\Phi$, if and only if the special vector
$v(t)$ is a singular vector.
\end{thm}

\begin{proof}
The theorem follows from Lemma \ref{lem v sing}.
\end{proof}

\begin{thm}[\cite{V4}]
\label{first theorem}\qquad
\begin{enumerate}\itemsep=0pt
\item[$(i)$]
For any  $t\in U$,
\[
S^{(a)}(v(t),v(t))   =   (-1)^k\ \Hess (t)  .
\]
\item[$(ii)$] If $t^1, t^2 \in U$ are different isolated critical points of $\Phi$,
then the special singular vec\-tors~$v(t^1)$,~$v(t^2)$ are orthogonal,
\[
S^{(a)}\big(v\big(t^1\big),v\big(t^2\big)\big) = 0  .
\]
\end{enumerate}
\end{thm}

\subsection{Arrangements with normal crossings only}
\label{An arrangement with normal crossings only}

 An essential arrangement $\A$ is with normal crossings only,
if exactly $k$ hyperplanes meet at every vertex of $\A$.
Assume that $\A$ is an essential arrangement with normal crossings only.

A subset $\{j_1,\dots,j_p\}\subset J$ will be called independent if the hyperplanes
$H_{j_1},\dots,H_{j_p}$ intersect transversally.

A basis of $\OS^p(\A)$ is formed by
$(H_{j_1},\dots,H_{j_p})$ where
$\{{j_1} <\dots <{j_p}\}$  are independent ordered $p$-element subsets of
$J$. The dual basis of $\FF^p(\A)$ is formed by the corresponding vectors
$F(H_{j_1},\dots,H_{j_p})$.
These bases of $\OS^p(\A)$ and $\FF^p(\A)$ will be called standard.

In $\FF^p(\A)$ we  have
\begin{gather}
\label{skew}
F(H_{j_1},\dots,H_{j_p}) = (-1)^{|\sigma|}
F(H_{j_{\sigma(1)}},\dots,H_{j_{\sigma(p)}})
\end{gather}
for any permutation  $\sigma \in S_p$.

For an independent subset $\{j_1,\dots,j_p\}$, we have
\[
S^{(a)}(F(H_{j_1},\dots,H_{j_p}) , F(H_{j_1},\dots,H_{j_p})) = a_{j_1}\cdots a_{j_p}
\]
and
\[
S^{(a)}(F(H_{j_1},\dots,H_{j_p}) , F(H_{i_1},\dots,H_{i_k})) = 0
\]
for distinct elements of the standard basis.

\subsection[Real structure on $\OS^p(\A)$ and $\FF^p(\A)$]{Real structure on $\boldsymbol{\OS^p(\A)}$ and $\boldsymbol{\FF^p(\A)}$}
\label{sec real structure}

We have def\/ined  $\OS^p(\A)$ and $\FF^p(\A)$ as vector spaces over $\C$. But one can def\/ine
 the corresponding spaces over the f\/ield $\R$ so that
$\OS^p(\A) = \OS^p(\A)_\R\otimes_\R \C$ and
$\FF^p(\A) = \FF^p(\A)_\R\otimes_\R \C$. If all weights~$a$ are real, then the dif\/ferential
$d^{(a)} : \OS^p(\A) \to \OS^{p+1}(\A) $ preserves the real subspaces and one can def\/ine
the subspace of singular vectors $\Sing \FF^k(\A)_\R \subset \FF^k(\A)_\R$ so that
$\Sing \FF^k(\A) = \Sing \FF^k(\A)_\R \otimes _\R\C$.

\subsection{A real arrangement with positive weights}
\label{sec real arr}

Let $t_1,\dots,t_k$ be standard coordinates on $\C^k$. Assume that  every polynomial
$f_j$, $j\in J$, has real coef\/f\/icients,
\[
f_j = z_j+b_j^1t_1 + \dots + b_j^kt_k ,
\]
where $z_j$, $b_j^i$ are real numbers.

Denote $U_\R = U\cap \R^k$. Let $U_\R = \cup_\al D_\al$ be the decomposition into the union
of connected components. Each connected component is a convex polytope.
It is known that the number of
bounded connected components equals $|\chi(U)|$, see \cite{Z}.

\begin{thm}[\cite{V2}]
\label{thm real crit pts}
Assume that  weights $(a_j)_{ j\in J}$ are positive. Then the union of
all critical points of the master function $\Phi_{\A,a}$
is contained in the union
of all bounded components of $U_\R$. Each bounded component contains
exactly one critical point. All critical points are nondege\-ne\-rate.
\end{thm}

\begin{cor}
\label{cor real basis}
Under assumptions of Theorem~{\rm \ref{thm real crit pts}} let
 $t^1,\dots,t^d \in U_\R$ be a list of all distinct critical points of the
master function $\Phi_{\A,a}$. Then the corresponding special vectors
$v(t^1),\dots,v(t^d)$ form a basis of
$\Sing \FF^k(\A)_\R$.
That basis is orthogonal with respect to the contravariant form~$S^{(a)}$.
\end{cor}

Note that the contravariant form on $\Sing \FF^k(\A)_\R$ is positive def\/inite.

\subsection{Resolution of a hyperplane-like divisor}

Let $Y$ be a smooth complex compact manifold of dimension $k$,
$D$ a divisor. The divisor $D$ is hyperplane-like
if $Y$ can be covered by coordinate charts such that in each chart
$D$ is a union of
hyperplanes. Such charts will be called linearizing.
Let $D$ be a hyperplane-like divisor, $U$~be
a linearizing chart. A local edge of $D$ in $U$ is
any nonempty irreducible intersection in $U$ of hyperplanes of~$D$ in~$U$.
An edge of $D$ is the
maximal analytic continuation in $Y$ of a local edge.
Any edge is an immersed submanifold
in~$Y$. An edge is called dense if it is locally dense.
For $0\leq i\leq k-2$, let $\mc L_i$
be the collection of all dense edges of $D$ of dimension $i$. The
following theorem is essentially contained in Section 10.8 of \cite{V1}.

\begin{thm}[\cite{STV}]
\label{thm resolution}
 Let $W_0 = Y$. Let $\pi_1:W_1\to W_0$
be the blow up along points in $\mc L_0$.
In general, for $1\leq  s \leq k-1$, let
$\pi_s:W_s\to W_{s-1}$
 be the blow up along the proper transforms
of the $(s-1)$-dimensional dense edges in
$\mc L_{s-1}$ under $\pi_1 \cdots \pi_{s-1}$.
Let $\pi = \pi_1 \cdots \pi_{k-1}$.
Then $W = W_{n-1}$ is nonsingular and $\pi^{-1}(D)$ has normal crossings.
 \end{thm}

\section{A family of parallelly translated hyperplanes}
\label{sec par trans}

\subsection[An arrangement in  $\C^n\times\C^k$]{An arrangement in  $\boldsymbol{\C^n\times\C^k}$}
\label{An arrangement in}
Recall that $J=\{1,\dots,n\}$.
Consider $\C^k$ with coordinates $t_1,\dots,t_k$,
$\C^n$ with coordinates $z_1,\dots,z_n$, the projection
$\C^n\times\C^k \to \C^n$.

Fix $n$ nonzero linear functions on $\C^k$,
\[
g_j = b_j^1t_1+\dots + b_j^kt_k,\qquad j\in J,
\]
where $b_j^i\in \C$.
Def\/ine $n$ linear functions on $\C^n\times\C^k$,
\[
f_j = z_j+g_j = z_j + b_j^1t_1+\dots + b_j^kt_k,\qquad j\in J.
\]
In $\C^n\times \C^k$ def\/ine
 an arrangement
\[
\tilde \A = \{ \tilde H_j\ | \ f_j = 0, \ j\in J \} .
\]
Denote $\tilde U = \C^n\times \C^k - \cup_{j\in J} \tilde H_j$.

For every f\/ixed $z=(z_1,\dots,z_n)$ the arrangement $\tilde \A$
induces an arrangement $\A(z)$ in the f\/iber over $z$ of the projection. We
identify every f\/iber with $\C^k$. Then $\A(z)$ consists of
hyper\-pla\-nes~$H_j(z)$, $j\in J$, def\/ined in~$\C^k$ by the same equations
$f_j=0$. Denote
\begin{gather*}
U(\A(z)) = \C^k - \cup_{j\in J} H_j(z)  ,
\end{gather*}
the complement to the arrangement $\A(z)$.

In the rest of the paper we assume that for any $z$ the arrangement $\A(z)$ has a vertex. This means that
the span of $(g_j)_{j\in J}$ is $k$-dimensional.

A point $z\in\C^n$ will be called {\it good} if $\A(z)$ has normal
crossings only.  Good points form the complement in $\C^n$ to the union
of suitable hyperplanes called the discriminant.

\subsection{Discriminant}
\label{Discr}

The collection $(g_j)_{j\in J}$ induces a
matroid structure on $J$.  A subset $C=\{i_1,\dots,i_r\}\subset J$ is
a~circuit  if $(g_i)_{i\in C}$ are linearly dependent but any
proper subset of $C$ gives linearly independent~$g_i$'s.

For a circuit $C=\{i_1,\dots,i_r\}$,   let
$(\la^C_i)_{i\in C}$ be a nonzero collection of complex numbers such that
$\sum_{i\in C}
\la^C_ig_i = 0$. Such a collection  is unique up to
multiplication by a nonzero number.

For every circuit $C$ f\/ix such a collection
and denote $f_C = \sum_{i\in C} \la^C_iz_i$.
The equation $f_C=0$ def\/ines a hyperplane $H_C$ in
$\C^n$.
It is convenient to assume that $\la^C_i=0$ for $i\in J-C$ and write
$f_C = \sum_{i\in J} \la^C_iz_i$.

For any $z\in\C^n$, the hyperplanes $(H_i(z))_{i\in C}$ in $\C^k$ have nonempty
intersection if and only if $z\in H_C$. If $z\in H_C$, then the
intersection has codimension $r-1$ in~$\C^k$.

Denote by $\frak C$ the set of all circuits in $J$.
Denote  $\Delta = \cup_{C\in \frak C} H_C$.

\begin{Lemma}
The arrangement $\A(z)$ in $\C^k$
has normal crossings only, if and only if $z\in \C^n-\Delta$.
\end{Lemma}

\begin{rem}
\label{rem g real}\rm
If all linear functions
 $g_j, j\in J$, are real, then for any circuit $C\in \frak C$ the numbers
$(\la^C_i)_{i\in C}$ can be chosen to be real. Therefore,
in that case every hyperplane $H_C$ is real.
\end{rem}

\subsection{Good f\/ibers}
\label{sec Good fibers}

For any
$z^1, z^2\in \C^n-\Delta$, the spaces $\FF^p(\A(z^1))$, $\FF^p(\A(z^2))$
 are canonically identif\/ied. Namely, a~vector $F(H_{j_1}(z^1),\dots,H_{j_p}(z^1))$ of the f\/irst space
is identif\/ied  with
the vector
 $F(H_{j_1}(z^2),\dots$, $H_{j_p}(z^2))$ of the second.

Assume that weights $a=(a_j)_{j\in J}$ are given and all of them are nonzero. Then each
arrangement $\A(z)$  is weighted.
The identif\/ication of spaces $\FF^p(\A(z^1))$,
$\FF^p(\A(z^2))$ for $z^1,z^2\in\C^n-\Delta$ identif\/ies the corresponding subspaces
$\Sing\FF^k(\A(z^1))$, $\Sing\FF^k(\A(z^2))$ and contravariant forms.

For a point $z\in\C^n-\Delta$, denote $V=\FF^k(\A(z))$, $\Sing V=\Sing\FF^k(\A(z))$.
The triple $(V, \Sing V, S^{(a)})$ does not depend on  $z\in\C^n-\Delta$
under the above identif\/ication.

\subsection{Bad f\/ibers}
\label{sec Bad fibers}

Points of  $\Delta\subset \C^n$ will be called bad.

Let $z^0\in\Delta$ and $z\in \C^n-\Delta$.
By def\/inition, for any $p$
the space $\OS^p(\A(z^0))$ is obtained from $\OS^p(\A(z))$ by
adding new relations.
Hence $\OS^k(\A(z^0))$ is canonically identif\/ied
with a~quotient space of $V^* = \OS^k(\A(z))$ and
$\FF^p(\A(z^0))$ is canonically identif\/ied
with a subspace
of $V = \FF^p(\A(z))$.

Let us consider $\FF^k(z^0)$ as a subspace of $V$. Let
$S^{(a)}\vert_{\FF^k(z^0)}$ be the restriction of the contravariant form on $V$ to that subspace.
Let $S^{(a)}(z^0)$ be the contravariant form on $\FF^k(\A(z^0))$ of the arrangement $\A(z^0)$.

\begin{Lemma}
\label{lem identif 2}
Under the above identifications,  $S^{(a)}\vert_{\FF^k(z^0)} = S^{(a)}(z^0)$.
\end{Lemma}

\section{Conservation of the number of critical points}
\label{sec Preservation}

Let $\A=(H_j)_{j\in J}$ be an essential arrangement in
$\C^k$ with weights $a$.  Consider its compactif\/ication in the
projective space $\Pee^k$ containing $\C^k$. Assign
the weight $a_\infty=-\sum_{j\in J} a_j$ to the hyperplane
$H_\infty=\Pee^k-\C^k$ and denote by  $\A^\vee$
the arrangement $(H_j)_{j\in J\cup \infty}$ in $\Pee^k$.

The weighted arrangement $(\A, a)$ will be called unbalanced
if the weight of any dense
edge of $\A^\vee$ is nonzero.

For example, if all weights $(a_j)_{j\in J}$ are positive, then the
weighted arrangement $(\A, a)$ is unbalanced.
Clearly, the
unbalanced weights form a Zarisky open subset in the space of all weights of $\A$.

\begin{Lemma}
\label{lem crit 1} If $(\A, a)$ is unbalanced, then
all critical points of the master function of the weighted arrangement
$(\A,a)$ are isolated.
\end{Lemma}

\begin{proof}
Let $\pi : W\to \Pee^k$ be the resolution (described in Theorem \ref{thm resolution})
of singularities of the divisor
$D=\cup_{j\in J\cup \infty}H_j$.
Let $\Phi_a$ be the master function of $(\A,a)$.
Then locally on $W$ the function
$\pi^{-1}\Phi_a$
has the form
\[
\pi^{-1}\Phi_a  =  \sum_{i=1}^m \al_i\log u_i  + \log \phi(u_1,\dots,u_k),
\qquad
\phi(0,\dots,0)\neq 0 .
\]
Here $u_1,\dots,u_k$ are local coordinates, $0\leq m\leq k$, the function
$\phi(u_1,\dots,u_k)$
is holomorphic at $u=0$,
the equation
$u_1\cdots u_m=0$ def\/ines $\pi^{-1}(D)$ in this  chart.
 If the image of a divisor $u_i=0$, $1\leq i\leq m$,
under the map $\pi$ is  an $s$-dimensional
piece of an $s$-dimensional dense edge of $\A^\vee$,
then $\al_i$ equals the weight of that edge. In particular, $\al_i$, $i=1,\dots,m$,
are all nonzero.

Let $U(\A)=\C^k-\cup_{j\in J} H_j$.
The critical point equations of $\pi^{-1}\Phi_a$ on $\pi^{-1}(U(\A))$ are
\begin{gather}
\label{Crit pt eqns}
\frac {\al_i }{u_i} + \frac{1}{\phi(u)}
\frac{\partial \phi}{\partial u_i}(u) = 0  ,
\qquad i=1,\dots,m ,
\\
\frac{1}{\phi(u)}
\frac{\partial \phi}{\partial u_i}(u) = 0  ,
\qquad i=m+1,\dots,k .
\notag
\end{gather}
If the critical set of $\pi^{-1}\Phi_a$ on  $\pi^{-1}(U(\A))$ is inf\/inite,
then it contains an algebraic curve. The closure of that curve
must intersect $\pi^{-1}(D)$. But equations~\eqref{Crit pt eqns} show that
this is impossible.
\end{proof}

Denote by $\mu(\A,a)$ the sum of Milnor numbers of all of the critical points of
$\Phi_a$ on $U(\A)$.

\begin{Lemma}\label{lem crit 2}
If $(\A,a)$ is unbalanced, then  $\mu(\A,a)=|\chi(U)|$.
\end{Lemma}

\begin{proof}
Assume that $a(s)$, $s\in [0,1]$, is a continuous family of unbalanced weights~of $\A$. Then $\mu(\A,a(s))$ does not depend on~$s$. Indeed,
using equations~\eqref{Crit pt eqns} one shows that the critical
points cannot approach $\pi^{-1}(D)$ as $s\in [0,1]$ changes.
For generic weights $a$ we have $\mu(\A,a(s))=|\chi(U)|$ by Theorem~\ref{thm V,OT,S}.
Hence, $\mu(\A,a)=|\chi(U)|$ for any unbalanced weights~$a$.
\end{proof}

\section{Hamiltonians of good f\/ibers}
\label{Hamiltonians of good fibers}

\subsection{Construction}
\label{Construction}
Consider the master function
\begin{gather*}
\Phi(z,t) = \sum_{j\in J} a_j \log f_j(z,t)
 \end{gather*}
as a function on $\tilde U\subset \C^n\times \C^k$.

Let $\kappa$ be a nonzero complex number.
The function $e^{\Phi(z,t)/\kappa}$
def\/ines a rank one local system~$\mc L_\kappa$
on $\tilde U$ whose horizontal sections
over open subsets of $\tilde U$
 are univalued branches of $e^{\Phi(z,t)/\kappa}$ multiplied by complex numbers.

For a f\/ixed $z$, choose any $\gamma\in H_k(U(\A(z)), \mc L_\kappa\vert_{U(\A(z))})$.
The linear map
\[
\{\gamma\}  : \   \OS^k(\A(z)) \to \C , \qquad
\om \mapsto \int_{\gamma} e^{\Phi(z,t)/\kappa} \om   ,
\]
is an element of $\Sing\FF^k(\A(z))$ by Stokes' theorem.

It is known that for generic  $\kappa$
any element of $\Sing\FF^k(\A(z))$ corresponds to a certain
$\gamma$ and  in that case the integration  identif\/ies $\Sing\FF^k(\A(z))$ and
$H_k(U(\A(z)), \mc L_\kappa\vert_{U(\A(z))})$, see~\cite{SV}.

The vector bundle
\[
\cup_{z\in \C^n-\Delta}\,H_k(U(\A(z)), \mc L_\kappa\vert_{U(\A(z))})  \to
 {} \ \C^n-\Delta
 \]
has a canonical  (f\/lat) Gauss--Manin connection.
The  Gauss--Manin connection
induces a f\/lat connection
on the trivial bundle
$\C^n\times \Sing V\to \C^n$
with singularities over the discriminant $\Delta \subset \C^n$.
That connection will be called the Gauss--Manin connection as well.

\begin{thm}
\label{thm ham normal}
The Gauss--Manin differential equations for horizontal sections of
the connection on $\C^n\times \Sing V\to \C^n$ have the form
\begin{gather*}
\kappa \frac{\der I}{\der z_j}(z) = K_j(z)I(z),
\qquad
j\in J,
\end{gather*}
where $I(z)\in \Sing V$ is a horizontal section,
 $K_j(z)$: $V\to V$, $j\in J,$ are suitable linear operators  preserving $\Sing V$ and
independent on
$\kappa$. For every $j$,
the operator $K_j(z)$ is a rational function of $z$ regular on $\C^n-\Delta$.
Each operator is symmetric
with respect to the contravariant form $S^{(a)}$.
\end{thm}

Theorem \ref{thm ham normal} is proved in Section~\ref{Proofs}.
A formula for $K_j(z)$ see in~\eqref{K_j}.

The f\/latness of the connection for all $\kappa$ implies
the commutativity of the operators,
\[
K_i(z)|_{\Sing V}K_j(z)|_{\Sing V}=K_j(z)|_{\Sing V}K_i(z)|_{\Sing V}
\qquad
\text{ for all} \ i,j \ \text{and}\ z\in\C^n-\Delta .
\]

Let $V^*$ be the space dual to $V$. If $M:V\to V$ is a linear operator, then
 $M^*:V^*\to V^*$  denotes the dual operator.
Let
$W \subset V^*$ be the image of $V$ under the map
$ V \to  V^*$ associated with the contravariant form and
$\sing W\subset W$ the image of $\sing V$.
The contravariant form induces on $W$ a nondegenerate
symmetric bilinear form, also denoted by~$S^{(a)}$.

\begin{Lemma}
\label{lem Main}
For   $z\in \C^n-\Delta$, the operators $K_j(z)^*$ preserve the subspaces $\sing W\subset W\subset V^*$. The operators
$K_j(z)^*|_{W} :W\to W$ are symmetric with respect to the contravariant form. The operators
$K_j(z)^*|_{\sing W} :\sing W\to \sing W$, $j\in J$, commute.
\end{Lemma}

For   $z\in \C^n-\Delta$, we def\/ine {\it the quantum integrable model assigned to}
$(\A(z), a)$ to be the collection
\begin{gather}
\label{Ham good fibers}
\big(\sing W;\ S^{(a)}|_{\sing W};\ K_1(z)^*|_{\sing W},\dots, K_n(z)^*|_{\sing W} :\  \sing W\to \sing W\big).
\end{gather}
The unital subalgebra of $\End(\sing W)$ generated by operators
$K_1(z)^*|_{\sing W} , \dots,
 K_n(z)^*|_{\sing W}$ will be called {\it the
algebra of geometric Hamiltonians of $(\A(z),a)$}.

If the contravariant form $S^{(a)}$ is nondegenerate  on $V$, then this model is isomorphic to the collection
\[
\big(\sing V;\ S^{(a)}|_{\sing V};\ K_1(z)|_{\sing V},\dots, K_n(z)|_{\sing V} : \ \sing V\to \sing V\big).
\]

It is clear that any weighted
essential arrangement with normal crossings only can be realized as a good f\/iber of such
a construction. Thus,  every weighted
essential arrangement with normal crossings only
is provided with a quantum integrable model.

\subsection{Key identity (\ref{formula H})}
\label{sec key identity}

For any circuit $C=\{i_1, \dots, i_r\}\subset J$, let us
def\/ine a linear operator
$L_C : V\to V$ in terms of the standard basis of
 $V$, see Section~\ref{An arrangement with normal crossings only}.

For $m=1,\dots,r$, def\/ine $C_m=C-\{i_m\}$.
Let $\{{j_1}<\dots <{j_k}\}
\subset J$ be an independent ordered subset and
$F(H_{j_1},\dots,H_{j_k})$
the corresponding element of the standard basis.
Def\/ine $L_C : F(H_{j_1},\dots,H_{j_k}) \mapsto 0$ if
$|\{{j_1},\dots,{j_k}\}\cap C| < r-1$.
If $\{{j_1},\dots,{j_k}\}\cap C = C_m$ for some
$1\leq m\leq r$, then using the skew-symmetry property~\eqref{skew}
we can write
\[
F(H_{j_1},\dots,H_{j_k})
 =
\pm\, F(H_{i_1},H_{i_2},\dots,\widehat{H_{i_{m}}},\dots,H_{i_{r-1}}H_{i_{r}},H_{s_1},\dots,H_{s_{k-r+1}})
\]
with $\{{s_1},\dots,{s_{k-r+1}}\}=
\{{j_1},\dots,{j_k}\}-C_m$.
Def\/ine
\begin{gather*}
L_C
 : \
F(H_{i_1},\dots,\widehat{H_{i_{m}}},\dots,H_{i_{r}},H_{s_1},\dots,H_{s_{k-r+1}})\\
\qquad \ {}
 \mapsto
 (-1)^m \sum_{l=1}^{r} (-1)^l a_{i_l}
F(H_{i_1},\dots,\widehat{H_{i_{l}}},\dots,H_{i_{r}},H_{s_1},\dots,H_{s_{k-r+1}}) .
\end{gather*}

\begin{Lemma}
\label{lem L_C}
The map $L_C$ is symmetric with respect to the
contravariant form.
\end{Lemma}

\begin{proof}\sloppy
For $l=1,\dots,r$,
 denote $F_l=F(H_{i_1},\dots,\widehat{H_{i_{l}}},\dots,H_{i_{r}},H_{s_1},\dots,H_{s_{k-r+1}})$.
It is clear that $L_C$ is symmetric with respect to the contravariant form if and only if
$S^{(a)}(L_CF_l,F_m)=S^{(a)}(F_l,L_CF_m)$ for all $1\leq l,m\leq r$.
But both sides of this expression are equal to
$(-1)^{l+m} a_{i_1}\cdots a_{i_r}a_{s_1}\cdots a_{s_{k-r+1}}$.
\end{proof}

On $\C^n\times\C^k$ consider the logarithmic dif\/ferential 1-forms
\[
\omega_j = \frac {df_j}{f_j}, \ j\in J,
\qquad
\omega_C = \frac {df_C}{f_C}, \qquad C\in \frak C.
\]
For any circuit $C=\{i_1,\dots,i_r\}$, we have
\[
\omega_{i_1} \wedge \dots \wedge \omega_{i_{r}} =
\omega_C \wedge \sum_{l=1}^{r} (-1)^{l-1}
\omega_{i_1} \wedge \dots \wedge \widehat{\omega_{i_{l}}} \wedge \dots \wedge
\omega_{i_{r}} .
\]

\begin{Lemma}
\label{lem 2}
We have
\begin{gather}
\sum_{{\rm independent } \atop \{j_1 < \dots < j_k\} \subset J }
\bigg( \sum_{j\in J} a_j
 \omega_j  \bigg) \wedge \omega_{j_1} \wedge \dots \wedge \omega_{j_k}
\otimes
F(H_{j_1}, \dots , H_{j_k}) \nonumber
\\
\qquad{} =
\sum_{{\rm independent } \atop \{j_1 < \dots < j_k\} \subset J }
\sum_{C\in \frak C}
\omega_{C} \wedge
 \omega_{j_1} \wedge \dots \wedge \omega_{j_k}
\otimes
L_C
F(H_{j_1}, \dots , H_{j_k}) .\label{formula H}
\end{gather}
\end{Lemma}

\begin{proof}
The lemma is a direct corollary of the def\/inition of maps $L_C$.
\end{proof}

Identity \eqref{formula H} is a key formula of this paper. Identity \eqref{formula H}
 is an analog of the key Theo\-rem~7.2.5 in~\cite{SV}
and it is a generalization
of the identity of Lemma~4.2 in~\cite{V4}.

\subsection{An application of the key identity (\ref{formula H}) -- proof of Theorem~\ref{thm ham normal}}
\label{Proofs}

Fix
$\kappa \in \C^\times$, $z\in \C^n-\Delta$, $\gamma \in  H_k(U(\A(z)), \mc L_\kappa\vert_{U(\A(z))})$.
Let $\{\gamma\} : \OS^k(\A(z)) \to \C$, $
\om \mapsto \int_{\gamma} e^{\Phi(z,t)/\kappa} \om$, be the corresponding element of~$V$.
We have
\[
\{\gamma\} =
\sum_{{\rm independent } \atop \{j_1 < \dots < j_k\} \subset J }
\left( \int_{\gamma} e^{\Phi(z,t)/\kappa}  \omega_{j_1} \wedge \dots \wedge \omega_{j_k}\right)
F(H_{j_1},\dots,H_{j_k})  .
\]
Let $z\mapsto \gamma(z) \in H_k(U(\A(z)), \mc L_\kappa\vert_{U(\A(z))})$ be a locally constant
section of the
Gauss--Manin connection.
Then
\[
\{\gamma(z)\} =
\sum_{{\rm independent } \atop \{j_1 < \dots < j_k\} \subset J }
 \left(
\int_{\gamma(z)} e^{\Phi(z,t)/\kappa}  \omega_{j_1} \wedge \dots \wedge \omega_{j_k}\right)
F(H_{j_1},\dots,H_{j_k}) .
\]

\begin{Lemma}
\label{lem diff}
The differential of the function $\{\gamma(z)\}$ is given by the formula
\[
\kappa
d\{\gamma(z)\}  =
\sum_{C\in \frak C}
 L_C\{\gamma(z)\}  \omega_{C}.
 \]
\end{Lemma}

\begin{proof}\sloppy
The lemma follows from identity~\eqref{formula H}
and the formula of dif\/ferentiation of an integ\-ral.
\end{proof}

\begin{Lemma}
\label{lem L_C preserve}
For every circuit $C$, the operator $L_C$ preserves the subspace
$\Sing V$.
\end{Lemma}

\begin{proof}
The values of the function $\{\gamma(z)\}$ belong to $\Sing V$. Hence, the values of its derivatives
belong to $\Sing V$. Now the
 lemma follows from Lemma \ref{lem diff}.
\end{proof}

Recall that $\omega_C = df_C/f_C$ and $f_C = \sum_{j\in J}\la^C_jz_j$. Denote
\begin{gather}
\label{K_j}
K_j(z)   =   \sum_{C\in \frak C}
 \frac{\la_j^C}{f_C(z)}  L_C  ,
\qquad
j\in J .
\end{gather}
Then
\begin{gather}
\label{LK}
\sum_{C\in \frak C}
\omega_{C} \otimes L_C = \sum_{j\in J} dz_j\otimes K_j(z) .
\end{gather}

\begin{Lemma}
\label{lem diff eqn}
Let $z\mapsto \gamma(z) \in H_k(U(\A(z)), \mc L_\kappa\vert_{U(\A(z))})$ be a locally constant
section of the
Gauss--Manin connection. Then
\begin{gather}
\label{GM}
\kappa \frac{\der}{\der z_j}\{\gamma(z)\} = K_j(z)\{\gamma(z)\} ,
\qquad
j\in J .
\end{gather}
\end{Lemma}

Lemmas \ref{lem L_C}, \ref{lem L_C preserve}, \ref{lem diff eqn} prove Theorem
\ref{thm ham normal}.

\subsection{Another application of the key identity (\ref{formula H})}
\label{More corollaries of key identity}

Recall that $\tilde U$ is the complement to the union of hyperplanes $(\tilde H_j)_{j\in J}$
in $\C^n\times \C^k$, see Sec\-tion~\ref{An arrangement in}.
Denote by $\C(z,t)_{\tilde U}$ the algebra of rational functions on $\C^n\times \C^k$ regular on
$\tilde U$.

For any basis vector $(H_{j_1},\dots,H_{j_k})$ of $V^*$, let us write
\begin{gather*}
\omega_{j_1} \wedge \dots \wedge \omega_{j_k} =
f_{j_1,\dots,j_k}(z,t) dt_1\wedge\dots\wedge dt_k
+ z\text{-part},
\end{gather*}
where $ f_{j_1,\dots,j_k} \in \C(z,t)_{\tilde U}$ and the $z$-part is a dif\/ferential form with zero
restriction to any f\/iber of the projection $\C^n\times \C^k \to \C^n$ (in
coordinates $t_1,\dots,t_k,z_1,\dots,z_n$,
that form  has at least one of $dz_1,\dots,dz_n$ as factors
 in each of its summands).
Def\/ine the canonical element $\tilde E \in \C(z,t)_{\tilde U} \otimes V$ by the condition
\[
\langle \tilde E, 1\otimes (H_{j_1},\dots,H_{j_k})\rangle =  f_{j_1,\dots,j_k} ,
\]
for any independent $\{j_1,\dots,j_k\}\subset J$.

\begin{thm}
\label{thm loc alg and Hams}
For any $j\in J$, there exist elements $h_1,\dots,h_k \in \C(z,t)_{\tilde U}\otimes V$ such that
\begin{gather*}
\left(1\otimes K_j(z)\right)  \tilde E(z,t)  = \left(\frac {a_j}{f_j(z,t)} \otimes 1\right)
 \tilde E(z,t) +
\sum_{i=1}^k\left(\frac{\der \Phi}{\der t_i}(z,t)\otimes 1\right)  h_i(z,t).
\end{gather*}
\end{thm}

\begin{proof}
We have
\begin{gather}
\label{om(a)}
\nu(a) = \sum_{i=1}^k \frac{\partial \Phi} {\partial t_i } (z,t)\, dt_i
+\sum_{j\in J} \frac {a_j}{f_j(z,t)}dz_j
\end{gather}
and
\begin{gather}
\label{vz}
\sum_{{\rm independent } \atop \{j_1 < \dots < j_k\} \subset J }\!\!\!
\omega_{j_1} \wedge \dots \wedge \omega_{j_k}
\otimes
F(H_{j_1}, \dots , H_{j_k}) =  \tilde E(z,t) \left(dt_1\wedge\dots\wedge dt_k \otimes 1\right) + z\text{-part} ,\!
\end{gather}
where $z$-part is a $V$-valued dif\/ferential $k$-form
with zero restriction to each f\/iber of the projection
$\C^n\times \C^k \to \C^n$. Then identity \eqref{formula H} and
formulas \eqref{om(a)}, \eqref{vz}, \eqref{LK}
imply the theorem.
\end{proof}

\subsection{Hamiltonians, critical points and the canonical element}
\label{sec Hams and cr pts}

Fix  $z\in \C^n-\Delta$.
Recall that in Section~\ref{Construction} we have def\/ined
 the quantum integrable model assigned to
$(\A(z), a)$ to be the collection
\[
\big(\sing W;\; S^{(a)}|_{\sing W};\; K_1(z)^*|_{\sing W},\dots, K_n(z)^*|_{\sing W} : \  \sing W\to \sing W\big).
\]

Let $p \in U(\A(z))$ be an isolated critical point of the master function $\Phi(z,\cdot)
: U(\A(z)) \to \C$.
Let $A_{p,\Phi}$ be the local algebra of the critical
point and  $[\ ] : \C(t)_{U(\A(z))}\to A_{p,\Phi}$ the canonical projection.
 Denote by $\Hess $ the  Hessian of $\Phi(z,\cdot)$ with respect to variables $t_1,\dots,t_k$.

Let $E \in \C(t)_{U(\A(z))}\otimes V$ be the canonical element associated with $\A(z)$,
see Section~\ref{Special}.
Denote by $[E]$ its  projection
to $A_{p,\Phi}\otimes V$.
By Lemma~\ref{lem v sing} we have $[E]   \in A_{p,\Phi}\otimes \Sing V$.
Let $ V \to W$ be the  map associate with the contravariant form and
$[\bs E]$ the image of $[E]$ under the induced map
$A_{p,\Phi}\otimes \Sing V \to A_{p,\Phi}\otimes \sing W$,
\begin{gather*}
[\bs E]   \in  A_{p,\Phi}\otimes \sing W .
\end{gather*}

\begin{thm}
\label{thm v good fiber}
We have
\begin{enumerate}\itemsep=0pt
\item[$(i)$]
$S^{(a)}([\bs E],[\bs E]) = (-1)^k[\Hess]$,

\item[$(ii)$]
$(1\otimes K_j(z)^*)[\bs E] =
([a_j/f_j(z,\cdot)]\otimes 1)[\bs E]$  for $j\in J$.
\end{enumerate}
\end{thm}

\begin{proof}
Part $(i)$ follows from Theorem \ref{first theorem}.
Part $(ii)$ follows from Theorem \ref{thm loc alg and Hams}.
\end{proof}

\begin{rem}\label{rem a_j/f_j form basis]}\rm
The elements
$[a_j/f_j(z,\cdot)]$,  $j\in J$, generate $A_{p,\Phi}$ due to
 Lemma~\ref{lem f_j generate} and the assumption   $(a_j\neq 0$ for all $j)$.
\end{rem}

\section{Asymptotic solutions
and eigenvectors of Hamiltonians}
\label{sec asympt}

\subsection{Asymptotic solutions, one variable}
Let $u$ be a variable, $W$ a vector space,
$M(u) \in \End(W)$ an endomorphism  depending holomorphically on $u$ at $u=0$.
Consider a dif\/ferential equation,
\begin{gather}
\label{1 eqn}
\kappa \frac{dI}{du}(u) = M(u)I(u)
\end{gather}
depending on a complex parameter $\kappa\in\C^*$.

Let $P(u)\in\C$, $(w_m(u)\in W)_{m\in \Z_{\geq 0}}$ be functions holomorphic
at $u=0$ and $w_0(0)\neq 0$. The series
\begin{gather}
\label{1 assymp}
I(u,\kappa)  =  e^{P(u)/\kappa} \sum_{m=0}^\infty  w_m(u) \kappa^m
\end{gather}
will be called
 an asymptotic solution to \eqref{1 eqn} if it satisf\/ies \eqref{1 eqn} to all orders in $\kappa$.
In particular, the leading order equation is
\begin{gather}
\label{lead eqn reg}
\frac {dP}{du}(u) w_0(u) = M(u) w_0(u) .
\end{gather}

Assume now that
\[
M(u) = \frac{M_{-1}}u + M_0 + M_1u + \cdots ,
\qquad
M_j \in \End(W),
\]
has a f\/irst order pole at $u=0$ and $I(u,\kappa)$ is
a series like in~\eqref{1 assymp}. The series $I(u,\kappa)$ will be called an asymptotic solution
to equation~\eqref{1 eqn} with such $M(u)$ if it satisf\/ies~\eqref{1 eqn} to all orders in
$\kappa$. In particular, the leading order equation is again equation~\eqref{lead eqn reg}.
Equation \eqref{lead eqn reg} implies
\begin{gather}
\label{lead eqn irreg}
w_0(0) \in {\rm ker}\, M_{-1},
\qquad
\frac{dP}{du}(0) w_0(0) = M_0w_0(0) + M_{-1} \frac{dw_0}{du}(0) .
\end{gather}

\subsection{Critical points of the master function and asymptotic solutions}
\label{sec crit and asympt}

Let us return to the situation of Section \ref{sec par trans}.

Let $t(z)$
 be a nondegenerate critical point of
$\Phi(z,\,\cdot\,) : U(\A(z))\to\C$. Assume that $t(z)$
depends on $z$ holomorphically  in a neighborhood of a point
$z^0\in \Bbb C^n$.
Fix a univalued branch of $\Phi$ in a neighborhood of $(z^0,t(z^0))$
(by choosing arguments of all of the logarithms).
Denote $\Psi(z) = \Phi(z,t(z))$.

Let $B\subset\C^k$ be a small ball with center at $t(z^0)$.
Denote
\[
B^- = \big\{ t\in B \ |\ {\rm{Re }}\,  \Phi\big(z^0,t\big(z^0\big)\big)
 > {\rm{Re }} \  \Phi\big(z^0,t\big)\big\}   .
\]
It is well known that $H_k(B,B^-;\Bbb Z)=\Bbb Z$, see
for example~\cite{AGV}.
There exist local coordinates $u_1,\dots ,u_k$ on $\Bbb C^k$ centered
at $t(z^0)$ such that
$\Phi(z^0,u) = -u^2_1 - \dots -u^2_k + {\text {const}}$.
Denote
\[
\delta  = \big\{ (u_1,\dots ,u_k)\in\Bbb R^k \ |\  u^2_1 +\dots +u^2_k  \leq   \epsilon\big\}  ,
\]
where $\epsilon$ is a small positive number. That $\delta$, considered as a $k$-chain,
generates
 $H_k(B,B^-;\Bbb Z)$.

Def\/ine an element $\{\delta\}(z,\kappa) \in V$ by the formula
\begin{gather}
\label{delta}
\{\delta\}(z,\kappa)   :\ V^*  \to  \C  ,
\qquad
\omega  \mapsto  \kappa^{-k/2}
\int_{\delta} e^{\Phi(z,t)/\kappa}  \om(z,t)  .
\end{gather}
Recall that any
element $\omega \in V^*$ is a linear combination of elements
$(H_{j_1},\dots,H_{j_k})$
 and such an element $(H_{j_1},\dots,H_{j_k})$
is identif\/ied with the dif\/ferential form
\[
\omega_{j_1} \wedge\dots\wedge\omega_{j_k} =
df_{j_1}(z,t)/f_{j_1}(z,t) \wedge\dots\wedge df_{j_k}(z,t)/f_{j_k}(z,t) .
\]
In \eqref{delta} we integrate over $\delta$ such a dif\/ferential form multiplied by
$e^{\Phi(z,t)/\kappa}$.

The element $\{\delta\}(z,\kappa)$ as a function of $z$, $\kappa$ is holomorphic if
$z$ is close to $z^0$ and $\kappa \neq 0$.

\pagebreak

 \begin{thm}  \quad
\label{thm asy}
\begin{enumerate}\itemsep=0pt
\item[$(i)$]
 Let $\kappa \in \Bbb R$ and $\kappa \to +0$. Then
the function $\{\delta\}(z,\kappa)$ has an
asymptotic expansion
\begin{gather}
\label{as expansion}
\{\delta\}(z,\kappa)  =   e^{\Psi(z)/ \kappa}  \sum^\infty_{m=0}
w_m(z) \kappa^m  ,
\end{gather}
where $(w_m(z)\in V)_{m\in \Z_{\geq 0}}$ are functions of $z$ holomorphic at $z^0$ and
\begin{gather*}
w_0(z)  =   \pm(2\pi)^{k/2}
\left((-1)^k \det_{1\leq i,j \leq k}
\left( \frac{\partial^2\Phi} {\partial t_i \partial t_j} \right)  (z,t(z)) \right)^{-1/2}
v(z,t(z))  .
\end{gather*}
Here $v(z,t(z))$ is the special vector associated
 with the critical point $(z,t(z))$ of the function
$\Phi(z,\,\cdot\,)$, see Section~{\rm \ref{Special}}. The sign $\pm$ depends on the choice of
the orientation of $\delta$.

\item[$(ii)$]
 The asymptotic expansion~\eqref{as expansion}
gives an asymptotic solution to the
Gauss--Manin diffe\-rential equations~\eqref{GM}.

\item[$(iii)$]
 The functions $(w_m(z))_{m\in \Z_{\geq 0}}$ take values in $\sing V$.
\end{enumerate}
\end{thm}

Part $(i)$ of the theorem is a direct corollary of the method of steepest
descent; see, for example, \S~11 in~\cite{AGV}.  Part~$(ii)$
follows from Lemma~\ref{lem 2} and formula of dif\/ferentiation of an integral.
Part~$(iii)$ follows from Stokes' theorem.

\begin{rem}\rm
The def\/inition of $\delta$ depends on the choice of local coordinates
$u_1,\dots,u_k$, but the asymptotic expansion~\eqref{as expansion} does not depend on the choice
of $\delta$ since the dif\/ference of the corresponding integrals is exponentially small.
\end{rem}

\begin{cor}
\label{cor bethe eigen}
Let $z^0$, $t=t(z)$, $\Psi(z)$ be the same  as in Theorem {\rm \ref{thm asy}}.
Assume that $z^0 \in \C^n-\Delta$.
In that case the operators $K_j(z^0)|_{\sing V} : \s V\to \s V$,
$j\in J$, are all well-defined, see~\eqref{K_j}, and
we have
\begin{gather*}
K_j\big(z^0\big)  v\big(z^0,t\big(z^0\big)\big)  =
\frac {\der \Psi}{dz_j}\big(z^0\big)  v\big(z^0,t\big(z^0\big)\big) ,
\qquad j\in J  ,
\end{gather*}
Thus, the special vector $v(z^0,t(z^0))$ is an eigenvector of the geometric Hamiltonians
$K_j(z^0)$.
\end{cor}

\begin{proof}
The corollary follows from equation \eqref{lead eqn reg}.
\end{proof}

Note that $\frac {\der \Psi}{dz_j}(z^0) = \frac{a_j}{f_j(z^0,t(z^0))}$.

\begin{rem}
\label{rem still is an asympt solution}
\rm The Gauss--Manin dif\/ferential equations \eqref{GM} have singularities over the discriminant
$\Delta\subset \C^n$. If  $z^0\in\Delta$, then expansion~\eqref{as expansion}
still gives an asymptotic solution to equations~\eqref{GM} and
that asymptotic solution is regular at~$z^0$.
\end{rem}

\section{Hamiltonians of bad f\/ibers}
\label{sec real posit bad}

\subsection{Naive geometric Hamiltonians}
\label{sec general remarks}
Let us return to the situation of Section~\ref{sec par trans}.
Let $z^0 \in \Delta$ and $z\in\C^n-\Delta$.
We have{\samepage
\begin{gather}
\Sing \FF^k(\A(z^0))
\subset \Sing \FF^k(\A(z)) \subset  \FF^k(\A(z)) ,\notag
\\
\Sing \FF^k(\A(z^0))
\subset
\FF^k(\A(z^0))
\subset \FF^k(\A(z)) ,
\label{subspaces}
\\
\Sing \FF^k(\A(z^0)) = \FF^k(\A(z^0)) \cap (\Sing \FF^k(\A(z))) ,
\notag
\end{gather}
 see Section \ref{sec Bad fibers}. Recall that $\FF^k(\A(z))$ was denoted by $V$.}

Consider the map $V\to V^*$ corresponding to the contravariant form.
In Section~\ref{Construction}
we denoted the images of $V$ and $\sing V$
by $W$ and $\Sing W$, respectively.
We denote the images of $\FF^k(\A(z^0))$
and $\sing \FF^k(\A(z^0))$ by $W(z^0)$ and $\sing W(z^0)$, respectively.
We have
\[
\Sing W\big(z^0\big)
\subset \Sing W \subset  W ,
\qquad
\Sing W\big(z^0\big)
\subset
W\big(z^0\big) \subset
W .
\]

Recall that $\Delta$ is the union of hyperplanes $H_C$, $C\in \frak C$. Denote
\begin{gather*}
\frak C_0 = \{ C\in \frak C\ | \ z^0\in H_C\} .
\end{gather*}
Consider the operator-valued functions
$K_j(z) : V \to V$, $j\in J$, given by formula~\eqref{K_j}. Denote
\begin{gather*}
K_j^0(z)   =   \sum_{C\in \frak C_0}
 \frac{\la_j^C}{f_C(z)}   L_C  ,
\qquad
K_j^1(z) = K_j(z) - K_j^0(z)  .
\end{gather*}
Each of the summands of $K_j^0(z)$ tends to inf\/inity as $z$ tends to $z^0$ in
$\C^n-\Delta$. The operator-valued function
$K_j^1(z)$ is regular at $z^0$.

The operators $K_j(z)^*$, $L^*_C$ preserve the subspaces $\sing W\subset W\subset V^*$ and are symmetric operators on~$W$
 with respect to the contravariant form  on~$W$.
The operators $K_j(z)^*$ restricted to $\Sing W$ commute.

The point $z^0\in\Delta$ def\/ines an edge $X_{z^0}$  of the arrangement $(H_C)_{C\in \frak C}$, where
$X_{z^0} = \cap_{C\in\frak C_0} H_C$.
Denote by $T_{z^0}$ the vector space of constant vectors f\/ields on $\C^n$
 which are tangent to
$X_{z^0}$,
\begin{gather*}
T_{z^0} = \left\{
 \xi = \sum_{j\in J}\xi_j \frac{\der}{\der z_j}\ | \ \xi_j\in\C, \
\xi (f_C) = 0 \ \text{for all}\ C\in \frak C_0\right\} .
\end{gather*}

\begin{Lemma}
\label{lem regularity of KZ}
For any $\xi\in T_{z^0}$,
\begin{enumerate}\itemsep=0pt
\item[$(i)$]

The linear operator
\begin{gather*}
 K_\xi(z) = \sum_{j\in J} \xi_j K_j(z)   :\  V   \to    V ,
\end{gather*}
considered as a function of $z$, is regular at $z^0$, moreover,
\[
K_\xi(z) = \sum_{j\in J} \xi_j K_j^1(z) .
\]
\item[$(ii)$]
The linear operator $K_\xi(z^0)$ preserves the subspace
$\FF^k(\A(z^0))\subset V$.
\item[$(iii)$]

The dual linear operator
\[
 K_\xi(z)^*   : \ V^*   \to   V^* ,
\]
considered as a function of $z$, is regular at $z^0$, moreover,
\[
K_\xi(z)^* = \sum_{j\in J} \xi_j K_j^1(z)^* .
\]
\item[$(iv)$]
The linear operator $K_\xi(z^0)^*$ preserves the subspace
$\Sing W(z^0)\subset V^*$.
\end{enumerate}
\end{Lemma}

\begin{proof}
Parts $(iii)$, $(iv)$ follow from parts $(i)$, $(ii)$. Part $(i)$ is clear. Part $(ii)$ follows from a~straightforward calculation.
\end{proof}

The operators $K_\xi(z^0)^*$ preserve the subspace $\sing W(z^0)$. The operators
\[
K_\xi(z^0)^*|_{\sing W(z^0)}   : \ \sing W(z^0)   \to
\sing W(z^0),
\qquad
\xi \in T_{z^0},
\]
form a commutative family of linear operators. The operators are
symmetric with respect to the contravariant form. These operators will be called {\it naive geometric Hamiltonians} on
$\sing W(z^0)$.

\subsection[Space $\FF^k(\A(z^0))$ and operators $L_C$]{Space $\boldsymbol{\FF^k(\A(z^0))}$ and operators $\boldsymbol{L_C}$}
\label{Space FF and operators L-C}

\begin{Lemma}
\label{lem F(z) and ker L-C}
\qquad
\begin{enumerate}\itemsep=0pt
\item[$(i)$]
The space $\FF^k(\A(z^0))$ lies in the kernel of $L_C : V\to V$ for any $C \in \frak C_0$.
\item[$(ii)$]

The space $W(z^0)$ lies in the kernel of $L_C^*|_W : W\to W$ for any $C \in \frak C_0$.

\item[$(iii)$]
For any $C\in\frak C_0$, the image of $L_C^*|_W$ is orthogonal to
$W(z^0)$ with respect to the contravariant form.
\end{enumerate}
\end{Lemma}

\begin{proof} Part $(i)$ follows from a straightforward easy calculation. Part $(ii)$ follows from
 part~$(i)$. Part~$(iii)$ follows from part~$(ii)$ and the fact that~$L^*_C$ is symmetric.
 \end{proof}

\subsection{Conjecture}
\label{Main Conjecture}

\begin{conj}
\label{main conj} Let $z^0\in\Delta$.
Assume that the contravariant form restricted to $\Sing W(z^0)$ is nondegenerate. Let  ${\rm pr} :
\Sing W\to \Sing W(z^0)$ be the orthogonal projection with respect to the contravariant form.
Then
the linear operators
\begin{gather*}
{\rm pr} K_j^1\big(z^0\big)^*|_{\Sing W(z^0)} : \ \Sing W\big(z^0\big)\to \Sing W\big(z^0\big)  ,
\qquad
j\in J ,
\end{gather*}
commute and are symmetric with respect to the contravariant form.
\end{conj}

For  $z^0\in \Delta$,
we def\/ine {\it the quantum integrable model assigned to}
$(\A(z^0), a)$ to be the collection
\begin{gather*}
\big(\Sing W\big(z^0\big),\ S^{(a)}|_{\Sing W(z^0)},
\nonumber
\\
\qquad
{\rm pr} K_j^1\big(z^0\big)^*|_{\Sing W(z^0)} : \ \Sing W\big(z^0\big)\to \Sing W\big(z^0\big),\
\text{ where}\
j\in J \big) .
\end{gather*}
The unital subalgebra of $\End(\Sing W(z^0))$ generated by operators
\[
{\rm pr}K_j^1\big(z^0\big)^*|_{\Sing W(z^0)} ,
\qquad
j\in J,
\]
 will be called {\it the
algebra of geometric Hamiltonians of $(\A(z^0),a)$}.

Note that the naive geometric Hamiltonians are elements of the algebra of geometric Hamiltonians, since
for any
$
\xi = \sum_{j\in J} \xi_j \frac{\der}{\der z_j} \in T_{X_{z^0}}$, we have
\[
K_\xi\big(z^0\big)^*|_{\Sing W(z^0)} = \sum_{j\in J}
\xi_j{\rm pr} K_j^1\big(z^0\big)^*|_{\Sing W(z^0)} .
\]

 In the next section we prove the conjecture under Assumption~\ref{i} of
 certain positivity and reality conditions, see Theorem~\ref{thm bar K_j}.
In Section \ref{More on Hamiltonians of bad fibers} more results in this direction will be obtained, see Theorems~\ref{thm BAD} and~\ref{thm ham of unb arr}.
For applications to the Gaudin model an equivariant version of the conjecture is needed, see
Sections~\ref{Arrangements with symmetries sec} and~\ref{sec Gaudin}.

\subsection[Positive  $(a_j)_{j\in J}$, real $(g_j)_{j\in J}$]{Positive  $\boldsymbol{(a_j)_{j\in J}}$, real $\boldsymbol{(g_j)_{j\in J}}$}
\label{sec  z^0 in Delta cap R^n }

\begin{assumption} \label{i}
 Assume  that all weights $a_j,$ $j\in J$, are positive and
all functions
$g_j = b_j^1t_1+\cdots + b_j^kt_k$, $ j\in J$,
have real coefficients $b_j^i$.
\end{assumption}

The space $V$ has a real structure,
$V=V_\R\otimes_\R\C$,  see Section \ref{sec real structure}.
Under Assumption \ref{i} all subspaces in \eqref{subspaces} are
real (can be def\/ined by real equations).
 The contravariant form
$S^{(a)}$
is positive def\/inite on $V_\R$ and is positive def\/inite on the real parts of all
of the
subspaces in \eqref{subspaces}.

Denote by ${\rm pr} : \Sing V \to \Sing \FF^k(\A(z^0))$ the orthogonal projection.

\begin{thm}
\label{thm bar K_j}
Assume that Assumption {\rm \ref{i}} is satisfied and  $z^0\in \Delta$.
Then the operators
\begin{gather*}
{\rm pr} K_j^1\big(z^0\big)|_{\Sing \FF^k(\A(z^0))}
: \
\Sing \FF^k\big(\A\big(z^0\big)\big) \to  \Sing \FF^k\big(\A\big(z^0\big)\big),
\qquad
j\in J ,
\end{gather*}
 commute and are symmetric with respect to the contravariant form.
\end{thm}

Theorem \ref{thm bar K_j} proves Conjecture \ref{main conj} under Assumption \ref{i}.

\subsection[Proof of Theorem \ref{thm bar K_j} for $z^0\in \Delta \cap \R^n$]{Proof of Theorem \ref{thm bar K_j} for $\boldsymbol{z^0\in \Delta \cap \R^n}$}
\label{Proof of Theorem {thm K_j} z real}

Assume that $z^0\in \Delta \cap \R^n$.
Let $r: (\C,\R,0)\to (\C^n,\R^n,z^0)$ be a germ of a holomorphic curve such that
$r(u)\in \C^n-\Delta$ for $u\neq 0$.
For  $u\in \R_{>0}$, the arrangement $\A(r(u))$ is real.
Denote  $U(r(u))_\R=(\C^k - \cup_{j\in J}H_j(r(u)))\cap \R^k$.
Let $U(r(u))_\R = \cup_{\al} D_\al(r(u))$ be the decomposition into the union
of connected components.
We label components so that for any
$\al$, the component $D_\al(r(u))$ continuously depends on $u>0$.
Let $A$ be the set of all $\al$ such that $D_\al(r(u))$ is bounded.
Let $A_1$ be the set of all $\al$ such that $D_\al(r(u))$
 is bounded and vanishes as $u\to +0$
(the limit of $D_\al(r(u))$ is not a domain of $\A(z^0)$). Let
$A_2$ be the set of all $\al$ such that $D_\al(r(u))$
is bounded and the limit of $D_\al(r(u))$ as $u\to 0$
is a domain of $\A(z^0)$. We have $A=A_1\cup A_2$ and $A_1\cap A_2=\varnothing$.

All critical points
of $\Phi(r(u),\,\cdot\,)$ lie
in $\cup_{\al\in A}D_\al(r(u))$. Each domain
$D_\al(r(u))$ contains a unique critical point
$(r(u),t(u)_\al)$ and that  critical point is nondegenerate. Denote by
$v(r(u),t(u)_\al) \in \Sing V$
the corresponding special
vector. That vector is an eigenvector of the geometric Hamiltonians,
\[
K_j(r(u))\,v(r(u),t(u)_\al)
=  \frac {a_j}{f_j(r(u),t(u)_\al)}
v(r(u),t(u)_\al) ,
\qquad
j\in J .
\]
If $\al \in A_2$, then all eigenvalues
$1/{f_j(r(u),t(u)_\al)}$, $j\in J$, are regular functions at $u=0$.
If $\al \in A_1$, then
there is an index $j\in J$ such that
$a_j/{f_j(r(u),t(u)_\al)}\to \infty$ as $u\to 0$.

\begin{Lemma}\label{lem lim}\sloppy
The span $\langle v(r(u),t(u)_\al) \rangle_{\al \in A_2}$ has
a limit as $u\to +0$.
That limit is
$\Sing \FF^k(z^0)\subset \Sing V$. Similarly,
the span $\langle v(r(u),t(u)_\al) \rangle_{\al \in A_1}$ has a
 limit as $u\to +0$.
That limit is
$(\Sing \FF^k(z^0))^\perp\subset \Sing V$ where ${}^\perp$ denotes the orthogonal complement.
\end{Lemma}

\begin{proof}
The lemma follows from Theorem \ref{thm real crit pts} and Corollary
\ref{cor real basis}.
\end{proof}

Assume now that a curve $r(u)=(z_1(u),\dots,z_n(u))$
is linear in  $u$. Then for any $j$ we have
 $K_j^0(r(u)) = N_j/u$ where $N_j :\Sing V
\to \Sing V$ is an operator independent of $u$.

\begin{Lemma}
\label{lem N_j}
The image of $N_j$ is a subspace of
$(\Sing \FF^k(z^0))^\perp$.
\end{Lemma}

\begin{proof}
The lemma follows from Theorem \ref{thm real crit pts} and Corollary
\ref{cor real basis}.
\end{proof}

By formula \eqref{lead eqn irreg} and Lemma \ref{lem N_j}, for any $\al \in A_2$ we have
\[
\frac {a_j}{f_j(r(0),t(0)_\al)} v(r(0),t(0)_\al) =
K_j^1(r(0)) v(r(0),t(0)_\al) + v_1 ,
\]
where $v_1 \in (\Sing \FF^k(z^0))^\perp$.
Thus,
\begin{gather}
\label{*}
{\rm pr} K_j^1(r(0))
\,v(r(0),t(0)_\al)   =
\frac 1{f_j(r(0),t(0)_\al)}  v(r(0),t(0)_\al) .
\end{gather}
Thus, all operators ${\rm pr} K_j^1(r(0))$ are diagonal in the basis
$(v(r(0),t(0)_\al))_{\al\in A_2}$
of
  $\Sing \FF^k(\A(z^0))$.
Equation \eqref{*} f\/inishes the proof of the commutativity of
operators ${\rm pr} K_j^1(z^0)$.

The symmetry of the operators ${\rm pr} K_j^1(z^0)$
with respect to the contravariant form
follows from the fact that operators ${\rm pr} K_j^1(z^0)$ are diagonal in the orthogonal basis
$(v(r(0),t(0)_\al))_{\al\in A_2}$.

\subsection[Proof of Theorem \ref{thm bar K_j} for any $z^0\in \Delta$]{Proof of Theorem \ref{thm bar K_j} for any $\boldsymbol{z^0\in \Delta}$}
\label{Proof of thm K_j}

Consider the arrangement $(H_C)_{C\in\frak C}$ in $\C^n$ and its arbitrary
edge $X$.
The arrangement $(H_C)_{C\in\frak C}$ is real, see Remark~\ref{rem g real}.
The edge $X$ is the complexif\/ication of $X\cap\R^n$.

Denote
\begin{gather}
\label{hat L}
\hat X = X - \cup_{C\in \frak C-\frak C_X}X\cap H_C .
\end{gather}
For any $z^1,z^2\in\hat X$, the subspaces
\[
\Sing \FF^k\big(\A\big(z^1\big)\big) \subset \Sing V,
\qquad
 \Sing \FF^k\big(\A\big(z^2\big)\big) \subset \Sing V
\]
coincide. Denote that subspace by
$\Sing \FF^k(\A(X))
\subset \Sing V$.

For $z\in \hat X$, the  operators ${\rm pr}K^1_j(z)  : \Sing \FF^k(\A(X))
 \to \Sing \FF^k(\A(X))$,
$j\in J$, depend on $z$ holomorphically. The operators commute and are symmetric
for $z\in \hat X \cap \R^n$, by reasonings in
Section~\ref{Proof of Theorem {thm K_j} z real}.
Hence they commute and are symmetric for all $z\in \hat X$.

\subsection{Critical points and eigenvectors}

\begin{thm}
\label{thm eigen posit}
Assume Assumption {\rm \ref{i}}.
Let $z^0\in \Delta$ and let
 $p$ be a critical point of
$\Phi(z^0,\,\cdot\,): U(\A(z^0))\to\C$.
Then the corresponding special vector $v(z^0,p) \in \Sing \FF^k(z^0)$ $($if nonzero$)$
is an eigenvector of the operators ${\rm pr} K_j^1(z^0)$, $j\in J$,
\begin{gather}
\label{eig pos}
{\rm pr} K_j^1\big(z^0\big)  v\big(z^0,p\big)  =
\frac {a_j}{f_j(z^0,p)}  v\big(z^0,p\big).
\end{gather}
\end{thm}

\begin{proof} If $z^0\in \Delta\cap \R^n$, then
the  theorem is just a restatement of formula~\eqref{*}.

Assume that $z^0$ is an arbitrary point of $\Delta$. Then there exists an
edge $X$ of the arrangement
 $(H_C)_{C\in\frak C}$ such that $z^0\in \hat X$, see~\eqref{hat L}.
For $z^0\in \hat X$, all objects in formula~\eqref{eig pos} depend on $z^0$ algebraically.
Hence, the fact, that formula~\eqref{eig pos} holds for all critical points if
$z^0\in \hat X\cap \R^n$, implies Theorem~\ref{thm eigen posit} for any
$z^0\in \hat X$.
\end{proof}

\subsection{Hamiltonians, critical points and the canonical element}
\label{sec Hams and cr pts, bad fibers}

Let Assumption \ref{i} be satisf\/ied.
Fix  $z^0\in \Delta$.
We have
def\/ined
 the quantum integrable model assigned to
$(\A(z^0), a)$ to be the collection
\begin{gather*}
\big(\Sing \FF^k\big(\A\big(z^0\big)\big);\ S^{(a)}|_{\Sing \FF^k(\A(z^0))};
\\
\qquad {\rm pr} K_j^1\big(z^0\big)|_{\Sing \FF^k(\A(z^0))} : \
\Sing \FF^k\big(\A\big(z^0\big)\big) \to \Sing \FF^k\big(\A\big(z^0\big)\big),\ \text{where}\ j\in J \big) ,
\end{gather*}
see Section \ref{sec  z^0 in Delta cap R^n }.

Let $p \in U(\A(z^0))$ be an isolated critical point of the function $\Phi(z^0,\,\cdot\,)
: U(\A(z^0)) \to \C$.
Let $A_{p,\Phi}$ be the local algebra of the critical
point  and  $[ \ ] : \C(t)_{U(\A(z^0))}\to A_{p,\Phi}$ the canonical projection.
 Denote by $\Hess $ the  Hessian of $\Phi(z^0,\,\cdot\,)$.

Let $E \in \C(t)_{U(\A(z^0))} \otimes \FF^k(\A(z^0))$ be the canonical
element associated with $\A(z^0)$,
see Section~\ref{Special}. Denote by $[ E]$ the projection of the canonical
element to $A_{p,\Phi}\otimes \FF^k(\A(z^0))$.
By Lemma~\ref{lem v sing}, we have
\[
[ E] \in    A_{p,\Phi}\otimes \Sing \FF^k\big(\A\big(z^0\big)\big) .
\]

\begin{thm}
\label{thm v bad fiber}
We have
\begin{enumerate}\itemsep=0pt
\item[$(i)$]
$S^{(a)}([E],[E]) = (-1)^k[\Hess]$,

\item[$(ii)$]
$(1\otimes {\rm pr}K^1_j(z^0))[ E] =
([a_j/f_j(z,\cdot)]\otimes 1)[ E]$ for $j\in J$.
\end{enumerate}
\end{thm}

\begin{proof}
Part $(i)$ follows from Theorem~\ref{first theorem}.
Part $(ii)$ follows from Theorem~\ref{thm eigen posit}.
\end{proof}

\section{Geometric interpretation of the algebra of Hamiltonians}
\label{Geometric interpretation of the algebra of Hamiltonians}

\subsection{An abstract setting}
\label{An abstract setting}
Let $k < n$ be positive integers and $J=\{1,\dots,n\}$ as before.

Let $F$ be a germ of a holomorphic function at a  point $p\in\C^k$. Assume that $p$
is an isolated critical point of $F$
with Milnor number $\mu_p$.
Let $A_{p,F}$ be the local algebra of the critical
point and $(\,,\,)_p$ the residue
bilinear form on $A_{p,F}$, see \eqref{Gr form}.
 Denote by
$[\HessF]$ the projection to $A_{p,F}$ of the germ
$\det_{1\leq l,m\leq k}(\der^2F/\der t_l\der t_m)$.

Let $h_1,\dots,h_{\mu_p}$ be a $\C$-basis of $A_{p,F}$.
Let $g_1,\dots,g_n \in A_{p,F}$ be a collection of elements such that
the unital subalgebra of $A_{p,F}$ generated by $g_1,\dots,g_n$
equals $A_{p,F}$.

Let $W$ be a vector space with a
symmetric bilinear form $S$. Let
$M_j : W\to W$, $j\in J$, be
a~collection of commuting symmetric linear operators,
\begin{gather*}
M_iM_j=M_jM_i,\qquad
 S(M_ju,v) = S(u,M_jv)
\qquad
{\rm for\ all}\ i,j\in J \ {\rm and}\ u,v\in W.
\end{gather*}

Assume that an element
\begin{gather*}
w = \sum_{l=1}^{\mu_p} h_l\otimes w_l   \in   A_{p,F} \otimes W
\end{gather*}
is given such that
\begin{gather}
\label{ass 2}
\sum_{l=1}^{\mu_p} h_l\otimes M_jw_l  =   \sum_{l=1}^m g_jh_l\otimes w_l ,
\qquad j\in J ,
\\
\label{ass 1}
 \sum_{l,m=1}^{\mu_p} S(w_l,w_m) h_l h_m   =  (-1)^k[\HessF] .
\end{gather}

Denote by $Y\subset W$ the vector subspace generated by $w_1,\dots,w_{\mu_p}$.
By property \eqref{ass 2}, every~$M_j$, $j\in J$, preserves~$Y$.
Denote by $A_Y$ the unital subalgebra of~$\End(Y)$ generated by~$M_j\vert_Y$, $j\in J$.
 The subspace $Y$ is an $A_Y$-module.
Def\/ine a linear map
\begin{gather*}
\al : \ A_{p,F} \to Y,
\qquad
f \mapsto \sum_{l=1}^{\mu_p} (f,h_l)_p w_l .
\end{gather*}

\begin{thm}
\label{model thm}
\quad
\begin{enumerate}\itemsep=0pt
\item[$(i)$]
The map $\al : A_{p,F} \to Y$ is an isomorphism of vector spaces.
The form $S$ restricted to $Y$ is nondegenerate.
\item[$(ii)$]
 The map $g_j \mapsto M_j\vert_{Y}$, $j\in J$,
extends uniquely to an algebra isomorphism
$\beta : A_{p,F} \to A_Y$.
\item[$(iii)$] The maps $\al$, $\beta$ give an isomorphism of the regular representation
of $A_{p,F}$ and the $A_Y$-module $Y$, that is $M_j\al(f) = \al(g_jf)$
for any $f\in A_{p,F}$ and $j\in J$.
\item[$(iv)$]
Define the value $w(p)$ of $w$ at $p$ as the image of $w$ under the natural projection
$A_{p,F} \otimes W \to  A_{p,F}/\frak m_p \otimes W = W$. Then $w(p) = \al(\HessF)/\mu_p$ and
the
value $w(p)$ is nonzero. The vector $w(p)$ is the only $($up to proportionality$)$
common eigenvector of the operators $M_j\vert_Y$, $j\in J$, and we have
$M_jw(p) = g_j(p)w(p)$.
\end{enumerate}
\end{thm}

This theorem is an analog of Theorem~5.5 and Corollary~5.6 in~\cite{MTV5}. The proof is analogous to the  proofs in~\cite{MTV5}.

\subsection{Proof of Theorem \ref{model thm}}
\label{ Proof of Theorem model thm}

\begin{Lemma}
\label{hom lem}
We have $M_j\al(f) = \al(g_jf)$ for any $f\in A_{p,F}$, $j\in J$.
\end{Lemma}

\begin{proof}
We have $M_j\al(f) = \sum\limits_{l=1}^{\mu_p} (f,h_l)_p M_jw_l$. By~\eqref{ass 2},
that is equal to
 $\sum\limits_{l=1}^m (f,g_jh_l)_p w_l =
\sum\limits_{l=1}^m (g_jf, h_l)_p w_l = \al(g_jf)$.
\end{proof}

Def\/ine a bilinear form $(\,,\,)_S$ on $A_{p,F}$,
\begin{gather*}
(f,g)_S =
S(\al(f),\al(g)) =
\sum_{l,m=1}^{\mu_p} S(w_l,w_m) (f,h_l)_p(g,h_m)_p .
\end{gather*}
\begin{Lemma}
\label{lem ()S} We have $(fg,h)_S = (f,gh)_S$ for all
$f,g,h\in A_{p,F}$.
\end{Lemma}

\begin{proof}
Since $g_j$, $j\in J$, generate $A_{p,F}$ it is enough to show that
$(fg_j,h)_S =  (f,g_jh)_S$ for all
$f,h\in A_{p,F}$, $j\in J$.
We have $(fg_j,h)_S =
S(\al(fg_j),\al(h))
= S(M_j\al(f), \al(g)) = S(\al(f),M_j\al(g)) $ $=
S(\al(f),\al(g_jh))= (f,g_jh)_S.$
\end{proof}

\begin{Lemma}\label{lem F el}
There exists a unique element
$s\in A_{p,F}$ such that $(f,g)_S = (sf,g)_p$ for all $f,g\in A_{p,F}$.
\end{Lemma}

\begin{proof}
Consider the linear function $A_{p,F}\to \C$, $f\mapsto (1,f)_S$. Since the bilinear
form $(\,,\,,)_p$ is nondegenerate there exists a unique
$s\in A_{p,F}$ such that
$(1,f)_S=(s,f)_p$ for any $f\in A_{p,F}$. Hence for any $f,g\in A_{p,F}$ we have
$(f,g)_S=(1,fg)_S=(s,fg)_p= (sf,g)_p$.
\end{proof}

\begin{Lemma}
\label{lem trace}
For any $f\in A_{p,F}$, the trace of the linear operator $L_f:A_{p,F} \to A_{p,F}$, $h\mapsto fh$,
is given by the formula ${\rm tr}\,L_f = (f,[\HessF])_p$.
\end{Lemma}
\begin{proof}
We have ${\rm tr}\,L_1= \mu_p=(1,[\HessF])_p$ and
${\rm tr}\,L_f= 0 =(f,[\HessF])_p$ for any $f\in\frak m_p$. This proves the lemma.
\end{proof}

Let $h_1^*,\dots,h^*_{\mu_p}$ be a $\C$-basis of $A_{p,F}$ dual to
$h_1,\dots,h_{\mu_p}$
with respect to the form $(\,,\,)_p$. Then $\sum\limits_{l=1}^{\mu_p} ([\HessF],h_l^*)_ph_l=
[\HessF]$.
Indeed, for  any $f\in A_{p,F}$ we have $\sum\limits_{l=1}^{\mu_p} (f,h_l^*)_ph_l=f$.

\begin{Lemma}
\label{lem gg=H}
We have $\sum\limits_{l=1}^{\mu_p} h_lh_l^* = [\HessF]$.
\end{Lemma}

\begin{proof}
For any $f\!\in\! A_{p,F}$, we
have ${\rm tr}\, L_f\! = \sum\limits_{l=1}^{\mu_p} (h_l^*,fh_l)_p = \sum\limits_{l=1}^{\mu_p} (h_l^*h_l,f)_p$
and ${\rm tr}\,L_f\! = (f,[\HessF])_p$. This proves the lemma.
\end{proof}

\begin{Lemma}
\label{lem ()S nond}
Let $s\in A_{p,F}$ be the element defined in Lemma {\rm \ref{lem F el}}. Then $s$
has the following two properties:
\begin{enumerate}\itemsep=0pt
\item[$(i)$]

the element $s$   is invertible and projects to $(-1)^k$ in $\C= A_{p,F}/\frak m_p$,

\item[$(ii)$] the form $(\,,\,)_S$ is nondegenerate.
\end{enumerate}

\end{Lemma}
\begin{proof}
To prove the lemma it is enough to show that
 $(-1)^k[\HessF]= s[\HessF]$.
Indeed, on one hand we have
\begin{gather*}
(f,g)_S = \sum_{l,m=1}^{\mu_p} S(w_l,w_m) (f,h_l)_p(g,h_m)_p.
\end{gather*}
 On the other hand we have
\begin{gather*}
(f,g)_S = (sf,g)_p = \sum_{l=1}^{\mu_p} (sf,h_l)_p(g,h_l^*)_p.
\end{gather*}
Hence, $\sum\limits_{l,m=1}^{\mu_p} S(w_l,w_m) h_l \otimes h_m =
 \sum\limits_{l=1}^{\mu_p} sh_l\otimes h_l^*$ in $A_{p,F}\otimes A_{p,F}$.
Therefore,
$\sum\limits_{l,m=1}^{\mu_p} S(w_l,w_m) h_l h_m =
 \sum\limits_{l=1}^{\mu_p} sh_lh_l^*$. By Assumption~\eqref{ass 1} and Lemma~\ref{lem gg=H}
we obtain $(-1)^k[\HessF]= s[\HessF]$.
\end{proof}

Let us prove Theorem \ref{model thm}.
Assume that $\sum\limits_{l=1}^{\mu_p} \la_lw_l = 0$ with $\la_l\in\C$. Denote
$h = \sum\limits_{l=1}^{\mu_p} \la_l h_l^*$.  Then $\al(h) = 0$ and $(f,h)_S = S(\al(f),
\al(h)) = 0$ for all $f\in A_{p,F}$. Hence $h=0$ since $(\,,\,)_S$ is
nondegenerate.  Therefore, $\la_l = 0$ for all $l$ and the vectors
$v_1,\dots,v_{\mu_\T}$ are linearly independent. This proves part~$(i)$ of Theorem~\ref{model thm}. Parts~$(ii)$ and $(iii)$ follow from Lemma \ref{hom lem}.

\subsection{Remark on maximal commutative subalgebras}
\label{Remark on maximal commutative subalgebras}

Let $A$ be a commutative algebra with unity element 1.
Let $B$ be the subalgebra of $\End(A)$ generated by all multiplication operators
$L_f:A \to A,\,h\mapsto fh$, where $f\in A$.

\begin{Lemma}
\label{lem on max abstract}
The subalgebra $B$ is a maximal commutative subalgebra of $\End(A)$.
\end{Lemma}

\begin{proof}
Let $T\in\End(A)$ be such that $[T,L_f]=0$ for all $f\in A$. Then $T = L_{T(1)}$.
\end{proof}

\begin{cor}
\label{cor abstract on max}
Under assumptions of Section {\rm \ref{An abstract setting}}, the algebra
$A_Y$ is a maximal commutative subalgebra of $\End(Y)$.
\end{cor}

\subsection{Interpretation of the algebra of Hamiltonians of good f\/ibers}
Under notations  of Section \ref{sec par trans} f\/ix a point $z\in\C^n-\Delta$.
Recall that in formula \eqref{Ham good fibers} we have def\/ined
 the quantum integrable model assigned to
$(\A(z), a)$ to be the collection
\[
\big(\sing W;\ S^{(a)}|_{\sing W};\ K_1(z)^*|_{\sing W},\dots, K_n(z)^*|_{\sing W} : \ \sing W\to \sing W\big) .
\]

Let $p \in U(\A(z))$ be an isolated critical point of the master
function
$\Phi(z,\,\cdot\,) : U(\A(z)) \to \C$.  Let $A_{p,\Phi}$ be the local
algebra of the critical point and $(\,,\,)_p$
the residue bilinear form on $A_{p,\Phi}$.
Let $[\bs E] \in A_{p,\Phi}\otimes \sing W$
be the element corresponding to the canonical element, see
Section~\ref{sec Hams and cr pts}.

Def\/ine a linear map
\begin{gather*}
\al_p : \ A_{p,\Phi} \to \sing W,
\qquad
g \mapsto  (g,[\bs E])_p .
\end{gather*}
Denote $Y_p$ the image of $\al_p$.

\begin{thm}\label{thm geom interpr good fibers}\quad
\begin{enumerate}\itemsep=0pt
\item[$(i)$]
We have $\ker\,\al_p = 0$.

\item[$(ii)$]
The operators $ K_j(z)^*$ preserve $Y_p$.
Moreover, for  any $j\in J$, $g\in A_{p,\Phi}$, we have
\[
\al_p(ga_j/[f_j(z,\cdot)]) = K_j(z)^*\al_p(g).
\]
\item[$(iii)$]
Define the value $[\bs E](p)$
of $[\bs E]$ at $p$ as the image of $[\bs E]$ under the natural projection
$A_{p,\Phi} \otimes \sing W \to  A_{p,\Phi}/\frak m_p \otimes \sing W = \sing W$.
Then the value $[\bs E](p)$  is nonzero.
The vector $[\bs E](p)$ is the only $($up to proportionality$)$
common eigenvector of the operators $K_j(z)^*|_{Y_p} : Y_p \to Y_p$,  $j\in J$,
and we have
\[
K_j(z)^*[\bs E](p) = \frac{a_j}{f_j(z,p)}[\bs E](p) .
\]
\end{enumerate}
\end{thm}

\begin{proof}
By Theorem \ref{thm v good fiber} and Remark~\ref{rem a_j/f_j form basis]}
the objects $\sing W, S^{(a)}|_{\sing W},$ $ K_j(z)^*|_{\sing W}, [{\bs v}],$
 $ [a_j/f_j(z,\cdot)]$,
$j\in J$, satisfy the assumptions of Theorem~\ref{model thm}. Now Theorem~\ref{thm geom interpr good fibers} follows
from Theorem \ref{model thm}.
\end{proof}

\begin{thm}
\label{thm good fibers (,)=S}
The linear map $\al_p$ identifies the contravariant form on $Y_p$ and the residue form $(\,,\,)_p$ on
$A_{p,\Phi}$ multiplied by $(-1)^k$,
\begin{gather}
\label{S=(,)}
S^{(a)}(\al_p(f),\al_p(g)) = (-1)^k(f,g)_p
\end{gather}
for any $f,g\in A_{p,\Phi}$.
\end{thm}

\begin{proof}
If the Milnor number of $p$ is one, then the theorem follows from
Lemma~\ref{lem ()S nond}. If the Milnor number is greater than one, the theorem follows by continuity from the case of
the Milnor number equal to one, since all objects involved depend continuously
on the weights $a$ and parameters $z$.

Note that the theorem says that the element $s$ of Lemma~\ref{lem ()S nond} in our situation equals $(-1)^k$.
\end{proof}

\begin{rem}
\label{rem integral identity} \rm
In formula~\eqref{S=(,)}, each of $\al_p(f)$, $\al_p(g)$ is given by the Grothendieck residue,
so each of~$\al_p(f)$,~$\al_p(g)$
is a $k$-dimensional integral. The quantity $(f,g)_p$ is also a $k$-dimensional integral. Thus formula
\eqref{S=(,)} is an equality relating a bilinear expression in $k$-dimensional integrals to an individual
$k$-dimensional integral.
\end{rem}

Denote by $A_{Y_p}$ the unital subalgebra of $\End(Y_p)$ generated by
$ K_j(z)^*\vert_{Y_p}$, $j\in J$.

\begin{cor}\label{cor hams good fiber}\quad
\begin{enumerate}\itemsep=0pt
\item[$(i)$]
The map $[a_j/f_j(z,\cdot)] \mapsto  K_j(z)^*\vert_{Y_p}$, $j\in J$,
extends uniquely to an algebra isomorphism
$\beta_p : A_{p,\Phi} \to A_{Y_p}$.
\item[$(ii)$]
The maps $\al_p$, $\beta_p$ give an isomorphism of the regular representation
of $A_{p,\Phi}$ and
the $A_{Y_p}$-module  $Y_p$, that is $\beta_p(h)\al_p(g) = \al_p(hg)$
for any $h,g\in A_{p,\Phi}$.
\item[$(iii)$] The algebra $A_{Y_p}$ is a maximal commutative subalgebra of
$\End(Y_p)$.
\item[$(iv)$] All elements of the  algebra $A_{Y_p}$ are symmetric operators with respect to the contravariant form $S^{(a)}$.
\end{enumerate}
\end{cor}

\begin{thm}
\label{thm many point}
Let $p_1,\dots,p_d$ be a list of all
distinct isolated critical points  of $\Phi(z,\,\cdot\,)$.
Let $Y_{p_s}=\al_{p_s}(A_{p_s,\Phi})\subset \sing W$, $s=1,\dots,d$, be the corresponding subspaces.
Then the sum of these subspaces is direct. The subspaces are orthogonal.
\end{thm}

\begin{proof}
It follows from Theorem \ref{thm geom interpr good fibers} that for any $s=1,\dots,d$ and $j\in J$ the
operator
 $K_{j}(z)^* - a_j/f_j(z,p_s)$ restricted to $Y_{p_s}$ is nilpotent. We also know that
 the numbers
$a_j/f_j(z,p_s)$ separate the points $p_1,\dots,p_d$.
These observations imply Theorem~\ref{thm many point}.
\end{proof}

\begin{cor}
\label{cor sum of Milnor}
The sum of Milnor numbers of the critical points $p_1,\dots,p_d$ is not greater than
the rank of the contravariant form $S^{(a)}\vert_{\Sing W}$.
\end{cor}

Denote  $Y = \oplus_{s=1}^d Y_{p_s}$. Denote
by $A_Y$ the unital subalgebra of $\End(Y)$ generated by
$ K_j(z)^*\vert_{Y}$, $j\in J$.
Consider the isomorphisms
\begin{gather*}
\al =\oplus_{s=1}^d\al_{p_s}  :\
\oplus_{s=1}^d
A_{p_s, \Phi}
\to  \oplus_{s=1}^d Y_{p_s}  ,
\qquad
\beta=\oplus_{s=1}^d\beta_s   :\ \oplus_{s=1}^dA_{p_s, \Phi}  \to  \oplus_{s=1}^d A_{Y_{p_s}} .
\end{gather*}

\begin{cor}
\label{cor main}
We have
\begin{enumerate}\itemsep=0pt
\item[$(i)$]
$A_Y=\oplus_{s=1}^d A_{p_s,Y_{p_s}}$;
\item[$(ii)$]
$A_Y$ is a maximal commutative subalgebra of $\End(Y)$.
\item[$(iii)$]
The isomorphisms
$\al$, $\beta$
identify the regular representation
of the algebra $\oplus_{s=1}^dA_{p_s, \Phi}$ and the $A_Y$-module $Y$.
The isomorphism
$\al$
identifies the contravariant form on $Y$ and the residue form $(\,,\,) = \oplus_{s=1}^d (\,,\,)_{p_s}$ on
$\oplus_{s=1}^d A_{p_s, \Phi}$ multiplied by $(-1)^k$.

\item[$(iv)$]
In particular, if the dimension of $\sing W$  equals the sum of Milnor numbers
 $\sum\limits_{s=1}^d \mu_s$, then the module $\sing W$ over the unital subalgebra
of $\End (\sing W)$ generated by geometric
Hamiltonians $ K_j(z)^*|_{\sing W} : \sing W\to \sing W$,\ $j\in J$, is isomorphic to the regular
representation of the algebra $\oplus_{s=1}^dA_{p_s, \Phi}$.
\item[$(v)$]
If for $z\in \C^n-\Delta$ the arrangement $(\A(z),a)$ is unbalanced, then
the module $\sing W$ over the unital subalgebra
of $\End (\sing W)$ generated by geometric Hamiltonians $ K_j(z)^*|_{\sing W}$, $j\in J$, is isomorphic to the regular
representation of the algebra $\oplus_{s=1}^dA_{p_s, \Phi}$.
\end{enumerate}
\end{cor}

\begin{cor}
If for $z\in \C^n-\Delta$ the arrangement $(\A(z),a)$ is unbalanced, then
the contrava\-riant form is nondegenerate on $\Sing \FF^k(\A(z))$.
\end{cor}

\begin{proof} Indeed in this case the sum of Milnor numbers of critical points of the master function equals
$|\chi(U(\A(z)))|$ and equals $\dim  \Sing \FF^k(\A(z))$.
\end{proof}

\subsection[Interpretation of the algebra of Hamiltonians of bad f\/ibers if Assumption~\ref{i} is satisf\/ied]{Interpretation of the algebra of Hamiltonians of bad f\/ibers\\ if Assumption~\ref{i} is satisf\/ied}
\label{sec Algebra of Hamiltonians of bad fibers}

Let Assumption \ref{i} be satisf\/ied.
Fix  $z^0\in \Delta$.
We have
def\/ined
 the quantum integrable model assigned to
$(\A(z^0), a)$ to be the collection
\begin{gather*}
\big(\Sing \FF^k\big(\A\big(z^0\big)\big);\ S^{(a)}|_{\Sing \FF^k(\A(z^0))}; \\
\qquad
 \tilde K_1\big(z^0\big),\dots,\tilde K_n\big(z^0\big) : \ \Sing \FF^k\big(\A\big(z^0\big)\big)\to \Sing \FF^k\big(\A\big(z^0\big)\big)\big) ,
\end{gather*}
 where
$\tilde K_j(z^0)= {\rm pr} K_j^1(z^0)|_{\Sing\FF^k(\A(z^0))}$, see Theorem~\ref{thm bar K_j}.

Let $p \in U(\A(z^0))$ be an isolated critical point of the
master function $\Phi(z^0,\,\cdot\,)
: U(\A(z^0)) \to \C$.
Let $A_{p,\Phi}$ be the local algebra of the critical
point
and $(\,,\,)_p$
the residue bilinear form on $A_{p,\Phi}$.
Let $[E] \in A_{p,\Phi}\otimes \Sing \FF^k(\A(z^0))$
be the canonical element corresponding to the arrangements $\A(z^0)$, see
Section~\ref{sec Hams and cr pts, bad fibers}.

Def\/ine a linear map
\begin{gather*}
\al_p : \ A_{p,\Phi} \to \Sing \FF^k\big(\A\big(z^0\big)\big),
\qquad
g \mapsto  (g,[ E])_p .
\end{gather*}
Denote $Y_p$ the image of $\al_p$.

\begin{thm}
\label{thm geom interpr bad fibers}
Let Assumption {\rm \ref{i}} be satisfied. Then
\begin{enumerate}\itemsep=0pt
\item[$(i)$]
We have $\ker\,\al_p = 0$.
The isomorphism
$\al_p$
identifies the contravariant form on $Y_p$ and the residue form $(\,,\,)_p$ on
$A_{p, \Phi}$ multiplied by $(-1)^k$.

\item[$(ii)$]
 For any $j\in J$, the operator $ \tilde K_j(z^0)$ preserves $Y_p$.
Moreover, for  any
$j\in J$ and $g\in A_{p,\Phi}$, we have
$\al_p(ga_j/[f_j(z^0,\cdot)]) = \tilde K_j(z^0)\al_p(g)$.
\item[$(iii)$]
Define the value $[ E](p)$
of $[ E]$ at $p$ as the image of $[ E]$ under the natural projection
$A_{p,\Phi} \otimes \Sing \FF^k(\A(z^0)) \to  A_{p,\Phi}/\frak m_p \otimes \Sing \FF^k(\A(z^0)) = \Sing \FF^k(\A(z^0))$.
Then the value $[ E](p)$  is nonzero.
The vector $[ E](p)$ is the only $($up to proportionality$)$
common eigenvector of the operators $\tilde K_j(z^0)|_{Y_p} : Y_p\to Y_p$,  $j\in J$,
and we have
\[
\tilde K_j(z^0)[E](p) = \frac{a_j}{f_j(z^0,p)}[E](p) .
\]
\end{enumerate}
\end{thm}

Denote by $A_{Y_p}$ the unital subalgebra of $\End(Y_p)$ generated by
$\tilde K_j(z^0)\vert_{Y_p}$, $j\in J$.

\begin{cor}
\label{cor hams bad fiber}
Let Assumption {\rm \ref{i}} be satisfied. Then
\begin{enumerate}\itemsep=0pt
\item[$(i)$]
The map $[a_j/f_j(z^0,\cdot)] \mapsto \tilde K_j(z^0)\vert_{Y_p}$, $j\in J$,
extends uniquely to an algebra isomorphism
$\beta_p : A_{p,\Phi} \to A_{Y_p}$.
\item[$(ii)$]
The maps $\al_p$, $\beta_p$ give an isomorphism of the regular representation
of $A_{p,\Phi}$ and
the $A_{Y_p}$-module  $Y_p$, that is $\beta_p(h)\al_p(g) = \al_p(hg)$
for any $h,g\in A_{p,\Phi}$.
\item[$(iii)$] The algebra $A_{Y_p}$ is a maximal commutative subalgebra of
$\End(Y_p)$.
\end{enumerate}
\end{cor}

Recall that under Assumption~\ref{i}
all critical points of $\Phi(z^0,\,\cdot\,)$ are isolated, the sum of their
Milnor numbers equals $\dim \Sing \FF^k(\A(z^0))$ and the form $S^{(a)}\vert_{\Sing \FF^k(\A(z^0))}$ is nondegenerate.

\begin{thm}
\label{thm many point bad}
Let Assumption {\rm \ref{i}} be satisfied.
Let $p_1,\dots,p_d$ be a list of all
distinct critical points  of $\Phi(z^0,\,\cdot\,)$.
Let $Y_{p_s}=\al_{p_s}(A_{p_s,\Phi})\subset \Sing \FF^k(\A(z^0))$, $s=1,\dots,d$, be the corresponding subspaces.
 Then the sum of these subspaces is direct, orthogonal and equals $\Sing \FF^k(\A(z^0))$.
\end{thm}

Denote by
$A_{\A(z^0),a}$ the unital subalgebra of $\End(\Sing \FF^k(\A(z^0)))$ generated by
$\tilde K_j(z^0)$, $j\in J$. The algebra $A_{\A(z^0),a}$
is called the algebra of geometric Hamiltonian of the arrangement $(\A(z^0),a)$, see Section~\ref{Main Conjecture}.

Consider the isomorphisms
\begin{gather*}
\al  = \oplus_{s=1}^d\al_{p_s}  :\
\oplus_{s=1}^d
A_{p_s, \Phi}
\to  \Sing \FF^k\big(\A\big(z^0\big)\big)=\oplus_{s=1}^d Y_{p_s} ,
\\
\beta  = \oplus_{s=1}^d\beta_s   :\ \oplus_{s=1}^d
A_{p_s, \Phi}  \to  \oplus_{s=1}^d A_{Y_{p_s}} .
\end{gather*}

\begin{cor}
\label{cor main bad}
Let Assumption {\rm \ref{i}} be satisfied. Then
\begin{enumerate}
\item[$(i)$]
$A_{\A(z^0),a}=\oplus_{s=1}^d A_{Y_{p_s}}$.
\item[$(ii)$]
$A_{\A(z^0),a}$ is a maximal commutative subalgebra of $\End(\Sing \FF^k(\A(z^0)))$.
\item[$(iii)$]
The isomorphisms
$\al$, $\beta$
identify the regular representation
of the algebra $\oplus_{s=1}^dA_{p_s, \Phi}$ and the $A_{\A(z^0),a}$-module $\Sing \FF^k(\A(z^0))$.
The isomorphism
$\al$
identifies the contravariant form on $\Sing \FF^k(\A(z^0))$ and the residue form $(\,,\,) = \oplus_{s=1}^d (\,,\,)_{p_s}$ on
$\oplus_{s=1}^d A_{p_s, \Phi}$ multiplied by $(-1)^k$.
\end{enumerate}
\end{cor}

\section{More on Hamiltonians of bad f\/ibers}
\label{More on Hamiltonians of bad fibers}

\subsection{An abstract setting}
\label{An abstract setting, 2}
Let $k < n$ be positive integers and $J=\{1,\dots,n\}$ as before.

Let $B\subset\C^k$ be a ball with center at a point $p$.
Let $F_u$ be a
holomorphic function on $B$ depen\-ding holomorphically on a complex
parameter $u$ at $u=0$.  Assume that $F_0=F_{u=0}$ has a single critical point at $p$ with  Milnor number $\mu$.
 Let $\C(t)_B$ be the algebra of holomorphic functions  on~$B$,  $I_u\subset \C(t)_B$ the ideal generated by
$\partial F_u/\partial t_i$, $i=1,\dots,k$, and $A_u = \C(t)_B/I_u$.
Assume that $\dim_\C A_u$ does not depend on $u$ for $u$ in a neighborhood of $0$.
  Let $[\,]_u : \C(t)_B\to
A_u$ be the cano\-ni\-cal projection,
 $(\,,\,)_u$ the residue
bilinear form on $A_u$ and
$\HessF_u=\det_{1\leq l,m\leq k}(\der^2F_u/\der t_l\der t_m)$.

Let $h_1,\dots,h_{\mu}\in\C(t)_B$ be a collection of elements such that
for any $u$ the elements
$[h_1]_u,  \dots$, $[h_{\mu}]_u$ form a $\C$-basis of $A_u$.

Let $g_{1,u},\dots,g_{n,u} \in \C(t)_B$ be elements depending on $u$ holomorphically
at $u=0$ and
such that
for any $u$ (close to $0$) the unital subalgebra of $A_u$ generated by $[g_{1,u}]_u,
\dots,[g_{n,u}]_u$
equals~$A_u$.

Let $W$ be a vector space with a
symmetric bilinear form $S$. For $u\neq 0$, let
$M_{j,u} : W\to W$, $j\in J$, be
a collection of commuting symmetric linear operators,
\begin{gather*}
M_{i,u}M_{j,u}=M_{j,u}M_{i,u},\qquad
 S(M_{j,u}x,y) = S(x,M_{j,u}y)
\qquad
{\rm for\ all}\ i,j\in J \ {\rm and}\ x,y\in W.
\end{gather*}
We assume that every $M_{j,u}$ depends on $u$ meromorphically (for $u$ close to $0$) and has at
most simple pole at $u=0$,
\begin{gather}
\label{M -1 0}
M_{j,u} =\frac {M_{j}^{(-1)}}u + M_{j}^{(0)} +  M_{j;1}^{(1)}u +\cdots ,
\qquad
 M_{j}^{(i)} \in \End ( W) .
\end{gather}
Let $w_{1,u},\dots,w_{\mu,u}\in W$ be a collection of vectors depending on
$u$ holomorphically at $u=0$.

Consider the element
\[
[w]_u = \sum_{l=1}^\mu [h_l]_u\otimes w_{l,u} \in A_u \otimes W .
\]
Assume that for every nonzero $u$ (close to $0$) we have
\begin{gather}
\label{ass 22}
\sum_{l=1}^\mu [h_l]_u\otimes M_{j,u}w_{l,u}  =   \sum_{l=1}^\mu [g_{j,u}]_u[h_l]_u\otimes w_{l,u} ,
\qquad j\in J,
\end{gather}
and for every $u$  (close to $0$) we have
\begin{gather*}
 \sum_{l,m=1}^{\mu_p} S(w_{l,u},w_{m,u}) [h_l]_u [h_m]_u   =  (-1)^k[\HessF_u]_u .
\end{gather*}

For any $u$,
 denote by $Y_u\subset W$ the vector subspace generated by $w_{1,u},\dots,w_{\mu,u}$.
By pro\-per\-ty~\eqref{ass 22}, for  any nonzero $u$ (close to~$0$) every $M_{j,u}$, $j\in J$, preserves~$Y_u$.
For  any nonzero~$u$ (close to~$0$) denote by~$A_{Y_u}$ the unital subalgebra of $\End(Y_u)$ generated by $M_{j,u}\vert_{Y_u}$, $j\in J$.
 The subspace~$Y_u$ is an $A_{Y_u}$-module.

For any $u$ ($u=0$ included)
def\/ine a linear map
\begin{gather*}
\al_u : \ A_{u} \to Y_u,
\qquad
[f]_u \mapsto \sum_{l=1}^{\mu} ([f]_u,[h_l]_u)_u w_{l,u} .
\end{gather*}

\begin{thm}
\label{model thm 2}
For any $u$ $($in particular, for $u=0)$
the map $\al_u : A_{u} \to Y_u$ is an isomorphism of vector spaces.
The form $S$ restricted to $Y_u$ is nondegenerate.
\end{thm}

\begin{proof}
Def\/ine a bilinear form $(\,,\,)_{S,u}$ on $A_{u}$,
\begin{gather*}
([f]_u,[g]_u)_{S,u}  =
S(\al_u([f]_u),\al_u([g]_u))  =
\sum_{l,m=1}^{\mu}  S(w_{l,u},w_{m,l})  ([f]_u,[h_l]_u)_u ([g]_u,[h_m]_u)_u  .
\end{gather*}

\begin{Lemma}
\label{lem ()S 2} For any $u$, we have $([f]_u[g]_u,[h]_u)_{S,u}\! =\!
([f]_u,[g]_u[h]_u)_{S,u}$
 for all
$[f]_u,[g]_u,[h]_u \in A_u$.
\end{Lemma}

\begin{proof}
For $u\neq 0$, the statement  follows from
Lemma~\ref{lem ()S}. For $u=0$, the statement follows by continuity.
\end{proof}

The next two lemmas are similar to the corresponding analogs in
Section~\ref{ Proof of Theorem model thm}.

\begin{Lemma}
\label{lem F el 2}
There exists a unique element
 $[s]_0\in A_{0}$ such that $([f]_0,[g]_0)_{S,0} =  ([s]_0[f]_0,[g]_0)_0$
for all $[f]_0,[g]_0\in A_0$.
\end{Lemma}

\begin{Lemma}
\label{lem ()S nond 2}
Let $[s]_0\in A_{0}$ be the element defined in Lemma~{\rm \ref{lem F el 2}}. Then $[s]_0$
has the following two properties:
\begin{enumerate}\itemsep=0pt
\item[$(i)$]

the element $[s]_0$   is invertible and projects to $(-1)^k$ in $\C= A_{0}/\frak m_0$,
where $\frak m_0 \subset A_{0}$ is the maximal ideal,

\item[$(ii)$] the form $(\,,\,)_{S,0}$ is nondegenerate.
\end{enumerate}
\end{Lemma}

Lemma \ref{lem ()S nond 2} implies Theorem \ref{model thm 2}, cf.\ the end of Section
\ref{ Proof of Theorem model thm}.
\end{proof}

Def\/ine the value $[w]_0(p)$ of $[w]_0$ at $p$ as the image of $[w]_0$ under the natural projection
$A_{0} \otimes W \to  A_{0}/\frak m_0 \otimes W = W$.

\begin{cor}
\label{cor value at 0}
 The value $[w]_0(p)$ is nonzero.
\end{cor}

\begin{cor}
\label{cor wl0}
The space $Y_0$ is of dimension $\mu$ and $w_{1,0}, \dots, w_{l,0}$ is its basis.
\end{cor}

For any $[g]_0\in A_0$, denote by $L_{[g]_0}\in\End(A_0)$ the linear operator on $A_0$
of multiplication by~$[g]_0$.
For any $j\in J$, def\/ine a linear map $\bar L_{j,0} :Y_0\to Y_0$ by the formula
$\bar L_{j,0} = \al_0 L_{[g_{j,0}]_0}(\al_0)^{-1}$.
Denote by $A_{Y_0}$ the unital subalgebra of $\End(Y_0)$ generated by
$\bar L_{j,0}$, $j\in J$. Clearly, $A_{Y_0}$ is commutative.
The subspace $Y_0$ is an $A_{Y_0}$-module.

\begin{thm}
\label{model thm 3}
\quad
\begin{enumerate}\itemsep=0pt
\item[$(i)$]
The map $[g_{j,0}]_0 \mapsto \bar L_{j,0}$, $j\in J$,
extends uniquely to an algebra isomorphism
$\beta_0 : A_{0} \to A_{Y_0}$.
\item[$(ii)$]
the algebra $A_{Y_u}$ tends to the algebra $A_{Y_0}$ as $u\to 0$. More precisely,
for any $j\in J$ and $l=1,\dots,\mu$, we have $M_{j,u}w_{l,u} \to \bar L_{j,0}w_{l,0}$
as $u\to 0$.
\end{enumerate}
\end{thm}

\begin{proof}
Part $(i)$ is clear. Part $(ii)$ follows from \eqref{ass 22}.
\end{proof}

Let $\tilde Y \subset W$ be a vector subspace such that
\begin{enumerate}\itemsep=0pt
\item[$(a)$]
 $Y_0\subset \tilde Y$;
\item[$(b)$]
 the bilinear form $S$ restricted on $\tilde Y$ is nondegenerate;
 \item[$(c)$] for any $j\in J$,
 the subspace
$\tilde Y$ lies in the kernel of $M_j^{(-1)}$.
\end{enumerate}
For example, we can choose $\tilde Y=Y_{0}$.
Let ${\rm pr}_{\tilde Y} :W\to \tilde Y$ be the orthogonal projection.

\begin{thm}
\label{thm pr Y}
\quad
\begin{enumerate}\itemsep=0pt
\item[$(i)$]
For $j\in J$, let $M_{j}^{(0)}$ be the constant coefficient of the Laurent
expansion of $M_{j,u}$, see~\eqref{M -1 0}. Then
\begin{gather}
\label{good formula}
\bar L_{j,0} = {\rm pr}_{\tilde Y} M_{j}^{(0)}\vert_{Y_0} .
\end{gather}
In particular, that means that the operators ${\rm pr}_{\tilde Y} M_{j}^{(0)}\vert_{Y_0}$ do not depend on the choice of $\tilde Y$.

\item[$(ii)$]
Let $(\xi_j)_{j\in J}\subset \C$ be numbers such that $\sum_{j\in J}\xi_j M^{(-1)}_j=0$, then
\begin{gather*}
\sum_{j\in J}\xi_j\bar L_{j,0} = \sum_{j\in J}\xi_j M_{j}^{(0)}\vert_{Y_0} .
\end{gather*}
\end{enumerate}
\end{thm}

\begin{proof}
For $j\in J$, let $w_{l,u} = w_{l,0} + w_l^{(1)}u + \cdots$ be the Taylor expansion
of $w_{l,u}$. By part $(ii)$ of Theorem~\ref{model thm 3} we have
\begin{gather}
\label{KER}
w_{j,0} \in \ker M_j^{(-1)}
\qquad
{\rm and}
\qquad
\bar L_{j,0}w_{l,0}
 = M_{j}^{(0)}w_{l,0} +
M_{j}^{(-1)}w_{l,0}^{(1)} .
\end{gather}
The operator $M_{j}^{(-1)}$ is symmetric since $M_{j,u}$ is symmetric. The image of
$M_{j}^{(-1)}$ is orthogonal to the kernel of $M_{j}^{(-1)}$. Hence, formula~\eqref{good formula}
follows from formula~\eqref{KER}. Part $(ii)$ of the theorem also follows from formula~\eqref{KER}.
\end{proof}

\begin{cor}
\label{cor on COMMUT}
For any $i,j\in J$, the operators
${\rm pr}_{\tilde Y}  M_{i}^{(0)}\vert_{Y_0}$, $ {\rm pr}_{\tilde Y} M_{j}^{(0)}\vert_{Y_0} $ are symmetric and
commute.
\end{cor}

\begin{cor}
\label{cor on COMMUT MAX}
The unital subalgebra of ${}$ $\End (Y_0)$ generated by
the operators
${\rm pr}_{\tilde Y} M_{j}^{(0)}\vert_{Y_0}$, $j\in J$,
is a maximal
commutative subalgebra.
\end{cor}

\begin{proof}
The proof follows from remarks in Section~\ref{Remark on maximal commutative subalgebras}.
\end{proof}

\begin{cor}
\label{cor unique eigenve}
The vector $[w]_0(p)$ is the only $($up to proportionality$)$
common eigenvector of the operators ${\rm pr}_{\tilde Y} M_{j}^{(0)}\vert_{Y_0}$, $j\in J$,
and we have
${\rm pr}_{\tilde Y} M_{j}^{(0)}[w]_0(p) = g_{j,0}(p)[w]_0(p)$.
\end{cor}

Assume that the parameter $u$ is changed holomorphically,
$u= c_1v+c_2v^2+\cdots $ where $c_i\in\C$, $c_1\neq 0$, and~$v$ is a new parameter. Let
\[
M_{j,u(v)} =\frac {\tilde M_{j}^{(-1)}}v + \tilde M_{j}^{(0)} +  \tilde M_{j;1}^{(1)}v +\cdots ,
\qquad
\tilde M_{j}^{(i)} \in \End ( W) ,
\]
be the  new Laurent expansion.

\begin{Lemma}
\label{rem on invariance}
For any $j\in J$, we have
${\rm pr}_{\tilde Y} M_{i}^{(0)}\vert_{Y_0}={\rm pr}_{\tilde Y} \tilde M_{i}^{(0)}\vert_{Y_0}$ and the algebra
$A_{Y_0} \subset \End(Y_0)$ does not change under the reparametrization of $u$.
\end{Lemma}

\begin{proof}
One proof of the lemma follows from Theorem~\ref{model thm 3}. Another proof follows from the fact
that $\tilde M_{i}^{(-1)} =  M_{i}^{(-1)}/c_1$,
$\tilde M_{i}^{(0)} =  M_{i}^{(0)} - c_2 M_{i}^{(-1)}/c_1^2$.
\end{proof}

\subsection{Hamiltonians of bad f\/ibers}
\label{Algebra of bad 2}

Let us return to the situation of Sections \ref{sec par trans} and \ref{sec real posit bad}
and recall the previous constructions.

Let $z^0\in \Delta$.
Let $V\to V^*$ be the map associated with the contravariant form. Let $W$, $\sing W$, $W(z^0)$,
$\sing W(z^0)$ be the images of~$V$, $\sing V$, $\FF^k(\A(z^0))$, $\sing \FF^k(\A(z^0))$, respectively.
The contravariant form on~$V$ induces a nondegenerate
symmetric bilinear form on $W$ also denoted by~$S^{(a)}$.

For $z\in \C^n-\Delta$, we have linear operators
$K_j(z): V \to  V$, $j\in J$, where
$K_j(z)  =  \sum\limits_{C\in \frak C}  \frac{\la_j^C}{f_C(z)}  L_C$, see~\eqref{K_j}.
For
$\frak C_0 = \{ C\in \frak C\ | \ z^0\in H_C\}$, we def\/ine
\[
K_j^0(z)   =   \sum_{C\in \frak C_0}
 \frac{\la_j^C}{f_C(z)}   L_C  ,
\qquad
K_j^1(z) = K_j(z) - K_j^0(z) .
\]
The dual operators $K_j(z)^* : V^*\to V^*$ preserve the subspaces $\sing W\subset W\subset V^*$, commute on the subspace
$\sing W$ and are symmetric on $W$ with respect to the contravariant form.
The operators $L_C^* : V^*\to V^*$, $C\in\frak C$, preserve the subspaces $\sing W\subset W\subset V^*$.
The space $W(z^0)$ lies in the kernel of $L_C^*|_W : W\to W$ for any $C \in \frak C_0$.

Let $T_{z^0} = \{
\xi = \sum_{j\in J} \xi_j\frac{\der}{\der z_j}\ |\ \xi_j\in \C, \
\xi (f_C) = 0 \ \text{for all}\ C\in \frak C_0\}$.

Let $p \in U(\A(z^0))$ be an isolated critical point of the
master function $\Phi(z^0,\,\cdot\,)
: U(\A(z^0)) \to \C$.
Let $A_{p,\Phi}$ be the local algebra of the critical
point
and $(\,,\,)_p$
the residue bilinear form on $A_{p,\Phi}$.
Let $[ \ ] : \C(t)_{U(\A(z^0))} \to A_{p,\Phi}$ be the canonical projection
and $[\bs E] \in A_{p,\Phi}\otimes \sing W(z^0)$
the element corresponding to the canonical element.

Def\/ine a linear map
\begin{gather*}
\al_p : \  A_{p,\Phi} \to \sing W\big(z^0\big),
\qquad
g \mapsto  (g,[ \bs E])_p .
\end{gather*}
Denote $Y_p$ the image of $\al_p$.

\begin{thm}
\label{thm geom interpr bad fibers 2}
\quad
\begin{enumerate}\itemsep=0pt
\item[$(i)$]
We have
$\ker\,\al_p = 0$.
The isomorphism $\al_p$ identifies the contravariant form on~$Y_p$ and the residue form $(\,,\,)_p$ on $A_{p,\Phi}$ multiplied by $(-1)^k$.
In particular,
the contravariant form on~$Y_p$ is nondegenerate.

\item[$(ii)$]
Let $\tilde Y \subset \sing W$ be a vector subspace such that
\begin{enumerate}\itemsep=0pt
\item[$(a)$]
 $Y_p\subset \tilde Y$;
\item[$(b)$]
 the  contravariant form  restricted on $\tilde Y$ is nondegenerate;
 \item[$(c)$] for any $j\in J$,
 the subspace
$\tilde Y$ lies in the kernel of $L^*_C$, $C\in \frak C_0$.
\end{enumerate}
Let ${\rm pr}_{\tilde Y} :\sing W\to \tilde Y$ be the orthogonal projection.
Then for any $j\in J$ the operator ${\rm pr}_{\tilde Y}K^1_j(z^0)^*|_{Y_p}$ maps $Y_p$ to $Y_p$
and does not depend on the choice of $\tilde Y$.
The operators
 ${\rm pr}_{\tilde Y}K^1_j(z^0)^*|_{Y_p}:Y_p\to Y_p$, $j\in J$,
 commute and are symmetric with respect to the contravariant form on~$Y_p$.

\item[$(iii)$]
The unital subalgebra $A_{Y_p}\subset \End (Y_p)$ generated by
 ${\rm pr}_{\tilde Y}K^1_j(z^0)^*|_{Y_p}$, $j\in J$,
is a maximal commutative subalgebra.

\item[$(iv)$] The naive geometric Hamiltonians $K_\xi(z^0)^*$, $\xi\in T_{z^0}$, preserve the subspace $Y_p$ and the operators
$K_\xi(z^0)^*|_{Y_p}$ are elements of the subalgebra $A_{Y_p}$.

\item[$(v)$]
The value $[\bs E](p)$
of $[\bs E]$ at $p$  is nonzero.
The vector $[\bs E](p)$ is the only $($up to proportionality$)$
common eigenvector of the operators ${\rm pr}_{\tilde Y}K^1_j(z^0)^*|_{Y_p}$, $j\in J$,
and we have
\[
{\rm pr}_{\tilde Y}K^1_j\big(z^0\big)^*[\bs E](p) = \frac{a_j}{f_j(z^0,p)}[\bs E](p) .
\]

\item[$(vi)$]
For any $j\in J$, $g\in A_{p,\Phi}$, we have
$\al_p(ga_j/[f_j(z^0,\cdot)]) = {\rm pr}_{\tilde Y}K^1_j(z^0)^*\al_p(g)$.
The map $[a_j/f_j(z^0,\,)] \mapsto {\rm pr}_{\tilde Y}K^1_j(z^0)^*|_{Y_p}$, $j\in J$,
extends uniquely to an algebra isomorphism
$\beta_p : A_{p,\Phi} \to A_{Y_p}$.

\item[$(vii)$]
The isomorphisms $\al_p$, $\beta_p$ identify the regular representation of the algebra $A_{p,\Phi}$ and the
$A_{Y_p}$-module $Y_p$.
\end{enumerate}
\end{thm}

All statements of the theorem (but the second statement of part $(i)$) follow from the corresponding statements of Section~\ref{An abstract setting, 2}. The second statement of part $(i)$ has the same proof as Theorem~\ref{thm good fibers (,)=S}.

\begin{thm}
\label{thm many point bad 2}
Let $p_1,\dots,p_d$ be a list of all
distinct isolated critical points  of $\Phi(z^0,\,\cdot\,)$.
Let $Y_{p_s}=\al_{p_s}(A_{p_s,\Phi})\subset \sing W(z^0)$, $s=1,\dots,d$, be the corresponding subspaces.
Then the sum of these subspaces is direct and  orthogonal with respect to the contravariant form.
\end{thm}

\begin{cor}
\label{cor sum of Milnor 2}
The sum of Milnor numbers of the critical points $p_1,\dots,p_d$ is not greater than
the rank of the contravariant form $S^{(a)}\vert_{\Sing \FF^k(\A(z^0))}$.
\end{cor}

Denote $Y = \oplus_{s=1}^d Y_{p_s}$.
Let $\tilde Y \subset \sing W$ be a vector subspace such that
\begin{enumerate}\itemsep=0pt
\item[$(a)$]
 $Y\subset \tilde Y$;
\item[$(b)$]
 the contravariant form restricted on $\tilde Y$ is nondegenerate;
 \item[$(c)$] for any $j\in J$,
 the subspace
$\tilde Y$ lies in the kernel of $L_C^*$,  $C\in \frak C_0$.
\end{enumerate}
For example, we can choose $\tilde Y=Y$.
Let ${\rm pr}_{\tilde Y} :  \sing W\to \tilde Y$ be the orthogonal projection.

\begin{thm}
\label{thm BAD}
\quad
\begin{enumerate}\itemsep=0pt
\item[$(i)$]
For any $j\in J$ the operator
${\rm pr}_{\tilde Y}K^1_j(z^0)^*$ maps $Y$ to $Y$ and does not depend on the choice of~$\tilde Y$.
The operators
 ${\rm pr}_{\tilde Y}K^1_j(z^0)^*|_{Y}:Y \to Y$, $j\in J$,
 commute, preserve each of the subspaces~$Y_{p_s}$
 and are symmetric with respect to the contravariant form on~$Y$.

Denote by
$A_Y$ be the unital subalgebra of $\End(Y)$ generated by
${\rm pr}_{\tilde Y}K^1_j(z^0)^*\vert_{Y}$, $j\in J$.

\item[$(ii)$]

The naive geometric Hamiltonians $K_\xi(z^0)^*$, $\xi\in T_{z^0}$, preserve the subspace $Y$ and the operators
$K_\xi(z^0)^*|_{Y}$ are elements of the subalgebra $A_{Y}$.

\item[$(iii)$]
Consider the isomorphisms
\begin{gather*}
\al  = \oplus_{s=1}^d\al_{p_s}  :\
\oplus_{s=1}^d
A_{p_s, \Phi}
\to  \oplus_{s=1}^d Y_{p_s}  ,
\\
\beta  =  \oplus_{s=1}^d\beta_s   :\ \oplus_{s=1}^dA_{p_s, \Phi}  \to  \oplus_{s=1}^d A_{Y_{p_s}} .
\end{gather*}
Then
\begin{enumerate}\itemsep=0pt
\item[$(a)$]
$A_Y=\oplus_{s=1}^d A_{Y_{p_s}}$;
\item[$(b)$]
$A_Y$ is a maximal commutative subalgebra of $\End (Y)$;

\item[$(c)$]
 the isomorphisms
$\al$, $\beta$
identify the regular representation
of the algebra $\oplus_{s=1}^d A_{p_s, \Phi}$ and the $A_Y$-module $Y$;
the isomorphism
$\al$
identifies the contravariant form on $Y$ and the residue form $(\,,\,)
=\oplus_{s=1}^d (\,,\,)_{p_s}$ on  $\oplus_{s=1}^d A_{p_s, \Phi}$ multiplied by $(-1)^k$;

\item[$(d)$] in particular, if
the rank of the contravariant form $S^{(a)}\vert_{\Sing \FF^k(z^0)}$ equals the sum of Milnor numbers
 of the points $p_1,\dots,p_d$,
then $Y=\sing W(z^0)$ and the module $\sing W(z^0)$ over the unital subalgebra
of $\End (\sing W(z^0))$ generated by geometric Hamiltonians
${\rm pr}_{\tilde Y}K^1_j(z^0)^*\vert_{\sing W(z^0)}$, $j\in J$,
 is isomorphic to the regular
representation of the algebra $\oplus_{s=1}^dA_{p_s, \Phi}$;
\end{enumerate}
\end{enumerate}
\end{thm}

The theorem follows from the corresponding statements of Section~\ref{An abstract setting, 2}.

\begin{thm}
\label{thm ham of unb arr}
Assume that the arrangement $(\A(z^0),a)$ is unbalanced.
\begin{enumerate}\itemsep=0pt
\item[$(i)$]
Then the contravariant form is nondegenerate on $\sing W(z^0)$.

\item[$(ii)$]
Let ${\rm pr}_{\sing W(z^0)} : \sing W\to \sing W(z^0)$ be the orthogonal projection, 
$A_{\sing W(z^0)}$ be the unital subalgebra of $\End (\sing W(z^0))$
 generated by the operators
 ${\rm pr}_{\sing W(z^0)}K^1_j(z^0)^*|_{\sing W(z^0)}$, $j\in J$.
Then $A_{\sing W(z^0)}$ is commutative and its elements are symmetric with respect to the contravariant form on $\sing W(z^0)$.

\item[$(iii)$]

Let $p_1,{\dots},p_d$ be a list of all
distinct isolated critical points  of $\Phi(z^0,\,\cdot\,)$.
Then the $A_{\sing W(z^0)}\!$-module $\sing W(z^0)$  is isomorphic to the regular
representation of the algebra $\oplus_{s=1}^dA_{p_s, \Phi}$.
\end{enumerate}
\end{thm}

\begin{proof} If $(\A(z^0),a)$ is unbalanced, then the sum of Milnor numbers of the master function equals
$|\chi(U(\A(z^0)))|$ and equals $\dim  \Sing \FF^k(\A(z^0))$. This implies part $(i)$ of the theorem. Parts $(ii)$ and $(iii)$ follow from Theorem~\ref{thm BAD}.
\end{proof}

\subsection{Remark on critical points of real arrangements}
\label{Remarks on critical point of real arrangements}

Assume that $(g_j)_{j\in J}$ are real, see Remark \ref{rem g real}. Assume that $z^0\in \R^n\subset \C^n$.
Assume that the contravariant form is positive def\/inite on $\sing W(z^0)$. Let ${\rm pr}_{\sing W(z^0)}: \sing
W\to \sing W(z^0)$
be the orthogonal projection. Assume that the operators
\[
{\rm pr}_{\sing W(z^0)}K^1_j\big(z^0\big)^*|_{\sing W(z^0)} : \ \sing W\big(z^0\big) \to \sing W\big(z^0\big),
\qquad
j\in J,
\]
commute and are symmetric with respect to the contravariant form.

\begin{thm}
\label{thm on real points} Under these assumptions, any
 critical point  $p$ of the master function $\Phi(z^0,\cdot) : U(\A(z^0))\to U(\A(z^0))$,
is nondegenerate and the nonzero value $[\bs E](p)$ at $p$ of the canonical element
lies in the real part $\sing W(z^0)_\R$ of $\sing  W(z^0)$ $($up to multiplication by a nonzero complex number$)$.
\end{thm}

\begin{proof}
On one hand,
under assumptions of the theorem all the linear operators preserve $\sing W(z^0)_\R$ and can be diagonalized simultaneously.
 That means that any element of the algebra of geometric Hamiltonians $A_{\sing W(z^0)}$
(generated by operators
${\rm pr}_{\sing W(z^0)}K^1_j(z^0)^*|_{\sing W(z^0)}$, $j\in J$) is diagonalizable.
 On the other hand, if the Milnor number of $p$ is greater than one, then the local algebra
$\AT$ has nilpotent elements and, by Theorem~\ref{thm geom interpr bad fibers 2}, the algebra
 $A_{\sing W(z^0)}$  has nondiagonalizable elements.

The second part of the theorem is clear.
\end{proof}

 Theorem \ref{thm on real points} is in the spirit of the main theorem of \cite{MTV7} and Conjecture 5.1 in
\cite{MTV6}.  The main theorem of \cite{MTV7} says that certain Schubert cycles intersect transversally and all intersection
point are real. These two statements correspond to the two statements of Theorem \ref{thm on real points}.

\section{Arrangements with symmetries}
\label{Arrangements with symmetries sec}

\subsection{A family of prediscriminantal arrangements}
\label{sec prediscriminantal}

In this section we consider a special family of parallelly  translated hyperplanes, see
Section \ref{sec par trans}. The members of that family will be called prediscriminantal arrangements.

\begin{data} \label{ii} \rm
Let $\h^*$ be a complex vector space of dimension $r$ with a collection of
vectors $\al_1,\dots,\al_r, \La_1,\dots,\La_N \in \h^*$ and a symmetric bilinear form $(\,,\,)$.
We assume  that
$(\al_i,\al_i) \neq 0$ for  every $i=1,\dots,r$.

Let $\bs k=(k_1,\dots,k_r)$ be a collection of nonnegative integers. We denote
$k=\sum_ik_i$ and assume that $k>0$.  We assume that for every $b=1,\dots,N$ there exists $i$ such
that $(\al_i,\La_b)\neq 0$ and $k_i>0$.
\end{data}

Consider the expressions{\samepage
\begin{gather}
\label{expressions}
f_{(i),l,l'} = t^{(i)}_{l} - t^{(i)}_{l'} + z_{(i),l,l'}
\qquad
\text{such that}\ i=1, \dots,r\ {}
\text{and } \   1\leq l<l'\leq k_i ;
 \\
f_{(i,i'),l,l'} = t^{(i)}_{l} - t^{(i')}_{l'} + z_{(i,i'),l,l'}
\notag\\
\phantom{f_{(i,i'),l,l'} =}{}
\text{such that}\ 1 \leq i< i'\leq r,\
1\leq l\leq k_i , 1\leq l'\leq k_{i'} \
\text{and} \    (\al_i,\al_{i'})\neq 0;
\notag
\\
f_{(i,b),l}  =  -t^{(i)}_{l} + z_{(i,b),l}
\qquad
\text{such that}\ 1 \leq i \leq r,\ 1\leq l\leq k_{i},\
1\leq b\leq N \
\text{and} \ (\al_i,\La_{b})\neq 0.
\notag
\end{gather}}

\noindent
Let~$J$ denote the set of all low indices of the letters $f$ in these expressions.
 So~$J$ is the union of three nonintersecting subsets~$J_1$, $J_2$, $J_3$ where
  $J_1$ consists of triples $\{(i),l,l'\}$ from the f\/irst line of~\eqref{expressions},
  $J_2$ consists of four-tuples $\{(i,i'),l,l'\}$ from the second line,
  $J_3$ consists of triples $\{(i,b),l\}$ from the third line.
 Let~$n$ be the number of elements in~$J$.

Consider $\C^k$ with coordinates
\[
 t = \big(t_1^{(1)}, \dots , t_{k_1}^{(1)}, \dots, t_{1}^{(r)}, \dots, t_{k_r}^{(r)}\big).
\]
Consider $\C^n$ with coordinates $z= (z_j)_{j\in J}$ and $\C^n\times \C^k$ with coordinates $z$, $t$.
For any $j\in J$ the expression $f_j$ can be considered as a linear function on $ \C^n\times \C^k$.
We have $f_j=z_j + g_j$ where $g_j = t^{(i)}_{l} - t^{(i)}_{l'}$ if $j=\{(i),l,l'\}$, where
$g_j = t^{(i)}_{l} - t^{(i')}_{l'}$ if $j=\{(i,i'),l,l'\}$ and
$g_j = -t^{(i)}_{l}$ if $j=\{(i,b),l\}$.
The functions $g_j$, $j\in J$, can be considered as linear functions on~$\C^k$.

For $j\in J$, the equation $f_j(z,t)=0$ def\/ines a hyperplane $\tilde H_j \subset\C^n\times\C^k$ and
 we get an arrangement   $\tilde \A=\{\tilde H_j\ | \ j\in J\}$ in $\C^n\times\C^k$.

We assign (nonzero) weights $a_j$ to  hyperplanes of $\tilde C$ by putting
\begin{gather}
\label{weights discrim}
a_{(i),l,l'} = (\al_i,\al_i),
\qquad
a_{(i,i'),l,l'} = (\al_i,\al_{i'}),
\qquad
a_{(i,b),l} = -(\al_i,\La_{b}).
\end{gather}

The weighted arrangement $\tilde C$ is an example of a family of parallelly translated hyperplanes
considered in Sections \ref{sec par trans}--\ref{More on Hamiltonians of bad fibers}.

\subsection{Discriminantal arrangements}
\label{sec Discriminantal arrangements}

Let $X\subset \C^n$ be the subset def\/ined by the following equations:
\begin{gather*}
z_{(i),l,l'}  =  0 , \qquad
 i=1, \dots,r \
\text{and} \   1\leq l<l'\leq k_i ;
\notag
\\
z_{(i,i'),l,l'}  =  0,
\qquad
 1 \leq i< i'\leq r,\
1\leq l\leq k_i ,\ 1\leq l'\leq k_{i'} \
\text{and} \   (\al_i,\al_{i'})\neq 0;
\\
z_{(i,b),l}  =  z_{(i',b),l'},\qquad
  1 \leq i \leq r,\ 1\leq l\leq k_{i},\
 1 \leq i' \leq r,\ 1\leq l'\leq k_{i'},\
1\leq b\leq N ,
\notag
\\
 \phantom{z_{(i,b),l}  =  z_{(i',b),l'},} \qquad
 (\al_i,\La_{b})\neq 0,\ (\al_{i'},\La_{b})\neq 0.\notag
\end{gather*}
The subset $X$ is an $N$-dimensional af\/f\/ine space. We will use the
following coordinates
$x_1,\dots,x_N$ on $X$ def\/ined by the equations
 $x_b = z_{(i,b),l}$, where $1\leq b\leq N$, $1\leq i\leq r$, $1\leq l \leq k_i$ and $(\al_i,\La_b)\neq 0$.
 We will be interested in the open subset $U(X)\subset X$,
\[
U(X) =\{ z\in X\ |\ x_1(z),\dots,x_N(z)\ \text{are all distinct}\} .
\]
Let us consider the arrangement $\A(z)=(H_j(z))_{j\in J}$ in the f\/iber $\C^k$ of the projection
$\C^n\times \C^k \to \C^n$ over a point $z\in U(X)$ with coordinates $x_1(z),\dots,x_N(z)$.
Its hyperplanes are def\/ined by the equations:
\begin{gather*}
t^{(i)}_{l} - t^{(i)}_{l'}  =  0,
\qquad i=1, \dots,r\
\text{and}\   1\leq l<l'\leq k_i ;
\notag
\\
t^{(i)}_{l} - t^{(i')}_{l'}  =  0 ,
\qquad
 1 \leq i< i'\leq r,\
1\leq l\leq k_i , \ 1\leq l'\leq k_{i'} \
\text{and} \   (\al_i,\al_{i'})\neq 0;
\\
t^{(i)}_{l} - x_{b}(z)  =  0 ,
\qquad
1 \leq i \leq r,\ 1\leq l\leq k_{i},\
1\leq b\leq N \
\text{and} \ (\al_i,\La_{b})\neq 0.
\notag
\end{gather*}
The weights of these hyperplanes are def\/ined by formula \eqref{weights discrim}.
This weighted arrangement is called discriminantal, see \cite{SV, V1}.

\subsection{Symmetries of the family of prediscriminantal arrangements}
\label{Symmetries of the family of prediscriminantal arrangements}

The product of symmetric groups $S_{\bs k}=S_{k_1}\times \dots \times S_{k_r}$ acts on $\C^k$ by permuting
coordinates $t^{(i)}_l$ with the same upper index.
More precisely, a point $p\in \C^k$ with coordinates
$(t_1^{(1)}(p)$, \dots , $t_{k_1}^{(1)}(p)$, \dots, $t_{1}^{(r)}(p)$, \dots, $t_{k_r}^{(r)}(p))$
is mapped by an element $\sigma = (\sigma_1,\dots,\sigma_r)$ to the point with coordinates
$(t_{\sigma_a^{-1}(1)}^{(1)}(p), \dots , t_{\sigma_1^{-1}(k_1)}^{(1)}(p),
\dots, t_{\sigma_r^{-1}(1)}^{(r)}(p), \dots, t_{\sigma_r^{-1}(k_r)}^{(r)}(p))$.

The group $S_{\bs k}$ acts also on $\C^n$. Namely, an
 element $\sigma = (\sigma_1,\dots,\sigma_r)$ sends a coordina\-te~$z_{(i),l,l'}$ to the coordinate
$z_{(i),\sigma_i (l),\sigma_i(l')}$ if
$\sigma_i (l)<\sigma_i(l')$ and to
$-z_{(i),\sigma_i (l),\sigma_i(l')}$ if
$\sigma_i (l)>\sigma_i(l')$.
An element~$\sigma$ sends a coordinate
$z_{(i,i'),l,l'}$ to
$z_{(i,i'),\sigma_i (l),\sigma_{i'}(l')}$.
An element $ \sigma$ sends a coordi\-na\-te~$z_{(i,b),l}$ to the coordinate
$z_{(i,b),\sigma_i (l)}$.

Clearly every point of $X\subset \C^n$ is a f\/ixed point of the $S_{\bs k}$-action.

The actions of $S_{\bs k}$ on $\C^n$ and $\C^k$ induce an action on
$\C^n\times \C^k$.
The $S_{\bs k}$-action on $\C^n\times \C^k$ preserves the arrangement $\tilde \A$ and sends f\/ibers of the projection
$\C^n\times \C^k \to\C^n$ to f\/ibers.

The $S_{\bs k}$-action on  $\tilde \A$ corresponds to an
$S_{\bs k}$-action on the set  $J=J_1\cup J_2\cup J_3$. That action preserves the summands. For
$\sigma = (\sigma_1,\dots,\sigma_r)\in S_{\bs k}$
we have
\begin{gather*}
  \{(i),l,l'\}\mapsto \{(i),\min(\sigma_i (l),\sigma_i(l')), \max(\sigma_i (l),\sigma_i(l'))\},
\\
 \{(i,i'),l,l'\}\mapsto \{(i,i'),\sigma_i (l),\sigma_{i'}(l')\},
\qquad
\{(i,b),l\}
\mapsto \{(i,b),\sigma_i (l)\} ,
\end{gather*}
 where
$\{(i),l,l'\}\in J_1$,  $\{(i,i'),l,l'\}\in J_2$,  $\{(i,b),l\}\in J_3$.

Consider the discriminant $\Delta = \cup_{C\in \frak C} H_C  \subset   \C^n$,
see Section~\ref{Discr}. Here $\frak C$ is the
set of all circuits of the matroid
of the collection $(g_j)_{j\in J}$.
Clearly the $S_{\bs k}$-action on $\C^n$ preserves the
discriminant and permutes the hyperplanes $\{H_C\ |\ C\in \frak C\}$.

The action of $S_{\bs k}$ on the hyperplanes of the discriminant corresponds to the
following action  on  $\frak C$.
If $C=\{j_1,\dots, j_l\}\subset J$ is a circuit
and $\sigma\in S_{\bs k}$, then $\sigma (C)$ is the circuit $ \{\sigma (j_1),\dots, \sigma(j_l)\}$.

\subsection[The $S_{k}$-action on geometric Hamiltonians]{The $\boldsymbol{S_{k}}$-action on geometric Hamiltonians}
\label{The S-{ k}-action on Hamiltonians}

We use notations of Section \ref{sec Good fibers} and for $z\!\in \!\C^n-\Delta$ denote
$V\!=\FF^k(\A(z))$, $\Sing V\!=\Sing\FF^k(\A(z))$.
The triple $(V, \Sing V, S^{(a)})$ does not depend on
 $z\in\C^n-\Delta$ as explained in Section \ref{sec Good fibers}.

Fix an order on the set $J$. Recall that the standard basis of
$V^* = \OS^k(\A(z))$, associated with an order on $J$, is formed by elements
$(H_{j_1},\dots,H_{j_k})$ where
 $\{{j_1} <\dots <{j_k}\}$  runs through the set of all independent ordered $k$-element subsets of
$J$. The (dual) standard basis of $V$ is formed by the corresponding vectors
$F(H_{j_1}(z),\dots,H_{j_k}(z))$.
We have $F(H_{j_1}(z),\dots,H_{j_k}(z)) = (-1)^{|\mu|}F(H_{j_{\mu(1)}}(z),\dots,H_{j_{\mu(k)}}(z))$ for any
$\mu \in S_k$, see Section \ref{An arrangement with normal crossings only}.

The
$S_{\bs k}$-action on $J$ induces an action on $V$ and $V^*$. For $\sigma\in S_{\bs k}$, we have
\[
\sigma : \ F(H_{j_1}(z),\dots,H_{j_k}(z)) \mapsto
F(H_{\sigma(j_1)}(z),\dots,H_{\sigma(j_k)}(z))
\]
 and
$(H_{j_1}(z),\dots,H_{j_k}(z))$
$ \mapsto $
$(H_{\sigma(j_1)}(z),\dots,H_{\sigma(j_k)}(z))$.
The $S_{\bs k}$-action
on $V$ preserves the subspace $\Sing V$ and preserves the contravariant form $S^{(a)}$ on $V$.

Let $z^0\in U(X)$. Then $z^0$ is $S_{\bs k}$-invariant and the group  $S_{\bs k}$ acts on the f\/iber
$\C^k$ over $z^0$ and  on the weighted arrangement $(\A(z^0),a)$ in that f\/iber.
The subspaces $\FF^k(\A(z^0))$, $\sing \FF^k(\A(z^0))$ of $V$  are $S_{\bs k}$-invariant.

Let $V\to V^*$ be the map associated with the contravariant form. Let $W$, $\sing W$, $W(z^0)$,
$\sing W(z^0)$ be the images of $V$, $\sing V$, $\FF^k(\A(z^0))$, $\sing \FF^k(\A(z^0))$, respectively.
All these subspaces are $\Sk$-invariant.

An $S_{\bs k}$-action on a vector space def\/ines an $S_{\bs k}$-action on
linear operators on that space. For $\sigma \in S_{\bs k}$
and a linear operator $L$ we def\/ine $\sigma(L) = \sigma L \sigma^{-1}$.

In Section \ref{sec key identity} we have def\/ined operators
$L_C : V\to V$, $C\in \frak C$.
Clearly for any $\sigma \in S_{\bs k}$ and any $C\in \frak C$ we have $\sigma(L_C)=L_{\sigma(C)}$.

In Section \ref{sec key identity} we have considered dif\/ferential 1-forms on $\C^n\times \C^k$ which
were denoted by~$\omega_j$, $j\in J$, and $\omega_C$, $C\in \frak C$.
The $S_{\bs k}$-action on
$\C^n\times\C^k$ preserves this set of dif\/ferential 1-forms. Namely for any $j\in J$,
$ C\in \frak C$,  $\sigma \in S_{\bs k}$,
we have $\sigma : \omega_j\mapsto \omega_{\sigma(j)}$, $\omega_C \mapsto \omega_{\sigma(C)}$.

\begin{Lemma}
\label{lem on invariantce}
The following objects are $S_{\bs k}$-invariant:
\[
\sum_{j\in J} a_j \omega_j ,
\qquad
\sum_{C\in \frak C}
\omega_{C} \otimes L_C ,
\qquad
\sum_{{\rm independent } \atop \{j_1 < \dots < j_k\} \subset J }
 \omega_{j_1} \wedge \dots \wedge \omega_{j_k}
\otimes
F(H_{j_1}, \dots , H_{j_k}) .
\]
\end{Lemma}

By the def\/inition of    $K_j(z) :  V\to  V$,  $j\in J$, we have
\[
\sum_{C\in \frak C}
\omega_{C} \otimes L_C = \sum_{j\in J} dz_j \otimes K_j(z) ,
\]
see formula \eqref{LK}. The functions $K_j(z)$  are $\End(V)$-valued meromorphic functions on~$\C^n$.

Since $S_{\bs k}$ acts on $\C^n$ and $\End(V)$ it also acts on $\End(V)$-valued functions on
$\C^n$,   $\sigma : F(q) \mapsto \sigma F(\sigma^{-1}(q))\sigma^{-1}$ for $q\in\C^n$.
Lemma~\ref{lem on invariantce} allows us to describe the $S_{\bs k}$-action on functions~$K_j(z)$.

\begin{cor}
\label{cor action on H's}
An element $\sigma=(\sigma_1,\dots,\sigma_r)\in S_{\bs k}$ acts on functions  $K_j$, $j\in J$,
by the formulas:
\begin{gather*}
K_{(i),l,l'}(q)  \mapsto
K_{(i),\sigma_i (l),\sigma_i(l')}\big(\sigma^{-1}(q)\big),
\qquad
\text{if} \ \sigma_i (l)<\sigma_i(l'),
\notag
\\
K_{(i),l,l'}(q)  \mapsto  -   K_{(i),\sigma_i (l'),\sigma_i(l)}\big(\sigma^{-1}(q)\big),
\qquad
\text{if} \ \sigma_i (l)>\sigma_i(l'),
\notag
\\
K_{(i,i'),l,l'}(q)  \mapsto
 K_{(i),\sigma_i (l),\sigma_{i'}(l')}\big(\sigma^{-1}(q)\big) ,
\notag
\\
K_{(i,b),l}(q)  \mapsto
K_{(i,b),\sigma_i (l)}\big(\sigma^{-1}(q)\big) .
\end{gather*}
\end{cor}

Let $z^0\in U(X)$.
Recall that
$\frak C_0 = \{ C\in \frak C\ | \ z^0\in H_C\}$,
\begin{gather*}
K_j^0(z)  =  \sum_{C\in \frak C_0}
 \frac{\la_j^C}{f_C(z)}   L_C  ,
\qquad
K_j^1(z) = K_j(z) - K_j^0(z) .
\end{gather*}

\begin{cor}
\label{cor on K1j}
An element $\sigma\in \Sk$ acts on operators $K^1_j(z^0) :V\to V$, $j\in J$, by the formula
\begin{gather*}
\sigma\big(K^1_j\big(z^0\big)\big) = \pm K^1_{\sigma(j)}\big(z^0\big) ,
\end{gather*}
where the minus sign is chosen only if $j=\{(i),l,l'\}$ and $\sigma_i(l)>\sigma(l')$.
\end{cor}

\subsection[Functions  $K_{\der_{x_b}}(z)$]{Functions  $\boldsymbol{K_{\der_{x_b}}(z)}$}
\label{Symmetric naive Hamiltonians}

Let $z^0 \in U(X)$.
Recall that
$
T_{z^0} = \{
\xi = \sum_{j\in J} \xi_j\frac{\der}{\der z_j}\ |\ \xi_j\in \C, \
\xi (f_C) = 0 \ \text{for all}\ C\in \frak C_0\}$,
see Section \ref{sec general remarks}.
 Def\/ine the following constant vector f\/ields on $\C^n$,
\[
\der _{x_b} =
\sum _{1\leq i \leq r,\ 1\leq l\leq k_i,\ (\al_i,\La_b)\neq 0}
\frac{\der\phantom {aaa}}{\der z_{(i,b),l}}   ,
\qquad
b=1,\dots,N .
\]

\begin{Lemma}
\label{lem naiv symm hams}
The vector fields $\der _{x_b}$, $b=1,\dots,N$, are elements of\ {} $T_{z^0}$. The
$\End (V)$-valued functions
\begin{gather*}
K_{\der _{x_b}}(z) = \sum _{1\leq i \leq r,\ 1\leq l\leq k_i,\ (\al_i,\La_b)\neq 0}
K_{(i,b),l}(z) ,
\qquad
b=1,\dots,N,
\end{gather*}
are $S_{\bs k}$-invariant.
\end{Lemma}

By Lemma \ref{lem regularity of KZ}, the functions
$K_{\der _{y_b}}(z)$ are regular at $z^0$ and their dual operators preserve the subspace $\sing W(z^0) \subset V^*$.

\begin{cor}
For $z^0\in U(X)$, the operators
\begin{gather*}
K_{\der _{x_b}}\big(z^0\big)^*|_{\sing W(z^0)} : \ \sing W\big(z^0\big) \to \sing W\big(z^0\big) ,
\qquad
b=1,\dots,N ,
\end{gather*}
are $S_{\bs k}$-invariant.
\end{cor}

The operators $K_{\der _{x_b}}(z^0)^*|_{\sing W(z^0)}$
are naive  geometric Hamiltonians on $\sing W(z^0)$
in the sense of Section \ref{sec general remarks}.
They commute and they are symmetric operators with respect to the contravariant form.

\subsection[Naive geometric Hamiltonians on $\sing W^-(z^0)$]{Naive geometric Hamiltonians on $\boldsymbol{\sing W^-(z^0)}$}
\label{sec skew of W}

The space $W$ has the canonical direct sum decomposition into isotypical components correspon\-ding to
irreducible representations of $S_{\bs k}$. One of the isotypical components is
the component
\[
 W^- = \big\{ x\in W\ | \sigma (x) = (-1)^{|\sigma|} x\ \text{for any}\ \sigma\in S_{\bs k}\big\} ,
\]
corresponding to the alternating representation.
If $L :W\to W$ is an $S_{\bs k}$-invariant linear operator, then $L$ preserves the canonical
decomposition and, in particular, it preserves the subspace $W^-$.

Let $z^0\in U(X)$.
Then the subspace $\sing W(z^0)\subset W$ is an $\Sk$-submodule. We def\/ine
$ \sing W^-(z^0) = W^-\cap \sing W(z^0)$.
For any $b=1,\dots,N$, the operator $K_{\der_{x_b}}(z^0)^*$ preserves
$\sing W^-(z^0)$.
The operators
\[
K_{\der _{x_b}}\big(z^0\big)^*|_{\sing W^-(z^0)} : \ \sing W^-\big(z^0\big)\to \sing W^-\big(z^0\big)
\]
will be called {\it
naive geometric Hamiltonians} on $\sing W^-(z^0)$.

\subsection[$S_k$-symmetries of the canonical element]{$\boldsymbol{S_k}$-symmetries of the canonical element}
\label{sec Sk prop of speci element}

Let $z^0\in U(X)$.
Let $p \in U(\A(z^0))$ be an isolated critical point of the
master function $\Phi(z^0,\,\cdot\,)
: U(\A(z^0)) \to \C$. Let $O(p) = \{\sigma(p)\ | \ \sigma\in \Sk\}$ be the $\Sk$-orbit of~$p$.
The orbit consists of $k_1!\cdots k_r!$ points.
Let
$A_{\sigma (p),\Phi}$
be the local algebra of the critical
point $\sigma(p)$
and $(\,,\,)_{\sigma(p)}$
the residue bilinear form on $A_{\sigma(p),\Phi}$.
Let
$[ \ ]_{\sigma(p)} : \C(t)_{U(\A(z^0))} \to A_{\sigma(p),\Phi}$ be the canonical projection
and $[\bs E]_{\sigma(p)} \in A_{\sigma(p),\Phi}\otimes \sing W(z^0)$ the projection of the
the canonical element.

The group $\Sk$ acts on functions on $U(\A(z^0))$. If $g\in \C(t)_{U(\A(z^0))}$ and $\sigma\in \Sk$, then
$\sigma (g)(q) = g(\sigma^{-1}(q))$ for $q\in U(\A(z^0))$.

\begin{Lemma}
\label{lem Sk action on fj}
An element $\sigma\in \Sk$ acts on functions $f_j$, $j\in J$, by the formula
\begin{gather*}
\sigma(f_j) = \pm f_{\sigma(j)} ,
\end{gather*}
where the minus sign is chosen only if $j=\{(i),l,l'\}$ and $\sigma_i(l)>\sigma(l')$.
\end{Lemma}

The $\Sk$-action on $\C(t)_{U(\A(z^0))}$
induces an isomorphism $\sigma : A_{p,\Phi} \to A_{\sigma (p),\Phi}$, $[g]_p \mapsto [\sigma(g)]_{\sigma(p)}$.
Let us compare the   residue bilinear forms on $A_{p,\Phi}$
and  $A_{\sigma(p),\Phi}$ and projections of the canonical element to
$A_{p,\Phi}$
and  $A_{\sigma(p),\Phi}$.

\begin{Lemma}
\label{lem S and (,)}
For $f,g\in  A_{p,\Phi}$, we have
\[
(\sigma(f),\sigma(g))_{\sigma(p)}
=
(-1)^\sigma(f,g)_{p} \qquad \mbox{and} \qquad
 \sigma([\bs E]_{p}) = (-1)^\sigma [\bs E]_{\sigma(p)},
\]
where $(-1)^\sigma = \prod\limits_{i=1}^r(-1)^{\sigma_i}$.
\end{Lemma}

For $\sigma\in\Sk$,
let $\al_{\sigma(p)} : A_{\sigma(p),\Phi} \to \sing W(z^0)$ be the linear monomorphisms constructed  in Section~\ref{Algebra of bad 2}. Let $Y_{\sigma(p)}$ be the image of $\al_{\sigma(p)}$.
If $\mu$ is the Milnor number of $p$, then $\dim Y_{\sigma(p)}=\mu$.

\begin{Lemma}
\label{lem S and al}
For $f \in  A_{p,\Phi}$, we have
$\sigma(\al_p(f)) = \al_{\sigma(p)}(\sigma(f))$.
\end{Lemma}

\begin{cor}
\label{cor Y-p and S}
We have $\sigma(Y_p) = Y_{\sigma(p)}$.
\end{cor}

By Theorem \ref{thm many point bad 2}, the subspaces $Y_{\sigma(p)}$ are all orthogonal and the contravariant form
on
$Y_{O(p)} =\oplus_{\sigma\in \Sk} Y_{\sigma (p)}$ is
nondegenerate.
The group $\Sk$ acts on $Y_{O(p)}$.
Denote
\[
Y_{O(p)}^- = \big\{ v\in Y_{O(p)}\ |\ \sigma(v) = (-1)^\sigma v\
\text{for all} \ \sigma\in \Sk\big\} .
\]

\begin{cor}
\label{cor Y--}
We have
$Y_{O(p)}^- = \{ \sum_{\sigma \in \Sk} (-1)^\sigma \sigma(v)\ | \ v\in Y_p\}$
and $\dim Y_{O(p)}^- = \mu$.
\end{cor}

Let ${\rm Ant} = \sum_{\sigma\in \Sk} (-1)^\sigma \sigma : \sing W(z^0)\to \sing W^-(z^0)$ be the anti-symmetrization operator.
Denote
\[
[\bs E]^-_p = (1\otimes {\rm Ant}) [\bs E]_p .
\]
 Def\/ine a linear monomorphism
\[
\al_p^{-} : \ A_{p,\Phi} \to \sing W^-\big(z^0\big), \qquad
f \to (f,[\bs E]^-_p)_p\ .
\]

\begin{Lemma}
\label{lem gamma}
We have
\[
Y_{O(p)}^- = \{ \al_p^-(f)\ | \ f\in A_{p,\Phi}\}.
\]
The linear map $\al_p^-$ identifies the contravariant form on $Y_{O(p)}^-$ and the residue bilinear form on~$A_{p,\Phi}$ multiplied by $(-1)^k k_1!\cdots k_r!$. In particular, the contravariant form on
$Y_{O(p)}^-$ is nondegenerate.
We also have
\begin{gather}
\label{lenth -}
S^{(a)}([\bs E])^-_p,[\bs E]^-_p)  =  (-1)^k k_1!\cdots k_r!\big[\Hess \Phi \big(z^0,\,\cdot\,\big)\big]_p .
\end{gather}
\end{Lemma}

Let $\tilde Y \subset \sing W$ be a vector subspace such that
\begin{enumerate}\itemsep=0pt
\item[$(a)$]
 $Y_{O(p)}\subset \tilde Y$;
\item[$(b)$]
 the contravariant  form restricted on $\tilde Y$ is nondegenerate;
 \item[$(c)$] for any $j\in J$,
 the subspace
$\tilde Y$ lies in the kernel of $L_C^*$, $C\in \frak C_0$;
\item[$(d)$] the subspace $\tilde Y$ is $\Sk$-invariant.
\end{enumerate}
For example, we can choose $\tilde Y=Y_{O(p)}$.
Let ${\rm pr}_{\tilde Y} :\sing W\to \tilde Y$ be the orthogonal projection with respect to the contravariant form.

For every $j\in J$ and $\sigma\in \Sk$, the map  ${\rm pr}_{\tilde Y}K^1_j(z^0)^*: \sing W\to \sing W$ preserves $Y_{\sigma(p)}$
and ${\rm pr}_{\tilde Y}K^1_j(z^0)^*|_{\tilde Y}$ does not depend on the choice of $\tilde Y$, see  Theorem~\ref{thm BAD}.
We also have $\sigma : {\rm pr}_{\tilde Y}K^1_j(z^0)^* \mapsto
\pm {\rm pr}_{\tilde Y}K^1_{\sigma(j)}(z^0)^*$ where the minus sign is chosen only if $j=\{(i),l,l'\}$ and \mbox{$\sigma_i(l)>\sigma(l')$}.

\medskip

\noindent
{\bf Important def\/inition.}
Denote by $P^\Sk$ the algebra of polynomials with complex coef\/f\/icients in~$n$ variables $a_j/f_j$, $j\in J$, such that for any
$F(a_j/f_j,j\in J)\in P^\Sk$ the function
$F(a_j/f_j(z^0,\cdot)$, $j\in J)\in U(\A(z^0))$ is $\Sk$-invariant.

\begin{Lemma}
\label{PiSk iso}
The natural homomorphism $P^\Sk \to \AT$ is an epimorphism.
\end{Lemma}

Let $F(a_j/f_j(z^0,\cdot),j\in J)\in P^\Sk$.
Replace in $F$ each variable $a_j/f_j$ with the operator ${\rm pr}_{\tilde Y}K^1_j(z^0)^*$. Denote the resulting operator
on $\sing W$  by $F({\rm pr}_{\tilde Y}K^1_j(z^0)^*, j\in J)$.
This operator  preserves $Y_{\sigma(p)}$ for any $\sigma\in\Sk$ and its restriction to $Y_{\sigma(p)}$
does not depend on the choice of $\tilde Y$.
The operator $F({\rm pr}_{\tilde Y}K^1_j(z^0)^*, j\in J)$ is $\Sk$-invariant by Corollary \ref{cor on K1j}.
Hence, $F({\rm pr}_{\tilde Y}K^1_j(z^0)^*, j\in J)$ preserves $Y^-_{O(p)}$.

\begin{thm}
\label{thm F's}
\quad
\begin{enumerate}\itemsep=0pt
\item[$(i)$]
For $F\in P^\Sk$, the operators
\begin{gather*}
F\big({\rm pr}_{\tilde Y}K^1_j\big(z^0\big)^*, j\in J\big)|_{Y^-_{O(p)}}   :\ Y^-_{O(p)}   \to  Y^-_{O(p)}
\end{gather*}
commute and are symmetric with respect to the contravariant form.
\item[$(ii)$]
The map $F\mapsto F({\rm pr}_{\tilde Y}K^1_j(z^0)^*, j\in J)|_{Y^-_{O(p)}}$ induces an algebra monomorphism
\begin{gather*}
\beta_p^-: \ A_{p,\Phi} \to \End(Y^-_{O(p)}) .
\end{gather*}
\item[$(iii)$]
The image of this monomorphism, denoted by $A_{Y^-_{O(p)}}$, is a maximal commutative
subalgebra of\ $\End(Y^-_{O(p)})$.
\item[$(iv)$]
The naive geometric Hamiltonians $K_{\der_{x_b}}(z^0)^*|_{Y_{O(p)}^-}$,
$b=1,\dots,N$, are elements of  $A_{Y^-_{O(p)}}$.
\item[$(v)$]
The maps $\al_p^-$, $\beta_p^-$ define an isomorphism of the regular representation of $A_{p,\Phi}$ and
the $A_{Y^-_{O(p)}}$-module $Y^-_{O(p)}$.
The linear map $\al_p^-$ identifies the contravariant form on $Y_{O(p)}^-$ and the residue bilinear form on
$A_{p,\Phi}$ multiplied by $(-1)^k k_1!\cdots k_r!$.

\item[$(vi)$]
Define the value of $[\bs E]^-_p$ at $p$ as the image of the natural projection of
$[\bs E]^-_p$ to $A_{p,\Phi}/\frak m_p\otimes \sing W^-(z^0) = \sing W^-(z^0)$. Then the value  $[\bs E]^-_p(p)$  is nonzero and lies in
$Y^-_{O(p)}$.
The vector $[\bs E]_p^-(p)$ is the only $($up to proportionality$)$
common eigenvector of the operators $F({\rm pr}_{\tilde Y}K^1_j(z^0)^*, j\in J)|_{Y^-_{O(p)}}$, $F\in P^{S_{\bs k}}$,
we have
\[
F\big({\rm pr}_{\tilde Y}K^1_j\big(z^0\big)^*, j\in J\big)[\bs E]^-_p(p) = F\big(a_j/f\big(z^0,p\big), j\in J\big)[\bs E]^-_p(p) .
\]
\end{enumerate}
\end{thm}

\begin{thm}
\label{thm many point bad Sk}
Let $p_1,\dots,p_d$ be a list of
isolated critical points  of $\Phi(z^0,\,\cdot\,)$ such that
the orbits $O(p_1),\dots,O(p_d)$ do not intersect.
Let $Y_{p_s}^-=\al_{p_s}^-(A_{p_s,\Phi})\subset \sing W^-(z^0)$, $s=1,\dots,d$, be the corresponding subspaces.
Then the sum of these subspaces is direct and orthogonal with respect to the contravariant form.
\end{thm}

\begin{cor}
\label{cor sum of Milnor S}
The sum of Milnor numbers of the critical points $p_1,\dots,p_d$ is not greater than
the rank of the contravariant form $S^{(a)}\vert_{\sing W^-(z^0)}$.
\end{cor}

Let $\tilde Y \subset W$ be a vector subspace such that
\begin{enumerate}\itemsep=0pt
\item[$(a)$]
 $\oplus_{s=1}^d Y_{O(p_s)}\subset \tilde Y$;
\item[$(b)$]
 the contravariant  form restricted on $\tilde Y$ is nondegenerate;
 \item[$(c)$] for any $j\in J$,
 the subspace
$\tilde Y$ lies in the kernel of $ L^*_C$, $C\in \frak C_0$;
\item[$(d)$] the subspace $\tilde Y$ is $\Sk$-invariant.
\end{enumerate}
For example, we can choose $\tilde Y=\oplus_{s=1}^d Y_{O(p_s)}$.
Let ${\rm pr}_{\tilde Y} : \sing W\to \tilde Y$ be the orthogonal projection.

Denote $Y^- = \oplus_{s=1}^d Y_{O(p_s)}^-$.
\begin{thm}
\label{thm BAD S}
For $F\in P^\Sk$, the operators
\begin{gather*}
F\big({\rm pr}_{\tilde Y}K^1_j\big(z^0\big)^*, j\in J\big)|_{Y^-}   :\ Y^-  \to  Y^-
\end{gather*}
do not depend on the choice of $\tilde Y$. They
commute and are symmetric with respect to the contravariant form.

Denote by
$A_{Y^-}$ the unital subalgebra of $\End(Y^-)$ generated by
$F({\rm pr}_{\tilde Y}K^1_j(z^0)^*, j\in J)|_{Y^-}$, $F\in P^\Sk$.
The naive geometric Hamiltonians $K_{\der_{x_b}}(z^0)^*|_{Y^-}$,
$b=1,\dots,N$, are elements of the algebra $A_{Y^-}$.

Consider the isomorphisms
\begin{gather*}
\al^-  = \oplus_{s=1}^d\al_{p_s}^-  :\
\oplus_{s=1}^d
A_{p_s, \Phi}
\to  \oplus_{s=1}^d Y^-_{O(p_s)}  ,
\\
\beta^-  =  \oplus_{s=1}^d\beta_{p_s}^-   :\ \oplus_{s=1}^dA_{p_s, \Phi}  \to  \oplus_{s=1}^d A_{Y_{O(p_s)}^-} .
\end{gather*}
Then
\begin{enumerate}\itemsep=0pt
\item[$(a)$]
$A_{Y^-}=\oplus_{s=1}^d A_{Y^-_{O(p_s)}}$;
\item[$(b)$]
$A_{Y^-}$ is a maximal commutative subalgebra of $\End(Y^-)$;
\item[$(c)$] the isomorphisms
$\al^-$, $\beta^-$
identify the regular representation
of the algebra
 $\oplus_{s=1}^d A_{p_s, \Phi}$ and the $A_{Y^-}$-module $Y^-$;
the linear map $\al^-$ identifies the contravariant form on $Y^-$ and the residue bilinear form
$(\,,\,)=\oplus_{s=1}^d (\,,\,)_{p_s}$ on $\oplus_{s=1}^d A_{p_s,\Phi}$
multiplied by   $(-1)^k k_1!\cdots k_r!$;

\item[$(d)$] in particular, if
the rank of the contravariant form $S^{(a)}\vert_{\sing W^-(z^0)}$ equals the sum of Milnor numbers
 of the points $p_1,\dots,p_d$,
then the module $\sing W^-(z^0)=Y^-$ over the unital subalgebra
of $\End (\sing W^-(z^0))$ generated by geometric Hamiltonians
$F({\rm pr}_{\tilde Y}K^1_j(z^0)^*, j\in J)|_{Y^-}$,
$F\in P^\Sk$,
 is isomorphic to the regular
representation of the algebra $\oplus_{s=1}^dA_{p_s, \Phi}$.
\end{enumerate}
\end{thm}

\begin{cor}
\label{cor unbal equiv}
If $z^0\in U(X)$ and the arrangement $(\A(z^0),a)$ is unbalanced, then the ope\-ra\-tors
$F({\rm pr}_{\tilde Y}K^1_j(z^0)^*, j\in J)|_{\sing W^-(z^0)} : \sing W^-(z^0)\to \sing W^-(z^0)$, $
F\in P^{S_{\bs k}}$,
commute and are symmetric with respect to the contravariant form.
\end{cor}

\begin{proof}
In this case $Y^-=\sing W^-(z^0)$ and the corollary follows from Theorem~\ref{thm BAD S}.
\end{proof}

\section{Applications to the Bethe ansatz of the Gaudin model}
\label{sec Gaudin}

\subsection{Gaudin model}
\label{sec G model}

Let $\g$ be a simple Lie algebra over $\C$ with Cartan matrix
$(a_{i,j})_{i,j=1}^r$.
Let $\h \subset \g$  be a
 Cartan subalgebra.
Fix simple roots $\al_1, \dots , \al_r$ in $\h^*$
and a nondegenerate $\g$-invariant bilinear form $(\, ,\, )$ on $\g$.
 The form identif\/ies $\g$ and $\g^*$ and def\/ines a bilinear form on
 $\g^*$.
 Let $H_1, \dots , H_r   \in \h$ be the corresponding
coroots, $\langle\la , H_i\rangle = 2 (\la,\al_i)/ (\al_i,\al_i)$
 for $\la\in\h^*$.
In particular, $\langle \al_j , H_i \rangle = a_{i,j}$.

Let $E_1, \dots , E_r\, \in \n_+$, $H_1, \dots , H_r\, \in \h$, $F_1, \dots , F_r  \in \n_- $ be the Chevalley generators of $\g$,
\begin{gather*}
[E_i, F_j]  =   \delta_{i,j}  H_i ,
\qquad i, j = 1, \dots, r ,
\\
 [ h , h']   =   0 ,
\qquad
h, h' \in \h ,
\\
 [ h, E_i]  =  \langle \alpha_i, h \rangle  E_i ,
\qquad
h \in \h, \ i = 1, \dots, r ,
\\
 [ h, F_i]  =  - \langle \alpha_i, h \rangle  F_i ,
\qquad  h\in \h, \ i = 1, \dots, r ,
\end{gather*}
and   $(\mathrm{ad}\,{} E_i)^{1-a_{i,j}} E_j = 0$,
$(\mathrm{ad}\,   F_i)^{1-a_{i,j}} F_j = 0$,
for all $i\neq j$.

Let $(J_i)_{i\in I}$ be an orthonormal basis of $\g$,
$\Omega =  \sum_{i\in I} J_i\otimes J_i  \in \g \otimes \g$
the Casimir element.

For a $\g$-module $V$ and $\mu \in \h^*$
denote by $V[\mu]$ the weight subspace of $V$ of weight $\mu$ and by
$\s V[\mu]$ the subspace of singular vectors of weight $\mu$,
\[
\s V[\mu]  =  \{   v \in V\ |\ \n_+v = 0, \ hv = \langle \mu, h
\rangle v   \}   .
\]

Let $N>1$ be an integer and  $\bs \La = (\La_1, \dots , \La_N)$,
$\La_b \in \h^*$, a set of weights.
For $\mu \in \h^*$ let $V_{\mu}$
be the irreducible $\g$-module with highest weight~$\mu$.
Denote
\[
V_{\bs \La}=
V_{\La_1} \otimes \dots \otimes V_{\La_N}  .
\]
For $X \in \End (V_{\La_i})$,  denote by
\[
X^{(i)} = 1\otimes \dots \otimes 1 \otimes X \otimes 1
\otimes \dots \otimes 1
\in \End
(V_{\bs \La})
\]
the operator acting nontrivially on the $i$-th factor only. For $X = \sum_m X_m \otimes Y_m \in
\End (V_{\La_i} \otimes V_{\La_j})$,  we set
$X^{(i,j)} = \sum_m X^{(i)}_m \otimes Y^{(j)}_m  \in \End (\V)$.

Let $x^0 = (x^0_1, \dots , x^0_N)$ be a point of $\C^N$ with distinct
coordinates.  Introduce linear operators $\Kk_1(x^0), \dots , \Kk_N(x^0)$ on
$V_{\bs \La}$ by the formula
\[
 \Kk_b(x^0)  =   \sum_{c \neq b}
\frac{\Omega^{(b,c)}}{x^0_b - x^0_c}  , \qquad b = 1, \dots , N .
\]
The operators are called {\it the Gaudin Hamiltonians}. The
Hamiltonians commute,
\[
 [ \Kk_b(x^0), \Kk_c(x^0) ] = 0 \qquad \text{for all}\  b,c .
 \]
The Hamiltonians commute with the $\g$-action on
$V_{\bs \La}$. Hence they preserve the subspaces
$\Sing \V[\mu]\subset \V$ of singular vectors of a given weight~$\mu$.

Let $\tau : \g\to \g $ be the anti-involution sending
$E_i$, $H_i$, $F_i,$ to $F_i$, $H_i$, $E_i$, respectively,  for all~$i$.
Let $W$ be a highest weight $\g$-module with a
highest weight vector~$w$.
The Shapovalov form~$S$ on~$W$ is the unique
symmetric bilinear form such that
\[
S(w, w) = 1 ,
\qquad
S(gu, v) = S(u, \tau(g)v)
\]
for all $u,v \in W$ and $g \in \g$.

Fix highest weight vectors  $v_1, \dots, v_N$  of $ V_{\La_1}, \dots, V_{\La_N}$, respectively.
Def\/ine
a symmetric bilinear form on the tensor product $\V$
 by the formula
\[
S_{\bs \La}  =  S_1 \otimes \cdots \otimes S_N ,
\]
where $S_b$ is the Shapovalov form on $V_{\La_b}$.
The form $S_{\bs\La}$ is called the tensor Shapovalov form.

{\sloppy The Gaudin Hamiltonians are symmetric with respect to the tensor Shapovalov form,
$S_{\bs\La}(\Kk_b(y)u,v)= S_{\bs\La}(u,\Kk_b(y)v)$ for any $u,v\in \V$ and $b=1,\dots,N$,
see~\cite{RV}.

}

For any $\mu\in \h^*$, the Gaudin model on
$\Sing \V[\mu]\subset \V$ is the collection
\[
\big(\Sing \V[\mu]; \  S_{\bs\La}\vert_{\Sing \V[\mu]};\  \Kk_b(y)\vert_{\Sing \V[\mu]}
 : \ \Sing \V[\mu]
\to \Sing \V[\mu], \   b=1,\dots, N\big) ,
\]
 see \cite{G1, G2}.
The main problem for the Gaudin model is to f\/ind common eigenvectors and eigenvalues of the
Gaudin Hamiltonians.

\subsection{Master    function  and  weight function, \cite{SV}}
\label{master sec}

The eigenvectors of the Gaudin Hamiltonians are constructed by the
Bethe ansatz method.  We remind that construction in this section.

Fix a collection of  nonnegative integers $\bs k = (k_1, \dots , k_r)$.  Denote $k =
k_1+\dots +k_r$,
\[
\La_\infty = \sum_{b=1}^N\La_b -\sum_{i=1}^rk_i\al_i.
\]

Consider $\C^k$ with coordinates
\[
 t = \big(t_1^{(1)}, \dots , t_{k_1}^{(1)}, \dots, t_{1}^{(r)}, \dots, t_{k_r}^{(r)}\big).
\]
Def\/ine
the master function
\begin{gather}
\label{Gaudin master}
\Phi (x^0, t, \bs\La,\bs k) =
\sum_{i=1}^{r}\sum_{1\leq j<j'\leq k_i} (\al_i,\al_i) \log \big(t_j^{(i)}-t_{j'}^{(i)}\big)
\\
\phantom{\Phi (x^0, t, \bs\La,\bs k)}{} +
\sum_{1\leq i<i'\leq r}
\sum_{j=1}^{k_i}
\sum_{j'=1}^{k_{i'}}
(\al_i,\al_{i'}) \log
\big(t_j^{(i)}-t_{j'}^{(i')}\big) -
\sum_{i=1}^{r}\sum_{j=1}^{k_{i}}
\sum_{b=1}^{N}
{(\La_b,\al_i)} \log
\big(t_{j}^{(i)}-x^0_b\big) .
\notag
\end{gather}
We consider $\Phi$ as a function of $t$ depending on parameters $x^0$.

Denote by $U$ the set of all points $p\in \C^k$ such that for any $\log h$ entering
\eqref{Gaudin master} (with a~nonzero coef\/f\/icient) we have $h(p)\neq 0$.
The set~$U$ is the complement in~$\C^k$ to the union of hyperplanes. The coef\/f\/icients of the logarithms in~\eqref{Gaudin master} def\/ine weights of the hyperplanes. This weighted arrangement is called
 discriminantal. Discriminantal arrangements were
considered in Section~\ref{sec Discriminantal arrangements}.

Let us construct the weight function
$\omega  : \C^k   \to \nash$ introduced in \cite{SV}, cf.~\cite{RSV}.
 Let ${b} = (b_1, \dots, b_N)$
be a sequence of nonnegative integers with  $\sum\limits_{i=1}^N b_i = k$. The set of all
such sequences will be denoted by~$B$.

For $b\in B$, let $\{c_1^1, \dots, c_{b_1}^1, \dots, c_1^N, \dots, c_{b_N}^N\}$ be a set of letters.
Let $\Sigma({b})$ be the set of all bijections $\sigma$ from the set
$\{c_1^1, \dots, c_{b_1}^1, \dots, c_1^N, \dots, c_{b_N}^N\}$ to the set of variables
$\{t_1^{(1)},\dots, t_{k_1}^{(1)}$,
\dots, $t_1^{(r)},
\dots, t_{k_r}^{(r)} \}$. Denote $d(t_j^{(i)}) = i$, and $d_{\sigma, {b}} =
(d_1^1, \dots, d_{b_1}^1, \dots, d_1^N, \dots, d_{b_N}^N)$,
where $d_l^m = d(\sigma(c_l^m))$.

To each ${b} \in B$ and $\sigma \in \Sigma({b})$ we assign the vector
\[
F_{d_{\sigma,{b}}} v = F_{d_1^1} \dots F_{d_{b_1}^1} v_1\otimes \cdots \otimes F_{d^N_1} \cdots F_{d^N_{b_N}} v_N
\in \nash .
\]
and the rational function
\[
\omega_{\sigma, {b}} = \omega^{1}_{\sigma, {b}}\big(x^0_1\big)  \cdots
\omega^N_{\sigma, {b}}\big(x^0_N\big),
\]
with
\[
\omega^e_{\sigma, {b}}\big(x^0_e\big) = \frac{1}{(\sigma(c^e_1) - \sigma(c^e_2))\cdots (\sigma(c^e_{b_e -1}) - \sigma(c^e_{b_e}))
(\sigma(c^e_{b_e}) - x^0_e)}
\]
and $\omega^e_{\sigma, {b}}(x^0_e)=1$ if $b_e=0$.
Then the weight function is given by the formula
\begin{gather}
\label{weih Func}
\omega \big(x^0,t,\bs k\big) = \sum_{{b} \in B} \sum_{\sigma \in \Sigma({b})} \omega_{\sigma,{b}}\big(x^0,t\big) F_{d_{\sigma,{b}}} v .
\end{gather}
The weight function is a function of $t$ depending on parameters $x^0$.

\begin{Lemma}[Lemma~2.1 in~\cite{MV4}]
\label{well def}
The weight function is regular on $U$.
\end{Lemma}

\subsection{Bethe vectors}
\label{sec Bethe vectors}

\begin{theorem}[\cite{B,BF,RV}]
\label{cr bethe}
Let $\g$ be a simple Lie algebra.
Let $\bs \La =(\La_1,\dots,\La_N)$, $\La_b\in\h^*,$
be a collection of weights, $\bs k=(k_1,\dots,k_r)$ a collection of nonnegative integers.
Assume that $x^0 \in \C^N$ has distinct coordinates.
Assume that  $p \in \C^k$ is a critical point of the master function
$\Phi(x^0,  \cdot ,   \bs \La, \bs k) : U\to U$. Then
the vector $\omega(x^0,p,\bs k)$ $($if nonzero$)$ belongs to $\snash$ and
 is an eigenvector of the Gaudin Hamiltonians $\Kk_1(x^0), \dots , \Kk_N(x^0)$.
\end{theorem}

The theorem also follows
directly from Theorem  6.16.2 in~\cite{SV},
cf.\ Theorem~7.2.5 in~\cite{SV},\ see also Theorem~4.2.2 in~\cite{FSV}.

The vector $\omega(x^0,p,\bs k)$ is called the Bethe vector corresponding to the critical point $p$.

\begin{thm}
\label{thm length bethe vect}
If $p$ is an isolated critical point of the master function, then
\begin{gather*}
S_{\bs\La} \big(\omega\big(x^0,p,\bs k\big), \omega\big(x^0,p, \bs k\big)\big) =
\det \left(
\frac{\partial^2\Phi}{\partial t^{(i)}_j
\partial t^{(i')}_{j'}}
\right)\big(x^0,p,\bs\La,\bs k\big) .
\end{gather*}
In particular, if the critical point  is nondegenerate, then
the Bethe vector  is nonzero.
\end{thm}

This theorem is proved for $\g=\slr$ in \cite {MV4} and for any simple Lie algebra in~\cite{V4}.

\subsection{Identif\/ication of Gaudin and naive geometric Hamiltonians}
\label{sec ident}

Let us identify constructions and  statements of Theorems~\ref{cr bethe},~\ref{thm length bethe vect} and of
 Theorem~\ref{thm F's}. First of all let us def\/ine the discriminantal arrangement associated with the Gaudin model in Theorems~\ref{thm length bethe vect}, \ref{cr bethe}.

The Gaudin model is determined by a simple Lie algebra~$\g$, a nondegenerate
$\g$-invariant bilinear form $(\,,\,)$, simple roots $\al_1,\dots, \al_r$, highest weights $\La_1,\dots,\La_N$,
distinct complex numbers $x^0_1,\dots,x^0_N$ and a vector $\bs k=(k_1,\dots,k_r)$ of nonnegative integers.

Let us take these $\h^*$, $\al_1,\dots, \al_r$,  $\La_1,\dots,\La_N$,
 $(\,,\,)$,   $\bs k=(k_1,\dots,k_r)$  as Data~\ref{ii} to def\/ine a family of prediscriminantal arrangements.
 Let us choose a point $z^0\in U(X)$ by conditions $x_b(z^0) = x^0_b$ for $b=1,\dots,N$.
The corresponding weighted discriminantal arrangement $(\A(z^0), a)$ of
 Section~\ref{sec Discriminantal arrangements} is exactly the weighted arrangement def\/ined by the master
 function in~\eqref{Gaudin master}. In particular, we have
 $\Phi(x^0,\,\cdot\,,\bs \La,\bs k)=\Phi (z^0, \,\cdot\,)$ where $\Phi(z^0,\,\cdot\,)$ is the master function of the arrangement
 $(\A(z^0),a)$.

 In \cite{SV} an isomorphism $\gamma : W^-(z^0) \to \nash$ was constructed with the following properties $(i)$--$(iv)$.
 \begin{enumerate}\itemsep=0pt
 \item[$(i)$]
 The isomorphism $\gamma$
identif\/ies $\sing W^-(z^0)$ with $\snash$.

\item[$(ii)$]
The isomorphism $\gamma$
identif\/ies a  naive geometric Hamiltonian $K_{\der_{x_b}}^*|_{\sing W^-(z^0)}$
with the Gaudin Hamiltonian $ \Kk_b(x^0)|_\snash$ up to addition of a scalar operator.
More precisely, for any $b=1,\dots,N$, there is a number $c_b$ such that
\[
\gamma|_{\sing W^-(z^0)} \big(K_{\der_{x_b}}^*+c_b\big)|_{\sing W^-(z^0)}
=
\Kk_b\big(x^0\big)|_\snash
\gamma|_{\sing W^-(z^0)} .
\]

\item[$(iii)$]\sloppy
The isomorphism $\gamma$
identif\/ies the contravariant form on $W^-(z^0)$ (multiplied by
$(-1)^k k_1!\cdots k_r!$) with the tensor Shapovalov form on $\nash$,
\[
 S_{\bs\La}(\gamma(x),\gamma(y))  =  (-1)^k k_1!\cdots k_r!  S^{(a)}(x,y)
 \qquad
\text{for any}\
x,y\in \nash .
\]

\item [$(iv)$]
Let $\bs E \in \C(t)_{U(\A(z^0))}\otimes W(z^0)$ be the canonical element of the arrangement $\A(z^0)$.
Let ${\rm Ant} = \sum_{\sigma\in \Sk} (-1)^\sigma \sigma : W(z^0)\to  W^-(z^0)$ be the anti-symmetrization
operator. Consider the element $(1\otimes {\rm Ant})\bs E \in
\C(t)_{U(\A(z^0))}\otimes W^-(z^0)$.
Let $\omega(x^0,\cdot,\bs k)$ $ \in \C(t)_{U(\A(z^0))} \otimes \nash$ be the weight function
def\/ined by \eqref{weih Func}.
Then
\[
 (1\otimes \gamma)(1\otimes {\rm Ant}) \bs E  = k_1!\cdots k_r!  \omega .
 \]
  \end{enumerate}
 See Theorems 5.13, 6.16.2, 6.6 and 7.2.5 in~\cite{SV}, see also Section~5.5 in~\cite{V4}.

Now Theorem~\ref{thm length bethe vect} follows from items $(iv)$ and $(vi)$ of Theorem~\ref{thm F's}
and Theorem~\ref{cr bethe} follows from Lemma~\ref{lenth -} and
Theorem~\ref{first theorem}.

From statements $(i)$--$(iv)$ and Theorem \ref{thm F's} we get the following improvement of Theo\-rem~\ref{thm length bethe vect}.

\begin{thm}
For any simple Lie algebra $\g$, if $p$ is an isolated critical point of the master function $\Phi(x^0,\cdot,\bs\La,\bs k)$, then
the corresponding Bethe vector $\omega(x^0,p,\bs k)$ is nonzero.
\end{thm}

\subsection{Bethe algebra}
\label{sec Bethe algebra}

The subalgebra of
$\End(\Sing \V[\mu])$ generated by the Gaudin Hamiltonians can be extended to a~larger commutative
subalgebra called the Bethe algebra. A construction of the Bethe algebra for any simple Lie algebra~$\g$
is given in~\cite{FFR}.
As a result of that construction,
for any~$x^0$  one obtains a commutative subalgebra $\B(x^0)\subset (U\g)^{\otimes N}$ which
commutes with the diagonal subalgebra $U\g\subset (U\g)^{\otimes N}$.
To def\/ine the Bethe algebra of $\V$ or of $\Sing \V[\mu]$ one considers the
 image of~$\B(x^0)$ in~$\End (\V)$ or in~$\End(\Sing \V[\mu])$. The Gaudin Hamiltonians $\Kk_b(x^0)$
 are elements of the Bethe algebra of~$\V$ or of~$\snash$.

 A more straightforward construction of the Bethe algebra is known for the Gaudin model of~$\frak {gl}_{r+1}$, see \cite{T}.
 Below we give its description.

Let $e_{ij}$, $i,j=1\lc r+1$, be the standard generators of
$\glr$ satisfying the relations
$[e_{ij},e_{sk}]=\dl_{js}e_{ik}-\dl_{ik}e_{sj}$. Let $\h\subset\glr$ be
the Cartan subalgebra generated by $e_{ii}$,  $i=1\lc r+1$.
Let $\h^*$ be the dual space. Let
$\epsilon_i$, $i=1\lc r+1,$ be the basis of $\h^*$ dual to the basis
$e_{ii}$, $i=1\lc r+1$, of $\h$.
Let $\al_1,\dots,\al_{r} \in \h^*$ be simple roots, $\al_i
=\epsilon_i-\epsilon_{i+1}$. Let
$(\,,\,)$ be the standard scalar product on $\h^*$ such that
the basis $\epsilon_i$, $i=1\lc r+1$,
is orthonormal.

Let $\glrs=\glr\otimes\C[s]$ be the Lie algebra of $\glr$-valued
polynomials with the pointwise commutator. For $g\in\glr$, we set
$g(u)=\sum\limits_{i=0}^\infty (g\otimes s^i)u^{-i-1}$.
We identify $\glr$ with the subalgebra $\glr\otimes1$ of constant polynomials
in $\glrs$. Hence, any $\glrs$-module has a~canonical structure of
a $\glr$-module.

For each $a\in\C$, there exists an automorphism $\rho_a$ of $\glrs$,
 $\rho_a:g(u)\mapsto g(u-a)$. Given a~$\glrs$-module $W$, we denote by $W(a)$
the pull-back of $W$ through the automorphism $\rho_a$. As $\glr$-modules,
$W$ and $W(a)$ are isomorphic by the identity map.

We have the evaluation homomorphism,
${\glrs\to\glr}$,  ${g(u) \mapsto g u^{-1}}$.
Its restriction to the subalgebra $\glr\subset\glrs$ is the identity map.
For any $\glr$-module $W$, we denote by the same letter the $\glrs$-module,
obtained by pulling $W$ back through the evaluation homomorphism.

Given an algebra $A$ and an ${(r+1)\times (r+1)}$-matrix $C=(c_{ij})$ with entries in $A$,
we def\/ine its {\it row determinant} to be
\[
\rdet\, C =
\sum_{\si\in \Sigma_{r+1}} (-1)^\si c_{1\si(1)}c_{2\si(2)}\cdots c_{r+1\,\si(r+1)} .
\]

Def\/ine the universal dif\/ferential operator $\D_\B$
by the formula
\[
\D_\B= \rdet\left( \begin{matrix}
\der_u-e_{11}(u) & - e_{21}(u)& \dots & -e_{r+1\,1}(u)\\[3pt]
-e_{12}(u) &\der_u-e_{22}(u)& \dots & -e_{r+1\,2}(u)\\[1pt]
\dots & \dots &\dots &\dots \\[1pt]
-e_{1\,r+1}(u) & -e_{2\,r+1}(u)& \dots & \der_u-e_{r+1\,r+1}(u)
\end{matrix}\right).
\]
We have
\[
\D_\B= \der_u^{r+1}+\sum_{i=1}^{r+1} B_i \der_u^{{r+1}-i},
\qquad
B_i = \sum_{j=i}^\infty B_{ij} u^{-j},
\qquad
B_{ij}\in U\glrs  .
\]
The unital subalgebra of $U\glrs$ generated by $B_{ij}$,  $i=1\lc {r+1}$,
 $j\geq i$, is called the {\it Bethe algebra\/} and denoted by $\B$.

By \cite{T}, cf.~\cite{MTV1}, the algebra $\B$ is commutative,
and $\B$ commutes with the subalgebra $U\glr\subset U\glrs$.

As a subalgebra of $U\glrs$, the algebra $\B$ acts on any $\glrs$-module
$W$. Since $\B$ commutes with $U\glr$, it preserves the $\glr$ weight subspaces
of $W$ and the subspace
$\Sing W$ of $\glr$-singular vectors.

If $W$ is a $\B$-module, then the image of $\B$ in $\End(W)$ is called
the  Bethe algebra of $W$.

Let $\V=\otimes_{b=1}^n V_{\La_b}$ be a tensor product of irreducible highest weight $\glr$-modules.
For given $x^0=(x^0_1,\dots,x^0_N)$, consider $\V$ as the $\glrs$-module $\otimes_{b=1}^n
V_{\La_b}(x_b^0)$. This $\glrs$-module structure on $\V$ provides $\V$ with  a Bethe algebra,
a commutative subalgebra of $\End(\V)$. It is known that this Bethe algebra of~$\V$ contains the Gaudin Hamiltonians
$\Kk_b(x^0)$, $b=1,\dots,N$. In fact, the Gaudin Hamiltonians are suitably normalized residues of the
generating function~$B_2(u)$, see Appendix~B in~\cite{MTV1}.

\begin{thm}[\cite{MTV7}]\sloppy
\label{thm B has symm opers}
Consider $\V$ as the $\glrs$-module $\otimes_{b=1}^n
V_{\La_b}(x^0_b)$. Then any element \mbox{$B\in \B$} acts on $\V$ as a
symmetric operator with respect to the tensor Shapovalov form,
$S_{\bs\La}(Bu,v)= S_{\bs\La}(u,Bv)$ for any $u,v\in \V$.
\end{thm}

\subsection[$\glr$ Bethe algebra and critical points of the master function]{$\boldsymbol{\glr}$ Bethe algebra and critical points of the master function}
\label{sec crit pointS and glr bethe alg}

In \cite{MTV5} the following generalization of Theorem \ref{cr bethe} for
$\g=\glr$ was obtained.

Let $\g=\glr$.
A sequence of integers $\La=(\la_1,\dots,\la_{r+1})$ such that
$\la_1\ge\la_2\ge\cdots\ge\la_{r+1}\ge0$ is called a  partition with at most
${r+1}$ parts.
We identify partitions $\La$
with vectors $\la_1\epsilon_1+
\cdots + \la_{r+1}\epsilon_{r+1}$ of $\h^*$.

Let $\bs \La = (\La_1,{\dots}, \La_N)$ be a collection of partitions,
where $\La_b=(\la_{b,1},{\dots},\la_{b,r+1})$ and
\mbox{$\la_{b,r+1}\!=0$}.
Let $\bs k = (k_1,\dots,k_r)$ be nonnegative integers such that
\[
\La_\infty  =  \sum_{b=1}^N \La_b - \sum_{i=1}^{r} k_i\al_i
\]
is a partition.
We  consider the $\glr$ Gaudin model with parameters $x^0=(x^0_1,\dots,x_N^0)$
on $\snash$ where $\V = V_{\La_1}\otimes \dots\otimes V_{\La_N}$.

Consider the master function $\Phi(x^0\!,t,\bs\La,\bs k)$ def\/ined by~\eqref{Gaudin master}
and the weight function $\omega(x^0\!,t,\bs k)\!$ def\/ined by~\eqref{weih Func}.

Let $u$ be a variable.
Def\/ine polynomials $T_1,\dots, T_r\in \C[u]$,
$Q_1,\dots,Q_r\in \C[u,t]$,
\[
T_i(u)  =  \prod_{b=1}^N\,(u-x_b)^{(\La_b,\al_i)}  ,
\qquad
Q_i(u,  t)  =  \prod_{j=1}^{k_i}\big (u- t^{(i)}_j\big)  ,
\]
and the dif\/ferential operator
\begin{gather*}
\D_{\Phi} =
\left(\der_u - \log' \left( \frac { T_1\dots T_r } { Q_{1} }\right ) \right)
\\
\phantom{\D_{\Phi} = }{}
 \times\left( \der_u - \log' \left( \frac {Q_{1}T_2\cdots T_{r} } {Q_{2} }\right ) \right)
\cdots \left( \der_u - \log' \left( \frac {Q_{r-1} T_r}{ Q_r } \right) \right)
( \der_u - \log' ( Q_r ) ) ,
\end{gather*}
where $\der_u = d/du$ and $\log' f$ denotes $ (df/du)/f$. We have
\[
\D_{\Phi}
  =
\der_u^{N+1}+\sum_{i=1}^{N+1} G_i \der_u^{N+1-i},
\qquad
G_i = \sum_{j=i}^\infty G_{ij} u^{-j} ,
\]
where $G_{ij} \in \C[\bs t]$.

Let $p\in U$ be an isolated critical point of the master function
$\Phi(x,,\cdot\,,\bs\La,\bs k)$ with Milnor number~$\mu$. Let $A_{p,\Phi}$
be its local algebra. For  $f\in \C( t)_U$ denote by $[f]$ the
image of $f$ in $\AT$.
Denote
\[
[\D_{\Phi}]
  =
\der_u^{N+1}+\sum_{i=1}^{N+1}  [G_i] \der_u^{N+1-i},
\]
where
$[G_i] = \sum\limits_{j=i}^\infty [G_{ij}] u^{-j}$.
 Let
$
[\omega] \in \AT\otimes \nash
$
be the element induced by the weight function.
The element $[\omega]$ belongs to
$\AT\otimes \snash$, see \cite{SV}, and we
 have $S_{\bs \La}([\omega],[\omega]) =
[{\rm Hess}_t\Phi]$ , see \cite{MV4, V4}.

\begin{thm} [\cite {MTV1}]
\label{thm MtV}
For any $i=1,\dots,r+1$, $j\geq i$, we have
\[
(1\otimes B_{ij})  [\omega]
 =
([G_{ij}]_p \otimes 1)  [\omega]
\]
in   $\AT\otimes \snash$.
\end{thm}

This statement is the Bethe ansatz method to construct eigenvectors of
the Bethe algebra in the $\glr$ Gaudin model starting with a critical
point of the master function.

Let $g_1,\dots, g_{\mu}$ be a basis of $\AT$
considered as a $\C$-vector space. Write
$[\omega]_p =\sum_i g_i\otimes w_i$, with $ w_i \in \snash$.
Denote by $\mc Y_p \subset {\rm Sing}\,\nash$ the vector subspace spanned by
$w_1,\dots,w_{\mu}$. Let $(\,,\,)_p$ be the
  bilinear form on $A_{p,\Phi}$.
Def\/ine a linear map
\[
\al_p :\ A_{p, \Phi}  \to  \mc Y_p  ,
\qquad
f \mapsto (f,[\omega]_p)_p = \sum_{i=1}^{\mu}  (f,g_i)_p w_i .
\]

\begin{thm} [\cite{MTV5}]
\label{thm main Bethe}
\quad
\begin{enumerate}\itemsep=0pt
\item[$(i)$]
The subspace $\mc Y_p\subset \snash$ is a $\B$-submodule.
Let $A_{\mc Y_p} \subset \End\,(\mc Y_p)$ be the Bethe algebra of $\mc Y_p$.
Denote by $\bar B_{ij}$ the image in $A_{\mc Y_p}$ of generators $B_{ij}\in \B$.

\item[$(ii)$]
The map $\al_p : A_{p, \Phi} \to \mc Y_p$ is an isomorphism of vector spaces.
\item[$(iii)$]
The map $[G_{ij}]_p \mapsto \bar B_{ij}$ extends uniquely to an algebra
isomorphism $\beta_p: \AT \to A_{\mc Y_p}$.
\item[$(iv)$]
The isomorphisms $\al_p$ and $\beta_p$ identify the regular representation
of $\AT$ and the $\B$-mo\-du\-le~$\mc Y_p$,  that is, for any $f,g\in\AT$ we have
$\al_p(fg)=\beta_p(f)\al_p(g)$.
\item[$(v)$]
The value  of the weight function at $p$ is a nonzero vector of~$\snash$.
\item[$(vi)$]
Let $p_1,\dots,p_d$ be a list of
isolated critical points  of $\Phi(x^0,\,\cdot\,,\bs \La,\bs k)$
 such that
the orbits $O(p_1),\dots,O(p_d)$ do not intersect.
Let $\mc Y_{p_s}=\al_{p_s}(A_{p_s,\Phi})\subset \snash$, $s=1,\dots,d$, be the corresponding subspaces.
Then the sum of these subspaces is direct and orthogonal.

\item[$(vii)$]
Denote $\mc Y=\oplus_{s=1}^d\mc Y_{p_s}$. Denote $A_{\mc Y}$ the Bethe algebra of $\mc Y$.
Consider the isomorphisms
\begin{gather*}
\al = \oplus_{s=1}^d\al_{p_s}  :\
\oplus_{s=1}^d
A_{p_s, \Phi}
\to  \oplus_{s=1}^d \mc Y_{p_s}  ,
\\
\beta =  \oplus_{s=1}^d\beta_{p_s}   :\ \oplus_{s=1}^dA_{p_s, \Phi}  \to  \oplus_{s=1}^d A_{Y_{p_s}} .
\end{gather*}
Then
\begin{enumerate}\itemsep=0pt
\item[$(a)$]
$A_{\mc Y}=\oplus_{s=1}^d A_{\mc Y_{p_s}}$;
\item[$(b)$]
$A_{\mc Y}$ is a maximal commutative subalgebra of $\End(\mc Y)$;
\item[$(c)$] the isomorphisms
$\al$, $\beta$
identify the regular representation
of the algebra
 $\oplus_{s=1}^d A_{p_s, \Phi}$ and the $A_{\mc Y}$-module $\mc Y$.
\end{enumerate}
\end{enumerate}
\end{thm}

Let us identify constructions and  statements of Theorem~\ref{thm main Bethe} and of Theorems
 \ref{thm F's}, \ref{thm BAD S}.

 Consider the discriminantal arrangement $\A(z^0)$, $z^0\in U(X)$, def\/ined in Section~\ref{sec ident}
 and the isomorphism $\gamma|_{\sing W^-(z^0)} : \sing W^-(z^0) \to \snash$ of that section.
 According to statements $(i)$--$(iv)$ of Section \ref{sec ident}, we have
  $\gamma (Y^-_{O(p_s)}) = \mc Y_{p_s}$, where $Y^-_{O(p_s)}\subset \sing W^-(z^0)$ is the subspace in Theorems~\ref{thm F's},~\ref{thm BAD S} and $\mc Y_{p_s}\subset \snash$ is the subspace in Theorem~\ref{thm main Bethe}. Moreover,
 $\gamma$ identif\/ies the algebra $A_{Y^-_{O(p_s)}}$ of geometric Hamiltonians on
 $Y^-_{O(p_s)}$ with the Bethe algebra of~$\mc Y_{p_s}$.

Additional information that is given by Theorems~\ref {thm F's},~\ref{thm BAD S} is the following theorem.

\begin{thm}
\label{lem shap=res}
The monomorphism
$\al_p$ of Theorem~{\rm \ref{thm main Bethe}} identifies the tensor Shapovalov form on $\mc Y_p$ and the residue bilinear form
on $A_{p,\Phi}$.
\end{thm}

\begin{proof}
The theorem is a corollary of properties $(iii)$, $(iv)$ in Section~\ref{sec ident} and Lemma~\ref{lem gamma}.
\end{proof}

\subsection{Expectations}
\label{sec expectations}
One may expect that for the Gaudin models of any simple Lie algebras the associated Bethe algebra def\/ined in \cite{FFR}
coincides with the algebra of geometric Hamiltonians associated with the arrangement of the corresponding master function.
That topic will be discussed in a forthcoming paper.

Below we consider three examples in which we have a well-def\/ined algebra of geometric Hamiltonians on $\snash$.

\subsubsection{Example}
\label{example 2}

Consider the Gaudin model corresponding to the following data: the Lie algebra $\glr$,
the collection of dominant integral weights $\bs \La=(\La_1,\dots,\La_N)$ with $\La_b=(1,0,\dots,0)$ for all~$b$,
a~vector of nonnegative integers
$\bs k=(k_1,\dots,k_r)$ such that $\La_\infty$ is a partition,
 a collection of generic distinct complex numbers $x^0=(x_1^0,\dots,x^0_N)$.

 By~\cite{MV4}, for generic $x^0$ the associated master function has a collection of critical points $p_1,\dots,p_d$ such that
 the $\Sk$-orbits of these points do not intersect and the sum of Milnor numbers of these points equals the dimension of
 $\snash$. In this case Theorem \ref{thm BAD S} def\/ines a~ma\-xi\-mal commutative subalgebra $A_{\sing (W^-(z^0))}$
$ \subset $ $\End(\sing (W^-(z^0)))$ containing
 naive geometric  Hamiltonians. The isomorphism $\gamma$ sends
  $A_{\sing (W^-(z^0))}$ to a maximal commutative subalgebra
 of $\End(\snash)$ containing the Gaudin Hamiltonians. By \cite{MTV7,MTV8}, the Bethe algebra of $\snash$ is a maximal commutative subalgebra
 of $\End(\snash)$ and the Bethe algebra of $\snash $ is generated by the Gaudin Hamiltonians. Hence, in this case,
  $\gamma$ establishes an isomorphism
 of the algebra of geometric Hamiltonians $A_{\sing (W^-(z^0))}$ and the Bethe algebra of $\snash$.

\subsubsection{Example}
\label{example 3}

Consider the Gaudin model corresponding to the following data: the Lie algebra $\frak{gl}_2$,
a collection of weights $\bs \La=(\La_1,\dots,\La_N)$ with $\La_b= \la_b\al_1$ such that $\la_b\in\R_{<0}$ for all $b$,
a  nonnegative integer
$\bs k=(k_1)$,  a collection of distinct complex numbers $x^0=(x_1^0,\dots,x^0_N)$.

Let $(\A(z^0),a)$ be the associated discriminantal arrangement def\/ined in Section~\ref{sec ident}.
By our assumptions, the weights $a$ of the discriminant arrangement are all positive.
Corollary \ref{cor unbal equiv} def\/ines in this case a maximal commutative subalgebra of $\End (\sing W^-(z^0)))$ and the isomorphism
$\gamma$ sends this subalgebra to a maximal commutative subalgebra of $\End(\snash)$, which contains  Gaudin Hamiltonians.
One may expect that this commutative subalgebra of $\End(\snash)$ coincides with the Bethe algebra of $\snash$.

\subsubsection{Example}
\label{example 1}

Consider the Gaudin model corresponding to the following data: a simple Lie algebra $\g$,
a~collection of dominant integral weights $\bs \La=(\La_1,\dots,\La_N)$, a vector of nonnegative integers
$\bs k=(k_1,\dots,k_r)$ with $k_i\leq 1$ for all $i$, a collection of distinct complex numbers $x^0=(x_1^0,\dots,x^0_N)$.

{\sloppy
Let $(\A(z^0),a)$ be the associated discriminantal arrangement.
By our assumptions, the weights~$a$ of the discriminant arrangement are all negative and the group $\Sk$ is trivial.
In this case the isomorphism $\gamma$ identif\/ies $\sing \FF^k(\A(z^0))$ and the space $\snash$.
Theorem~\ref{thm ham of unb arr} def\/ines in this case a maximal commutative subalgebra of $\End (\sing \FF^k(\A(z^0)))$ and the isomorphism
$\gamma$ sends this subalgebra to a maximal commutative subalgebra of $\End(\snash)$, which contains the Gaudin Hamiltonians.
One may expect that this commutative subalgebra of $\End(\snash)$ coincides with the Bethe algebra of~$\snash$.

}

\subsection*{Acknowledgments}
\label{sec remarks}

The idea  that an analog of the Bethe ansatz construction
does exist for an arbitrary arrangement of hyperplanes
was formulated long time ago in \cite{V2}.
That program had been realized partially in \cite{V4}.

This paper is an extended exposition of my lectures at
Mathematical Society of Japan  Seasonal Institute   on Arrangements of Hyperplanes in
August of 2009. I thank  organizers  for invitation and  Hokkaido University for hospitality.
I thank for hospitality Universit\'{e} Paul Sabatier in Toulouse, where this paper had been f\/inished.
I thank E.~Mukhin, V.~Schechtman, V.~Tarasov, H.~Terao for discussions.

The author  was supported in part by NSF grant DMS-0555327.

\LastPageEnding

\end{document}